\RequirePackage{ifpdf}
\ifpdf 
\documentclass[pdftex]{sigma}
\else
\documentclass{sigma}
\fi

\numberwithin{equation}{section}
\numberwithin{lemma}{section}
\numberwithin{theorem}{section}
\numberwithin{definition}{section}
\numberwithin{proposition}{section}

\begin{document}

\allowdisplaybreaks

\renewcommand{\PaperNumber}{019}

\FirstPageHeading

\ShortArticleName{The Decomposition of Global Conformal Invariants: Some Technical Proofs.~I}

\ArticleName{The Decomposition of Global Conformal Invariants:\\ Some Technical Proofs.~I}

\Author{Spyros ALEXAKIS}

\AuthorNameForHeading{S. Alexakis}

\Address{Department of Mathematics, University of Toronto,  40 St. George Street, Toronto,   Canada}
\Email{\href{mailto:alexakis@math.toronto.edu}{alexakis@math.toronto.edu}}

\ArticleDates{Received April 01, 2010, in f\/inal form February 15, 2011;  Published online February 26, 2011}

\Abstract{This paper forms part of a larger work
 where we  prove a conjecture of Deser and Schwimmer
regarding the algebraic structure of ``global conformal
invariants''; these  are def\/ined to
be conformally invariant integrals of geometric scalars.
 The conjecture asserts that the integrand of
any such  integral can be expressed as a linear
combination of a local conformal invariant, a divergence and of
the Chern--Gauss--Bonnet integrand.}

\Keywords{conormal geometry; renormalized volume; global invariants;
Deser--Schwimmer conjecture}

\Classification{53B20; 53A55}

\section{Introduction}

The present paper complements the monograph~\cite{alexakis} and also~\cite{a:dgciI,a:dgciII,alexakis4,alexakis5}
in proving a conjecture of Deser and Schwimmer~\cite{ds:gccaad} regarding
the algebraic structure of ``global conformal invariants''\footnote{A
formulation of the conjecture which is closer to the mathematical terminology
 used here can also be found in~\cite{hening-sken}.}.
 Here we provide a
proof of certain lemmas which were announced in the second  chapter of~\cite{alexakis}.
The proofs of these claims  do not logically depend on material
appearing elsewhere in this entire work.

For the reader's convenience, we brief\/ly review the Deser--Schwimmer
conjecture.

 We recall that a global
conformal invariant is an integral of a natural scalar-valued function of
Riemannian metrics, $\int_{M^n}P(g)dV_g$, which
remains invariant under conformal re-scalings of the underlying
metric\footnote{See the introduction of \cite{alexakis}
for a detailed discussion on scalar Riemannian invariants and
on applications of the Deser--Schwimmer conjecture.}.
More precisely, $P(g)$ is assumed to be a linear combination, $P(g)=\sum_{l\in L} a_l C^l(g)$,
where each $C^l(g)$ is a complete contraction in the form:
\begin{gather}
\label{contraction}
{\rm contr}^l\big(\nabla^{(m_1)}R\otimes\dots\otimes\nabla^{(m_s)}R\big).
\end{gather}
Here each factor $\nabla^{(m)}R$ stands for the $m^{th}$ iterated
covariant derivative of the curvature tensor~$R$. $\nabla$ is the Levi-Civita
connection of the metric $g$ and $R$ is the curvature associated to this connection.
The contractions are taken with respect to the quadratic form $g^{ij}$.
In the present paper, along with~\cite{alexakis,alexakis5} we prove:

\begin{theorem}
\label{thetheorem}
Assume that $P(g)=\sum_{l\in L} a_l C^l(g)$, where each $C^l(g)$ is a
complete contraction in the form \eqref{contraction}, with weight $-n$.
Assume that for every closed Riemannian manifold $(M^n,g)$ and every $\phi\in C^\infty (M^n)$:
\[
\int_{M^n}P\big(e^{2\phi}g\big)dV_{e^{2\phi}g}=\int_{M^n}P(g)dV_g.
\]

We claim that $P(g)$ can then be expressed in the form:
\begin{gather}
\label{thedecomposition}
P(g)=W(g)+{\rm div}_iT^i(g)+ ({\rm const})\cdot\operatorname{Pfaf\/f}(R_{ijkl}).
\end{gather}
Here $W(g)$ stands for a local conformal invariant of weight $-n$  $($meaning
that $W(e^{2\phi}g)=e^{-n\phi}W(g)$ for every $\phi\in C^\infty (M^n))$,
${\rm div}_iT^i(g)$ is the divergence of a Riemannian vector field of
weight $-n+1$,   $\operatorname{Pfaf\/f}(R_{ijkl})$ is the Pfaffian of the curvature tensor, and $({\rm const})$ is a~constant.
\end{theorem}

 {\bf Remark.} The notion of a ``global conformal invariant'' is closely related to  the algebraic structure  of
  {\it conformal anomalies}\footnote{These are sometimes called Weyl anomalies.},
  which were the motivation of Deser and Schwimmer~\cite{ds:gccaad}
  (see also  especially \cite{hening-sken}); however the two notions are not identical: very broadly speaking,
  a~conformal anomaly represents the failure of an action that depends on a metric $g$
   to be conformally invariant. A well-studied such anomaly (which directly gives rise to the
   global conformal inva\-riants we study here)
   is considered in~\cite{hening-sken} in the context of the AdS-CFT
   correspondence; it is the renormalized
   volume $\cal{A}$  of a Poincar\'e--Einstein metric $(X^{n+1}, h)$
    with  conformal inf\/inity $(M^n,[g])$, when  $n$ is even. In
    that setting if $\dot{g}=2\phi g$, then $\dot{\cal{A}}(g)=\int_{M^n}P(g)\phi dV_g$.
    Furthermore as noted by Henningson--Skenderis \cite{hening-sken}, the integral $\int_{M^n}P(g)dV_g$
    is a global conformal invariant in the sense described above.
(This fact was noted in \cite{hening-sken} on general grounds;
Graham in \cite{graham:volren} showed the conformal invariance of
$\int_{M^n}P(g)dV_g$ by explicitly studying the conformal variation of the renormalized volume.)
    Thus, the decomposition (\ref{thedecomposition})
    that we derive for $P(g)$ corresponds to a decomposition of the integrand in
    the conformal anomaly; this conforms with the stipulation
    of \cite{hening-sken}\footnote{See also interesting work of
    Boulanger in \cite{boulanger} on Weyl anomalies
    which satisfy the Wess--Zumino consistency conditions.}.

We now digress in order to discuss the relation of this entire work (consisting
of the present paper and \cite{alexakis, alexakis4,alexakis5})
with classical and recent work on  local invariants in various geometries.

{\bf Local invariants and Fef\/ferman's program for the Bergman
and Szeg\"o kernels.} The theory of {\it local}
invariants of Riemannian structures
(and indeed, of more general geometries,
e.g.~conformal, projective, or CR)  has a long history.
The original foundations of this
f\/ield were laid in the work of Hermann Weyl
and \'Elie Cartan, see \cite{w:cg, cartan}.
The task of writing out local invariants of
a given geometry is intimately connected
with understanding which polynomials in
a space of tensors with  given symmetries
 remain invariant under the action of a Lie group.
In particular, the problem of writing down all
 local Riemannian invariants reduces to understanding
the invariants of the orthogonal group.

 In more recent times, a major program was initiated by C.~Fef\/ferman in \cite{f:ma}
aimed at f\/inding all scalar local invariants in CR geometry. This was motivated
by the problem of understanding the
local invariants which appear in the asymptotic expansions of the
Bergman and Szeg\"o kernels of strictly pseudo-convex CR manifolds,
 in a similar way to which Riemannian
invariants appear in the asymptotic expansion
of the heat kernel; the study of the local invariants
in the singularities of these kernels led to important breakthroughs
in~\cite{beg:itccg} and more recently by Hirachi in~\cite{hirachi1}.
 This program was later extended  to conformal geometry in~\cite{fg:ci}.
Both these geometries belong to a
broader class of structures, the
{\it parabolic geometries}; these are
structures which admit a principal bundle with
structure group a parabolic subgroup~$P$ of a semi-simple
Lie group~$G$, and a Cartan connection on that principle bundle
(see the introduction in~\cite{cg1}).
An important question in the study of these structures
is the problem of constructing all their local invariants, which
can be thought of as the {\it natural, intrinsic} scalars of these structures.

  In the context of conformal geometry, the f\/irst (modern) landmark
in understanding {\it local conformal invariants} was the work of Fef\/ferman
and Graham in 1985~\cite{fg:ci},
where they introduced their {\it ambient metric}. This allows one to
construct local conformal invariants of any order in odd
dimensions, and up to order~$\frac{n}{2}$ in even dimensions.
A natural question is then whether
{\it all} invariants arise via this construction.

The subsequent work of Bailey--Eastwood--Graham \cite{beg:itccg} proved that
this is indeed true in odd dimensions; in even dimensions,
they proved that the result holds
when the weight (in absolute value) is bounded
by the dimension. The ambient metric construction
in even dimensions was recently extended by Graham--Hirachi, \cite{grhir}; this enables them to
identify in a satisfactory manner {\it all} local conformal invariants,
even when the weight (in absolute value) exceeds the dimension.

 An alternative
construction of local conformal invariants can be obtained via the {\it tractor calculus}
introduced by Bailey--Eastwood--Gover in~\cite{bego}. This construction bears a strong
resemblance to the Cartan conformal connection, and to
the work of T.Y.~Thomas in 1934~\cite{thomas}. The tractor
calculus has proven to be very universal;
tractor bundles have been constructed \cite{cg1} for an entire class of parabolic geometries.
The relation between the conformal tractor calculus and the Fef\/ferman--Graham
ambient metric  has been elucidated in~\cite{cg2}.

The present work, while pertaining to the question above
(given that it ultimately deals with the algebraic form of local
{\it Riemannian} and {\it conformal} invariants),
nonetheless addresses a~dif\/ferent
{\it type} of problem. We here consider Riemannian invariants $P(g)$ for
which the {\it integral} $\int_{M^n}P(g)dV_g$ remains invariant
under conformal changes of the underlying metric. We then seek to understand
the algebraic form of the {\it integrand} $P(g)$,
ultimately proving that it can be de-composed
in the way that Deser and Schwimmer asserted.
It is thus not surprising that the prior work on
 the construction and understanding of local {\it conformal}
invariants plays a central role in this
endeavor, both in~\cite{alexakis} and in the present paper.

{\bf Index theory.} Questions similar to the Deser--Schwimmer
conjecture  arise naturally
in index theory; a good reference for such questions is~\cite{bgv}.
 For example, in the heat kernel proof of the index theorem (for Dirac operators)
by Atiyah--Bott--Patodi~\cite{abp},
the authors were led to consider integrals arising in the (integrated) expansion of
the heat kernel over Riemannian manifolds of general Dirac operators,
and sought to understand the
local structure  of the integrand\footnote{We note that
the geometric setting in \cite{abp} is more
general than the one in the Deser--Schwimmer conjecture:
 In particular one considers vector bundles, equipped with an auxiliary connection,
  over compact Riemannian manifolds; the
local invariants thus depend {\it both} on the
curvature of the Riemannian metric {\it and} the curvature of the connection.}.
In that setting, however, the fact that one deals with
a {\it specific}  integrand which arises in the heat kernel expansion
plays a key role in the understanding of its local  structure.
This is true both of the original proof of Patodi, Atiyah--Bott--Patodi \cite{abp}
and of their subsequent simplif\/ications and generalizations by Getzler,
Berline--Getzler--Vergne, see \cite{bgv}.

The closest analogous problem to the one considered here
is the work of Gilkey and Branson--Gilkey--Pohjanpelto~\cite{gilkey,bgp}.
In~\cite{gilkey}, Gilkey considered Riemannian invariants $P(g)$ for
which the {\it integral} $\int_{M^n}P(g)dV_g$ on any given
(topological) manifold $M^n$ has a given
value, {\it independent of the metric $g$}.
He proved that $P(g)$ must then be equal to a divergence,
plus possibly a~multiple of the Chern--Gauss--Bonnet integrand,
if the weight of $P(g)$ agrees with the dimension in absolute value.
In \cite{bgp} the authors considered the problem of Deser--Schwimmer
for locally conformaly f\/lat metrics and derived the same decomposition
(for {\it locally conformaly flat metrics}) as in~\cite{gilkey}.
Although these two results can be considered  precursors of ours, the methods
there are entirely dif\/ferent from the ones here; it is highly unclear whether
the methods of \cite{gilkey,bgp} could be applied to the problem at hand.

We next review the postponed claims that we will be proving here.

In Section~\ref{sequals2} we recall from~\cite{alexakis} that~$P(g)$ is thought of as a linear combination of complete
contractions involving factors of two types:
 iterated covariant derivatives of the Weyl tensor and iterated covariant
 derivatives of the Schouten tensor\footnote{We
refer the reader to the introduction of \cite{alexakis}
for a def\/inition of these classical tensors.
The Weyl tensor is the trace-free part of the curvature
tensor and is conformally invariant, while the Schouten tensor is
a trace-adjustment of the Ricci tensor
(and is not conformally invariant).}. In other words,
 we write $P(g)=\sum_{l\in L} a_l C^l(g)$, where each $C^l(g)$
is in the form~(\ref{drassi}) below. (Recall from the introduction in \cite{alexakis}
that all complete contractions in $P(g)$ are assumed to
have weight $-n$.)

In Section~\ref{sequals2} we prove certain special cases of
Propositions~2.3.1, 2.3.2  from Chapter~2 in \cite{alexakis} that were postponed to the
present paper; these are contained in the next lemma:

\begin{lemma}
\label{theatro}
Let $P(g)$ be a Riemannian invariant of weight ${-}n$
such that the integral $\int_{\!M^n}\!\! P(g)dV_g\!$ is a global
conformal invariant. Assume that $P(g)$ is in the form
$P(g)=\sum_{l\in L} a_l C^l(g)$, where each $C^l(g)$ is a complete contraction:
\begin{gather}
\label{drassi}
{\rm contr}\big(\nabla^{(m)}W\otimes\dots\otimes
\nabla^{(m')}W\otimes\nabla^{(a)}P\otimes\dots\otimes\nabla^{(a')}P\big).
\end{gather}
Let $\sigma>0$ be the minimum number of factors among the
complete contractions $C^l(g)$, $l\in L$. Denote by $L_\sigma\subset L$
the corresponding index set. We then claim that if
$\sigma\le 2$, there exists a  local conformal invariant~$W(g)$ and
 a Riemannian vector field $T^i(g)$ such that:
\begin{gather*}
P(g)={\rm div}_iT^i(g)+W(g)+\sum_{l\in L'} a_l C^l(g),
\end{gather*}
where each $C^l(g)$, $l\in L'$ in the r.h.s.\ is in the form \eqref{drassi} and
has at least $\sigma+1$ factors in total.
\end{lemma}

  In Section~\ref{appA} we prove the remaining claims for
Lemmas~2.5.4, 2.5.5 from Chapter~2 in~\cite{alexakis}.

{\bf Recall the setting and claim of Lemma 2.5.4  in \cite{alexakis}.}
We recall some def\/initions. Recall f\/irst that given a
$P(g)$ for which $\int_{M^n}P(g)dV_g$ is a global conformal invariant,
we let $I^s_g(\phi)$ to be the $s^{th}$ {\it conformal
variation} of $P(g)$: $I^s_g(\phi):=\frac{d^s}{dt^s}|_{t=0} e^{nt\phi}P(e^{2t\phi}g)$.
 $I^s_g(\psi_1,\dots,\psi_s)$ is obtained from  $I^s_g(\phi)$
via complete polarization:
\[
I^s_g(\psi_1,\dots,\psi_s):=
\frac{d^s}{d\lambda_1\cdots d\lambda_s}\Big|_{\lambda_1=0,
\dots,\lambda_s=0}I^s_g(\lambda_1\cdot\psi_1
+\cdots+\lambda_s\cdot \psi_s).
\]

We recall that since $\int_{M^n}P(g)dV_g$ is a
global conformal invariant,
\[
\int_{M^n} I^s_g(\psi_1,\dots,\psi_s)dV_g=0
\]
for {\it all} $\psi_1,\dots,\psi_s\in {\cal C}^\infty (M^n)$
and {\it all} metrics $g$. We also recall that by virtue
of the transformation laws of the curvature tensor and
the Levi-Civita connections under conformal transformations,
$I^s_g(\psi_1,\dots,\psi_s)$ can be expressed as
a linear combination of complete contractions in the form:
\begin{gather}
{\rm contr}\big(\nabla^{f_1\dots f_y}\nabla^{(m_1)}R_{ijkl}\otimes\dots
\otimes \nabla^{g_1\dots g_p}\nabla^{(m_r)}R_{ijkl}
\nonumber\\
\qquad{} \otimes \nabla^{y_1\dots
y_w}\nabla^{(d_1)}{\rm Ric}_{ij}\otimes\dots\otimes\nabla^{x_1\dots
x_p}\nabla^{(d_q)}{\rm Ric}_{ij}
\nonumber\\ \qquad{}\otimes \nabla^{a_1\dots a_{t_1}}\nabla^{(u_1)}\psi_1\otimes\dots\otimes \nabla^{c_1\dots
c_{t_s}}\nabla^{(u_s)}\psi_s\big),\label{linisymric}
\end{gather}
with the conventions introduced in  Chapter~2 in \cite{alexakis}\footnote{Note
in particular that $r$ is the number of factors $\nabla^{(m)}R_{ijkl}$,
$q$ is the number of factors $\nabla^{(d)}{\rm Ric}$,
and $s$ is the number of factors $\nabla^{(u)}\psi_h$.}. In each factor
$\nabla^{f_1\dots f_y}\nabla^{(m)}_{r_1\dots r_m}R_{ijkl}$,
 each of the upper indices ${}^{f_1},\dots ,{}^{f_y}$
contracts against one of the indices ${}_{r_1},\dots ,{}_l$, while
no two of the indices ${}_{r_1},\dots ,{}_l$ contract between
themselves. On the other hand, for each factor $\nabla^{y_1\dots
y_t} \nabla^{(u)}_{a_1\dots a_u}\psi_h$,  each of the
upper indices ${}^{y_1},\dots ,{}^{y_t}$ contracts against one of
the indices ${}_{a_1},\dots ,{}_{a_u}$. Moreover, none of the
indices ${}_{a_1},\dots ,{}_{a_u}$ contract between themselves.
For the factors $\nabla^{x_1\dots x_p}\nabla^{(u)}_{t_1\dots
t_u}{\rm Ric}_{ij}$, we impose the condition that each of the upper indices
${}^{x_1},\dots ,{}^{x_p}$ must contract against one of the lower
indices ${}_{t_1},\dots ,{}_{t_u},{}_i,{}_j$. Moreover, we impose
the restriction that none of the indices ${}_{t_1},\dots
,{}_{t_u},{}_i,{}_j$ contract between themselves (this assumption
can be made by virtue of the contracted second Bianchi identity).
We recall some important def\/initions:

\begin{definition}
 \label{recall.def}
We recall that for complete contractions in the form (\ref{linisymric}), $\delta$ stands for the
number of internal contractions\footnote{We recall that an ``internal contraction''
is a pair of indices that belong to the same factor
and contract against each other.}, where we are also
counting the internal contraction $({}^a,{}_a)$
in each factor $\nabla^{(p)}{\rm Ric}_{ij}=\nabla^{(p)}{R^a}_{iaj}$
and the two internal contractions in each factor $R={R^{ab}}_{ab}$.
We also recall that for such complete contractions, $|\Delta|$ stands for the total
number of factors~$\Delta\psi_h$ and~$R$ in~(\ref{linisymric}).
We also recall that a complete contraction in the form (\ref{linisymric})
is called ``acceptable'' if each factor $\psi_h$ is dif\/ferentiated at least twice.
Also, partial contractions in he form (\ref{linisymric}) are called
acceptable is each $\psi_h$ is dif\/ferentiated at least twice\footnote{For future reference,
we also recall that  a partial contraction with one free index is also called a vector f\/ield.}.
\end{definition}

Furthermore, we recall that among the complete contractions ${\rm contr}(\cdots)$
in $I^s_g(\psi_1,\dots,\psi_s)$ with the minimum number $\sigma$ of factors in total,
$\mu$ is the minimum number of internal contractions.

We also recall  a def\/inition we have introduced for
complete contractions of in the form~(\ref{linisymric}):

\begin{definition}
\label{proka}
Consider any
complete contraction $C_g(\psi_1,\dots,\psi_a)$ in the form (\ref{linisymric})
with $\sigma$ factors.
If $C_{g}(\psi_1,\dots ,\psi_s)$ has $q=0$ and
$\delta =\mu$ it will be called a target. If $C^l_{g}(\psi_1,\dots
,\psi_s)$ has $q=0$ and $\delta >\mu$, it will be called a
contributor.

 If $C_{g}(\psi_1,\dots ,\psi_s)$ has $q>0$ and
$\delta > \mu$ we call it 1-cumbersome. We call
$C_{g}(\psi_1,\dots ,\psi_s)$ 2-cumbersome if it has $q>0$ and
$\delta = \mu$ and the feature that each
 factor~$\nabla^{a_1\dots a_t}\nabla^{(p)}_{r_1\dots r_p}{\rm Ric}_{ij}$ has $t>0$ and
the index ${}_j$ is contracting against one of the indices
${}^{a_1},\dots ,{}^{a_t}$.

  Finally, when we say $C_g(\psi_1,\dots,\psi_s)$ is
``cumbersome'', we will mean it is either 1-cumbersome or
2-cumbersome.
\end{definition}

\begin{definition}
\label{proka'}
We denote by $\sum_{j\in J} a_j C^j_g(\psi_1,\dots ,\psi_s)$
a generic linear combination of contributors.
We denote by  $\sum_{f\in F} a_f
C^f_{g}(\psi_1,\dots ,\psi_s)$
a generic linear combination
of cumbersome complete contractions.
\end{definition}

We have then derived in Chapter 2 in \cite{alexakis} that
$I^s_{g}(\psi_1,\dots ,\psi_s)$ can be expressed in the form:
\begin{gather}
 I^s_{g}(\psi_1,\dots ,\psi_s)=\sum_{l\in L_\mu} a_l
C^{l,\iota}_{g}(\psi_1,\dots ,\psi_s)+ \sum_{j\in J} a_j
C^j_{g}(\psi_1,\dots ,\psi_s)
\nonumber\\
\phantom{I^s_{g}(\psi_1,\dots ,\psi_s)=}{} + \sum_{f\in F} a_f
C^f_{g}(\psi_1,\dots ,\psi_s)+({\rm Junk}). \label{polis}
 \end{gather}
(We recall from Chapter 2 in \cite{alexakis} that $\sum_{l\in L_\mu} a_l
C^{l,\iota}_{g}(\psi_1,\dots ,\psi_s)$ in the above is a {\it specific}
linear combination of targets, which is in precise correspondence
with a {\it specific} linear combination of terms in~$P(g)$;
the precise form of this correspondence is
not important for the present paper,
so we do not dwell on this further.)

We also recall from Chapter 2 in \cite{alexakis}
that for each $q$, $1\le q\le \sigma-s$, $F^q\subset F$
stands for the index set of complete contractions with precisely~$q$ factors $\nabla^{(p)}{\rm Ric}$ or~$R$. Recall that for each
complete contraction in the form (\ref{linisymric}) we have
denoted by~$|\Delta|$ the number of factors in one of the forms
$\Delta\psi_h$, $R$. For each index set $F^q$ above, let us denote
by $F^{q,*}\subset F^q$ the index set of complete contractions
with $|\Delta|\ge\sigma-2$, and $F^{*}=\bigcup _{q>0}F^{q,*}$.

{\bf Important remark.} We recall a remark made in
Chapter 2 in \cite{alexakis} where we noted that
any $C^f_g(\psi_1,\dots,\psi_s)$ with $\alpha>0$
factors $R$ (of the scalar curvature) will have $\delta\ge \mu+2\alpha$.

The missing claim needed to derive Lemma 2.5.4  in \cite{alexakis} is then the following:

\begin{lemma}
 \label{postpone1st}
 There exists a linear combination of vector fields
$($indexed in $H$ below$)$, each in the form \eqref{linisymric} with
$\sigma$ factors\footnote{The vector fields in question are partial
contractions in the form (\ref{linisymric})
with one free index and with weight $-n+1$.},
so that modulo complete contractions of length
$>\sigma$:
\begin{gather*}
\sum_{f\in F^{*}} a_f
C^f_g(\psi_1,\dots,\psi_s)-{\rm div}_i\sum_{h\in H} a_h
C^{h,i}_g(\psi_1,\dots,\psi_s)= \sum_{y\in Y} a_y
C^y_g(\psi_1,\dots,\psi_s),
\end{gather*}
where the complete contractions indexed in~$Y$ are in the form
\eqref{linisymric} with length~$\sigma$, and satisfy all the
properties of the sublinear combination $\sum_{f\in F}\cdots$ but
in addition have $|\Delta|\le\sigma-3$.
\end{lemma}

{\bf The setting and remaining claims for Lemma 2.5.5  in \cite{alexakis}.}
Recall that in the setting of Lemma 2.5.5, our point of
reference is a linear combination:
 \begin{gather*}
Y_{g}(\psi_1,\dots ,\psi_s)=\sum_{l\in L_\mu} a_l
C^{l,\iota}_{g}(\psi_1,\dots ,\psi_s)+ \sum_{j\in J} a_j
C^j_{g}(\psi_1,\dots ,\psi_s)+({\rm Junk}),
 \end{gather*}
with the same conventions as in the ones under~(\ref{polis}).
The remaining claims of Lemma 2.5.5 (the proof of which was deferred to this paper)
 are as follows:

Denote
by $L^{*}_\mu\subset L_\mu$, $J^{*}\subset J$ the index sets
of complete contractions in $Y_g$
with $|\Delta|\ge\sigma-2$, among the complete contractions
indexed in $L_\mu$, $J$ respectively.

\begin{lemma}
\label{violi}
 We claim that there
exists a linear combination of vector fields $($indexed in
$H$ below$)$, with length $\sigma$, in the form \eqref{linisymric}
  without factors $\nabla^{(p)}{\rm Ric}$, $R$  and with $\delta=\mu$ so
that:
\begin{gather*}
 \left[\sum_{l\in L^{*}_\mu} a_l
C^{l,i_1\dots i_\mu}_g(\psi_1,\dots,\psi_s)-{\rm div}_i \sum_{h\in
H}
a_hC^{h,i|i_1\dots i_\mu}_g(\psi_1,\dots,\psi_s)\right]
\nabla_{i_1}\upsilon\cdots\nabla_{i_\mu}\upsilon
\nonumber\\
\qquad{} =\sum_{l\in
\overline{L}} a_l C^{l,i_1\dots i_\mu}_g(\psi_1,\dots,\psi_s)
\nabla_{i_1}\upsilon\cdots\nabla_{i_\mu}\upsilon,
\end{gather*}
where the complete contractions indexed in $\overline{L}$ are in
the form~\eqref{linisymric} with no factors $\nabla^{(p)}{\rm Ric}$ or~$R$ and with $|\Delta|\le\sigma-3$.
\end{lemma}

 In the setting
$L^{*}_\mu=\varnothing$, what remains to be shown to complete the
proof of Lemma 2.5.5 in~\cite{alexakis} is the following:

\begin{lemma}
\label{postpone2}
 Assume that $L^{*}_\mu=\varnothing$, and $J^*$ is as above.
 We then claim that there exists a linear combination of
vector fields $($indexed in $H$ below$)$ so that:
\begin{gather*}
\sum_{j\in J^{*}} a_j C^j_g(\psi_1,\dots,\psi_s)-{\rm div}_i
\sum_{h\in H} a_h C^{h,i}_g(\psi_1,\dots,\psi_s)=\sum_{y\in Y'}
a_y C^y_g(\psi_1,\dots,\psi_s),
\end{gather*}
where the complete contractions indexed in $Y'$ are in the form
\eqref{linisymric} with length $\sigma$, with no factors
$\nabla^{(p)}{\rm Ric}$ or $R$ and have $\delta\ge\mu+1$ and in addition
satisfy $|\Delta|\le\sigma-3$.
\end{lemma}

{\bf The setting and claims of Lemma 2.5.3  in \cite{alexakis}.}
Recall that $P(g)=\sum_{l\in L} a_l C^l(g)$ is assumed to be a
linear combination of
complete contractions in the form (\ref{drassi}). Recall that
 $\sigma$ is the minimum number of factors (in total)
 among all complete contractions $C^l(g)$; $L_\sigma\subset L$ is the corresponding index set.
Also, $s>0$ is the maximum number of factors $\nabla^{(a)}P$ among
the complete contractions $C^l(g)$ indexed in $L_\sigma$; we
denote the corresponding index set by $\Theta_s\subset L'$.
We have def\/ined $P(g)|_{\Theta_s}:=\sum_{l\in \Theta_s} a_l C^l(g)$.

{\bf Special def\/inition.}   If
 $s=\sigma-2$ then $P(g)|_{\Theta_s}$ is ``good'' if the only complete
 contraction in $P(g)|_{\Theta_s}$ with $\sigma-2$ factors $P^a_a$
 is of the form $({\rm const})\cdot {\rm contr}(\Delta^{\frac{n}{2}-\sigma-2}\nabla^{il}
 W_{ijkl}\otimes\nabla^{i'l'}{{W_{i'}}^{jk}}_{l'}\otimes(P^a_a)^{\sigma-2})$
 (when $\sigma<\frac{n}{2}-1$)
 or $({\rm const})\cdot {\rm contr}(\nabla^lW_{ijkl}\otimes\nabla_{l'}W^{ijkl'}\otimes
 (P^a_a)^{\sigma-2})$ when $\sigma=\frac{n}{2}-1$. If $s=\sigma-1$, then
$P(g)|_{\Theta_s}$ is ``good'' if all complete contractions in $P(g)|_{\Theta_s}$
have $\delta_W+\delta_P=\frac{n}{2}-1$.\footnote{In other words, if there are complete contractions in
$P(g)|_{\Theta_s}$ with $\delta_W+\delta_P<\frac{n}{2}-1$
 then $P(g)|_{\Theta_s}$ is ``good'' if no complete contractions
in $P(g)|_{\Theta_s}$
have $\sigma-2$ factors $P^a_a$.}

Lemma 2.5.3  in  \cite{alexakis} claims:
\begin{lemma}
\label{tistexnes}
There exists a divergence ${\rm div}_iT^i(g)$ so that
\[
P(g)|_{\Theta_s}-{\rm div}_iT^i(g)=\sum_{l\in \Theta'_s} a_l
C^l(g)+\sum_{t\in T} a_t C^t(g).
\] Here each $C^t(g)$ is in the
form \eqref{drassi} and has fewer than $s$ factors $\nabla^{(p)}P$.
The complete contractions indexed in $\Theta'_s$ are in the form
\eqref{drassi} with $s$ factors $\nabla^{(p)}P$ and moreover this
linear combination is ``good''.
\end{lemma}

\section{Proof of Lemma \ref{theatro}}
\label{sequals2}

 We f\/irst observe that the claim is trivial when $\sigma=1$.
In that case (modulo applying a curvature identity and introducing
longer correction terms), the sublinear combination $P(g)|_1$ will
be $P(g)|_1=({\rm const})\cdot \Delta^{\frac{n}{2}-1}R$ ($R$ is the
scalar curvature). Thus, $P(g)|_1$ can be written as $({\rm const}) \nabla^i(\nabla_i\Delta^{\frac{n}{2}-2}R)$.

  The case $\sigma=2$ is dealt with by explicitly constructing divergences of
vector f\/ields and long  calculations of {\it one} local conformal invariant,
using the Fef\/feman--Graham ambient metric, \cite{fg:ci, fg:latest}\footnote{We refer
the reader to Chapter~2 in \cite{alexakis} for a review
 of the ambient metric and and of the algorithm we employ for computations.}. We f\/irst
consider the terms in $P(g)|_2$ with two factors
$\nabla^{(m)}P_{ab}$. We will show that we can explicitly  construct a divergence ${\rm div}_i
T^i(g)$ so that:
\begin{gather}
\label{simple.claim}
P(g)|_2-{\rm div}_i
T^i(g)=({\rm const})|\nabla^{(\frac{n}{2}-2)}P|^2+\sum_{h\in H} a_h
C^h(g)+\sum_{w\in W} a_w C^w(g),
\end{gather} where the terms indexed in $H$
have one factor $\nabla^{(m)}W_{ijkl}$ and one factor
$\nabla^{(m')}P_{ab}$. The terms indexed in $W$ have two factors
$\nabla^{(m)}W_{ijkl}$. (The above holds modulo longer complete
contractions, as usual).

We explain how (\ref{simple.claim}) is proven in detail, since
the main idea will be used repeatedly throughout this section. Let us
f\/irst recall a few classical identities, which can be found in \cite{alexakis}.

{\bf Useful formulae.} Firstly, antisymmetrizing the indices~${}_{c}$,~${}_a$ in a
factor $\nabla^{(m+1)}_{r_1\dots r_mc}P_{ab}$,\footnote{This is the $(m+1)$st iterated
covariant derivative of the Schouten tensor, see \cite{alexakis}
 for details.} gives rise to a Weyl tensor:
\begin{gather}
\label{cotton}
\nabla^{(m+1)}_{r_1\dots r_mc}P_{ab}-\nabla^{(m+1)}_{r_1\dots r_ma}P_{cb}=
\frac{1}{3-n}\nabla^{(m)}_{r_1\dots r_ms}{W_{cab}}^s.
\end{gather}

We also recall that the indices ${}_i$, ${}_j$ and ${}_k$, ${}_l$ in each
 tensor $\nabla^{(m)}_{r_1\dots r_M}W_{ijkl}$ are anti-symmetric.
 Finally, we recall the ``fake'' second Bianchi identities
 from Chapter 2 in~\cite{alexakis}. These are substitutes for the second Bianchi identity for the tensor
 $\nabla_aW_{ijkl}$.

 Now, consider any complete contraction $C(g)$ involving exactly two factors,
  $\nabla^{(m)}P_{ab}$ and $\nabla^{(m')}P_{a'b'}$. So
   $C(g)={\rm contr}(\nabla^{(m)}P_{ab}\otimes\nabla^{(m')}P_{a'b'})$. We f\/irstly show that by
    subtracting an explicitly constructed divergence\footnote{As allowed in
    the Deser--Schwimmer conjecture.}
    ${\rm div}_iT^i(g)=\sum_{h=1}^K {\rm div}_i C^{K,i}(g)$ we can write:
    \begin{gather}
    \label{simpl1}
    C(g)-{\rm div}_iT^i(g)=C'(g)
    \end{gather}
    (modulo terms with more than two factors),
where $C'(g)$ is some complete contraction involving two factors
 $T_1=\nabla^{(\frac{n}{2}-2)}P_{ab}$, $T_2=\nabla^{(\frac{n}{2}-2)}P_{a'b'}$
 with the additional property that each of the $\frac{n}{2}$ indices in the f\/irst
  factor contract against an index in the second factor,
  and vice versa.

\begin{proof}[Proof of (\ref{simple.claim}).] We construct the divergences needed for (\ref{simpl1}).
Let us suppose that $C(g)$ has~$M$ pairs of indices $({}_s,{}_t)$ which
belong to the same factor and contract against each other; we call such
contractions ``internal contractions''. If $M=0$ $C(g)$
is in the desired form; thus we may assume that $M>0$. Firstly, by
 using the second contracted Bianchi identity $\nabla_cP^a_a=\nabla_aP^a_c$,
 we may assume that if one of the factors
  $\nabla^{(m)}P_{ab}$ in $C(g)$ has $m>0$, then the indices ${}_a$, ${}_b$
  are {\it not} contracting against each other.
After this, we construct the divergence by an iterative procedure. Pick out
 any pair of indices ${}_s$, ${}_t$ in $C(g)$ which contract against each
 other and belong to the same factor; assume without loss of generality  that ${}_s$ is a derivative
 index\footnote{This assumption can be made
 by virtue of the previous sentence and since we are dealing
  with complete contractions of weight $-n$.}. We then construct a
  partial contraction $C^{1,i}(g)$ out of $C(g)$ by {\it erasing}
  the (derivative) index ${}_s$ and
  making the index ${}_t$ into a {\it free index} ${}^i$.
Observe then that (modulo complete contractions with at least three factors),
 $C(g)-{\rm div}_iC^{1,i}(g)=\overline{C}(g)$, where
 $\overline{C}(g)$ has $M-1$ ``internal contractions''.
 Iterating this step $M-1$ more times, we derive~(\ref{simpl1}).

Now, consider $C'(g)$ and consider
the indices ${}_a$, ${}_b$ in the factor $T_1=\nabla^{(\frac{n}{2}-2)}P_{ab}$.
They contract against two indices (say ${}_c$, ${}_d$) in the factor $T_2$.
We then apply the curvature identity repeatedly
 and also the identity (\ref{cotton}) to arrange that the indices ${}_a$,~${}_b$ in $T_1$
 contract against the indices~${}_c$,~${}_d$ in $T_2$. The correction terms that
 arise by the application of the curvature identity have three factors.
 The correction terms arising from (\ref{cotton})
 will be in the generic form $\sum_{h\in H} a_h C^h(g)$
 as described below~(\ref{simple.claim}).
 \end{proof}

  Now, we will show that $({\rm const})=0$ in (\ref{simple.claim}). We derive
this easily. Consider $I^2_g(\phi)$ $(:=\frac{d^2}{dt^2}|_{t=0}[e^{nt\phi}P(e^{2t\phi}g)])$
 and apply the super divergence
formula to the above\footnote{See the algorithm at  the end of \cite{a:dgciI}.}.
 We derive that $({\rm const})\cdot
|\nabla^{(\frac{n}{2})}\phi|^2=0$ (modulo longer terms).
Therefore, $({\rm const})=0$.

  Thus, we may assume without loss of generality  that
\[
P(g)|_2=\sum_{h\in H} a_h
C^h(g)+\sum_{w\in W} a_w C^w(g),
\]
with the same conventions introduced under (\ref{simple.claim}).

We next claim that we can explicitly construct a
divergence ${\rm div}_i T'^i_g$ such that:
\begin{gather*}
\sum_{h\in H} a_h
C^h(g)-{\rm div}_i T'^i(g)=\sum_{w\in W} a_w C^w(g),
\end{gather*} where the terms
indexed in $W$ are of the generic form described above.
The divergence needed for the above is constructed
by the same technique as for (\ref{simpl1}): In each $C^h(g)$ we
iteratively pick out the internal contractions in the factor $\nabla^{(m)}W_{ijkl}$,
erase one (derivative) index in that contraction,
(thus obtaining a partial contraction with one free index)
and subtract the divergence of that partial
contraction\footnote{Sometimes, by abuse of language,
we will refer to this subtraction of an
explicit divergence as an ``integration by parts''.}.
After repeating this process
 enough times so that  no internal contraction is
  left in the factor $\nabla^{(m)}W_{ijkl}$, we end up with a formula:
 \begin{gather*}
 \sum_{h\in H} a_hC^h(g)-{\rm div}_i T'^i(g)=\sum_{r\in R}a_r C^r(g).
 \end{gather*}
Here  the complete contractions $C^r(g)$ in the r.h.s.\ are in the form
${\rm contr}(\nabla^{(m)}W_{ijkl}\otimes\nabla^{(t)}P_{ab})$, where
 {\it all} $m+4$ indices in the factor $\nabla^{(m)}W_{ijkl}$
 contract against the other factor $\nabla^{(t)}P_{ab}$. Now, since the indices~${}_a$,~${}_b$ in Schouten tensor are symmetric and the indices~${}_i$,~${}_j$ and~${}_k$,~${}_l$ are antisymmetric,
   it follows from (\ref{cotton}) and from the curvature identity that:
\[
\sum_{r\in R}a_r C^r(g)=\sum_{w\in W} a_w C^w(g),
\]
modulo complete contractions with three factors, which arise due to the curvature identity.

 So we may
assume without loss of generality  that $P(g)|_2$ consists of terms with two factors
$\nabla^{(m)}W_{ijkl}$.

We can then again explicitly construct  a divergence ${\rm div}_iT''^i(g)$, in
order to write (modulo longer complete contractions):
\[
P(g)|_2-{\rm div}_iT''^i(g)=\sum_{w\in W'} a_w C^w_g,
\]
where the terms in the r.h.s.\ have the additional feature that none
of the two factors $\nabla^{(m)}W_{abcd}$ have internal
contractions. The divergence ${\rm div}_iT''^i(g)$ is constructed by the iterative
procedure used to prove (\ref{simpl1}).
 So we may assume without loss of generality  that all
terms in $P(g)|_2$ have this property.

  Finally, by just keeping track of the
correction terms of the form $\nabla^sW_{sdfh}\otimes g$
in the ``fake'' second Bianchi identities\footnote{These
are presented in Chapter 2 in \cite{alexakis}.},
 we will prove that we can can explicitly
 construct a divergence ${\rm div}_iT^i(g)$ such that:
\begin{gather}
\label{namen}
 P(g)|_2=({\rm const})'C^{*}(g)+{\rm div}_iT^i(g),
\end{gather}
where $C^{*}(g)$ is the complete contraction
$|\nabla^{(\frac{n}{2}-2)}W_{abcd}|^2$.

\begin{proof}[Proof of (\ref{namen}).] By virtue of the anti-symmetry of the indices~${}_i$,~${}_j$
 and~${}_k$, ${}_l$ and of the f\/irst Bianchi identity in the Weyl tensor $W_{ijkl}$, we
see that (modulo introducing correction terms
with three factors), $P(g)|_2$
can be expressed in the form:
\begin{gather*}
P(g)|_2=a\cdot\big|\nabla^{(\frac{n}{2}-2)}W_{abcd}\big|^2+ b\cdot
C^2(g)+c\cdot C^3(g),
\end{gather*}
where $C^2(g)$ is the complete contraction:
\[
{\rm contr}\big(\nabla_{r_1\dots r_{\frac{n}{2}-3} s}W_{tjkl}\otimes
\nabla^{r_1\dots r_{\frac{n}{2}-3} t}W^{sjkl}\big)
\] while $C^3(g)$ is
the complete contraction:
\[
{\rm contr}\big(\nabla_{r_1\dots r_{\frac{n}{2}-4} su}W_{tjky}\otimes
\nabla^{r_1\dots r_{\frac{n}{2}-4} ty}W^{sjku}\big).
\]

 We then only have to apply the ``fake'' second Bianchi
identities from Chapter~2 in \cite{alexakis} to derive that we can write:
\begin{gather}
\label{korusxades}
C^2(g)=\frac{1}{2}\big|\nabla^{(\frac{n}{2}-2)}W_{abcd}\big|^2+
\frac{1}{n-3}\big|\nabla^{(\frac{n}{2}-3)}\nabla^sW_{sjkl}\big|^2.
\end{gather}
By virtue of (\ref{korusxades}) we easily derive that we can
construct a divergence ${\rm div}_iT^i(g)$ so that:{\samepage
\begin{gather*}
C^2(g)-{\rm div}_i
T^i(g)=\frac{n-3}{2(n-4)}\big|\nabla^{(\frac{n}{2}-2)}W_{abcd}\big|^2
\end{gather*}
(notice the constant is strictly positive)\footnote{The
 divergence is constructed in the same way as in the proof of (\ref{simpl1}). We just
 ``integrate by parts'' the two internal
 contractions ${}^s,{}_s$ in the two factors in $|\nabla^{(\frac{n}{2}-3)}\nabla^sW_{sjkl}|^2$.}.}

In order to deal with the complete contraction $C^3(g)$, we will
note another useful identity. Let $\overline{C}^3(g)$ stand for
the complete contraction $({\nabla_{r_1\dots
r_{\frac{n}{2}-4}a}}^lW_{ijkl}\otimes \nabla^{r_1\dots
r_{\frac{n}{2}-4}is}{W^{ajk}}_s)$; we also denote by $C'^3(g)$ the
complete contraction:   $(\nabla_{r_1\dots
r_{\frac{n}{2}-3}s}{W_{ijk}}^s \otimes\nabla^{r_1\dots
r_{\frac{n}{2}-3}t} {W^{ijk}}_t)$. We then calculate:
\begin{gather}
\label{korusxades2} \overline{C}^3(g)=\frac{1}{2}C'^3(g).
\end{gather}

Thus, using the above we derive as before
 that we can explicitly construct a divergence
${\rm div}_i T'^i(g)$ such that:
\begin{gather}
\label{korusxades4} C^3(g)-{\rm div}_i T'^i(g)=\frac{1}{2}C^2(g).
\end{gather}

 Therefore, we derive that modulo subtracting a divergence
${\rm div}_iT^i(g)$ from $P(g)|_2$, we may assume that $P(g)|_2$ is in
the form $P(g)|_2=({\rm const})\cdot |\nabla^{(\frac{n}{2}-2)}W|^2$.
\end{proof}

{\bf A construction in the Fef\/ferman--Graham ambient metric.}
 The reader is referred to the discussion on the Fef\/ferman--Graham ambient metric
from Subsection~3.2.1 in \cite{alexakis}. We consider the
local conformal invariant:
$\tilde{C}(g)={\Delta}^{\frac{n}{2}-2}_{\tilde{g}}|\tilde{R}|_{\tilde{g}}^2$.
Here we think of $|\tilde{R}|^2_{\tilde{g}}$ as the product
$\tilde{g}^{a\alpha}
\tilde{g}^{b\beta}\tilde{g}^{c\gamma}\tilde{g}^{d\delta}
\tilde{R}_{abcd}\tilde{R}_{\alpha\beta\gamma\delta}$,
where $\tilde{R}$ is the ambient curvature tensor, $\tilde{g}$
is the ambient metric tensor and $\Delta_{\tilde{g}}$
 is the Laplace--Beltrami operator associated to $\tilde{g}$.

 We will show that there exists an explicit divergence
${\rm div}_i T^i(g)$, so that
modulo terms of length $\ge 3$:
\begin{gather}
\label{motlagh} \tilde{C}(g)=({\rm const})C^{*}(g)+{\rm div}_i T^i(g).
\end{gather}
(Here the constant $({\rm const})$ is non-zero.) If we can prove the
above then in view of (\ref{namen})  our claim will clearly follow.

\begin{proof}[Proof of (\ref{motlagh}).] We recall certain facts about the
ambient metric The reader is referred to~\cite{fg:latest},
or Chapter~2 in~\cite{alexakis}. Firstly recall that given any
point $P\in M$ and coordinates $\{x^1,\dots, x^n\}$, there is a canonical
coordinate system $\{t=x^0,x^1,\dots, x^n,\rho\}$
around the image $\tilde{P}$ of $P$ in the
ambient manifold $(\tilde{G}, \tilde{g})$. Recall that the local conformal invariant is evaluated
at the point $\tilde{P}$, and $t(\tilde{P})=1,\rho(\tilde{P})=0$.

Recall also that the vector f\/ields $X^0,X^1,\dots, X^n,X^\infty$
are the coordinate vector f\/ields $\frac{\partial}{\partial x^0}, \dots ,$
$\frac{\partial}{\partial x^n}, \frac{\partial}{\partial\rho}$.
Further, 
when we give values $0,1,\dots, n,\infty$ to indices of tensors
 that appear further down, these values correspond to the coordinate frame above.

Now, we recall some basic facts regarding the Taylor expansion
 of the ambient metric~$\tilde{g}$ at~$\tilde{P}$.
 Firstly, that $\partial_\infty
\tilde{g}^{ab}=-2P^{ab}$ if $1\le a,b\le n$ and $\partial_\infty
\tilde{g}^{\infty\infty}=-2$. These are the only non-zero
components of the matrix $\partial_\infty \tilde{g}^{cd}$ (with
raised indices).

 Now, in order to state our next claim, we recall that
$\tilde{\Gamma}_{ab}^c$ stand for the Christof\/fel symbols of the
ambient metric. We recall that if $0\le a,b,c\le n$ then
$\partial_\infty \tilde{\Gamma}_{ab}^c=F_{ab}^c(R)$, where the
expression in the r.h.s.\ stands for a tensor (in the indices~${}_a$,~${}_b$,~${}^c$) involving {\it at least one $($possibly internally
contracted$)$ factor of the curvature tensor}. On the other hand, we
also recall that if $1\le a,b,c\le n$ then
$\tilde{\Gamma}_{ab}^c=\Gamma_{ab}^c$ (the r.h.s.\ stands for the
Christof\/fel symbol of the metric $g$). Finally, if $1\le a,b\le n$
then $\tilde{\Gamma}_{ab}^\infty=-g_{ab}$. Furthermore, we will use the fact that
 $\partial_\infty \tilde{g}_{ab}=2P_{ab}$ and  also the
formula~(3.21) from~\cite{fg:latest}:
\[
\partial^{s}_\infty \tilde{g}_{ij}=\frac{2}{(4-n)\cdots (2s-n)}
\big[\Delta^{s-1}P_{ij}-\Delta^{s-2}\nabla_{ij}P^a_a\big] +Q(R)
\] (for
$s\ge 2$; all other components $\partial^{(s)}_{\infty\dots\infty}\tilde{g}_{Ab}$
vanish if $A=0,\infty$ and $s\ge 2$; $b=0,1,\dots, ,n,\infty$)\footnote{$Q(R)$ stands for a linear combination of
partial contractions involving at least two curvature terms, as
usual.}. We calculate that for $1\le i,j,k,l\le n$ and $1\le
\alpha\le \frac{n}{2}-2$:
\begin{gather}
\label{robingr1}
\partial^\alpha_\infty \tilde{R}_{ijkl}=(-1)^{\alpha-1} \frac{1}{(n-3)(n-4)\cdots (n-2\alpha)}
\Delta^{\alpha-1}\big[\nabla^{t}_jW_{tikl}-\nabla^{t}_iW_{tjkl}\big] +Q(R)
\end{gather}
(where if $\alpha=1$ then the constant above is $\frac{1}{n-3}$).
Furthermore for $0\le\alpha\le \frac{n}{2}-3$:
\begin{gather}
\label{robingr2}
\partial^\alpha_\infty \tilde{R}_{\infty jkl}=(-1)^{\alpha-1} \frac{1}{(n-3)(n-4)\cdots(n-2-2\alpha)}
\Delta^\alpha \nabla^sW_{sjkl}+Q(R)
\end{gather}
(where if $\alpha=0$ then the constant above is $\frac{1}{n-3}$).
Moreover for $0\le\alpha\le \frac{n}{2}-4$:
\begin{gather}
\label{robingr3}
\partial^\alpha_\infty \tilde{R}_{\infty jk\infty}=(-1)^{\alpha} \frac{1}{(n-3)(n-4)\cdots (n-4-2\alpha)}
\Delta^\alpha \nabla^{il}W_{ijkl}+Q(R).
\end{gather}
(The left-hand sides of the above are also known as {\it Graham's extended obstruction tensors}, see~\cite{graham:obstrten}
for a detailed study of these tensors.)

Next, two calculations. Let $F(\tilde{g})$ be any Riemannian invariant in the ambient metric $\tilde{g}$,
with weight $w$. (In particular, $F(\tilde{g})$ will have homogeneity~$-w$ in~$t$.)
 Then we calculate:
\begin{gather}
\label{olozontano}
\partial^k_\infty \Delta_{\tilde{g}}[F(\tilde{g})]=\Delta_g L(g)+
(2w+n-2k)\partial^{k+1}_\infty F(\tilde{g})+Q(R).
\end{gather}
$L(g)$ is a Riemannian invariant of the metric $g$, {\it not the ambient metric $\tilde{g}$},
with weight $w-2k$. $Q(R)$ is a linear combination of complete contractions in the
iterated covariant derivatives of the curvature tensor, and
each complete contraction involves at least two such curvature terms.

A word regarding the derivation of the above: By the form of the ambient metric
(see in particular page 20 in \cite{fg:latest}), we derive
that at any point on the ambient manifold,
\begin{gather*}
\Delta_{\tilde{g}}F(\tilde{g})=\sum_{i,j=1}^n\tilde{g}^{ij}\nabla^{(2)}_{ij}
F(\tilde{g})+2t^{-1}\tilde{\nabla}^{(2)}_{0\infty}F(\tilde{g})
-2t^2\rho\tilde{\nabla}^{(2)}_{\infty\infty}F(\tilde{g})
\nonumber\\
\phantom{\Delta_{\tilde{g}}F(\tilde{g})}{}
=\sum_{i,j=1}^n\tilde{g}^{ij}\left(\partial^{(2)}_{ij}
-\sum_{k=1}^n\tilde{\Gamma}_{ij}^k\partial_k\right)
F(\tilde{g})-\left[\sum_{i,j=1}^n\tilde{g}^{ij}\tilde{\Gamma}_{ij}^0\partial_0+
\sum_{i,j=1}^n\tilde{g}^{ij}\tilde{\Gamma}_{ij}^\infty\partial_\infty\right]F(\tilde{g})
\nonumber\\
\phantom{\Delta_{\tilde{g}}F(\tilde{g})=}{}
{} +2wt^{-1}\partial_\infty F(\tilde{g})
-2t^2\rho\partial^{(2)}_{\infty\infty}F(\tilde{g}).
\end{gather*}
Thus, if we take the $k^{\rm th}$ derivative $\partial^{(k)}_{\infty\dots \infty}$ of the above
equation and then evaluate at $t=1$, $\rho=0$, we obtain the r.h.s.\
of (\ref{olozontano}) as follows: the term $\Delta_g L(g)$
arises when all $k$ derivatives~$\partial_\infty$ hit the factor $F(\tilde{g})$ in
$\sum_{i,j=1}^n\tilde{g}^{ij}\partial^{(2)}_{ij}F(\tilde{g})$.
The coef\/f\/icient $2w$ arises from
the term $\partial^{(2)}_{0\infty} F(\tilde{g})$
in~$\tilde{\nabla}^{(2)}_{0\infty}F(\tilde{g})$, due to the
homogeneity $w$ of $F(\tilde{g})$ in $x^0=t$. The coef\/f\/icient $+n$ arises from the
term $-\sum_{i,j=1}^n\tilde{g}^{ij}\tilde{\Gamma}_{ij}^\infty\partial_\infty$.
Finally, the term $-2k$ arises when exactly one of the
 $k$ derivatives~$\partial_\infty$ hits the coef\/f\/icient~$-2t^2\rho$ of the expression $\partial^{(2)}_{\infty\infty}F(\tilde{g})$.
All other terms that arise are in the form~$Q(R)$.

  Now, we can iteratively apply the above
 formula to obtain a useful expression for
$\Delta^{\frac{n}{2}{-}2}_{\tilde{g}}\!|\tilde{R}_{ijkl}|^2_{\tilde{g}}$;
we consider $\Delta_{\tilde{g}}(\Delta_{\tilde{g}}(\dots (|\tilde{R}_{ijkl}|^2_{\tilde{g}})\dots )$,
and we replace each of the $\Delta_{\tilde{g}}$'s according to (\ref{olozontano}),
{\it from left to right}. The resulting equation is:
\begin{gather*}
\tilde{C}(g)=({\rm const})'\cdot
\partial^{\frac{n}{2}-2}_\infty
(|\tilde{R}_{ijkl}|^2_{\tilde{g}})+{\rm Cubic}(R)+\Delta_g L(g).
\end{gather*}
Here the constant $({\rm const})'$ is {\it non-zero};
this is because at each application of the identity
(\ref{olozontano}) we have $-w+2k+2=n$,
and $w$ takes on the values $-n+2, \dots, -4$; thus each factor in the product is non-zero.
 ${\rm Cubic}(R)$ is a linear
combination of Riemannian invariants with at least three factors.
$L(g)$ is a Riemannian invariant (of the base metric~$g$)
of weight $-n+2$.

 Now, using formulas (\ref{robingr1}),
(\ref{robingr2}), (\ref{robingr3}) and also the
formula $\partial_\infty \tilde{g}^{\infty\infty}=-2$, we derive
that:
\begin{gather}
 \partial^{\frac{n}{2}-2}_\infty
\big(|\tilde{R}_{ijkl}|^2_{\tilde{g}}\big)= (-1)^{\frac{n}{2}}\Bigg\{ 2\cdot
2\frac{1}{(n-3)(n-4)\cdots
4}W_{ijkl}\Delta^{\frac{n}{2}-3}\nabla^i_tW^{tjkl}
\nonumber\\
\qquad{}
+2\sum_{x=1}^{\frac{n}{2}-3}{{\frac{n}{2}-2}\choose{x}}\left[\frac{1}{(n-3)(n-4)\cdots
(n-2x)}\right]
\nonumber\\
\qquad{}
\times \left[\frac{1}{(n-3)(n-4)\cdots(4+2x)}\right]
 \Delta^x\nabla^{as}W_{sjkl}\Delta^{\frac{n}{2}-2-x}\nabla_{at}W^{tjkl}
+\sum_{f\in F} a_f C^f(g)
\nonumber\\
\qquad{}
+\sum_{x=0}^{\frac{n}{2}-3} 2\cdot 4{{\frac{n}{2}-2}\choose{1}}\cdot
 {{\frac{n}{2}-3}\choose{x}}\cdot
\left[\frac{1}{(n-3)(n-4)\dots (n-2x-2)}\right]
\nonumber \\
 \qquad{}
 \times  \left[\frac{1}{(n-3)(n-4)\dots (4+2x)}\right]
 \Delta^x\nabla^sW_{sjkl}\Delta^{\frac{n}{2}-3-x}\nabla_tW^{tjkl}
\nonumber\\
\qquad{}
+\sum_{x=0}^{\frac{n}{2}-4}2\cdot
2{{\frac{n}{2}-2}\choose{1}}2\cdot2\cdot{{\frac{n}{2}-3}\choose{1}}\cdot
{{\frac{n}{2}-4}\choose{x}}
\nonumber\\
\qquad{}
\times \left[\frac{1}{(n-3)(n-4)\dots
(n-4-2x)}\right]\cdot\left[\frac{1}{(n-3)(n-4)\cdot
(4+2x)}\right]
\nonumber\\
\qquad{}
\times\Delta^x\nabla^{st}W_{sjtl}\Delta^{\frac{n}{2}-4-x}\nabla_{s't'}W^{s'jt'l}\Bigg\}+{\rm Cubic}(R).\label{aggelopoulou}
 \end{gather}
Here ${\rm Cubic}(R)$ stands for a generic linear combination of
partial contractions with at least three factors of the form
$\nabla^{(m)}R_{abcd}$.
 The terms indexed in $F$ are complete contractions in
the form:
\[
{\rm contr}\big(\Delta^q \nabla^{\beta a}W_{abcd}\otimes\Delta^{q'}\nabla^{b\alpha}{W_{\alpha\beta}}^{cd}\big).
\]
We will show below that by subtracting a divergence we can ``get
rid'' of such complete contractions, modulo introducing correction
terms with at least three factors.

\begin{proof}[Mini-proof of (\ref{aggelopoulou}).] This equation
follows by an iterated application of the Leibniz rule\footnote{Recall
that we think of $|\tilde{R}|^2_{\tilde{g}}$ as the product
$\tilde{g}^{a\alpha}
\tilde{g}^{b\beta}\tilde{g}^{c\gamma}\tilde{g}^{d\delta}
\tilde{R}_{abcd}\tilde{R}_{\alpha\beta\gamma\delta}$.}. Each
derivative can hit either one of the curvature tensors
$\tilde{R}_{ijkl}(\tilde{g})$, {\it or} one of the metric tensors~$\tilde{g}^{ab}$ (with raised indices).

Now, we have sums $\sum_{x=\cdots}^{\dots}$ in (\ref{aggelopoulou}), to which we will refer
to as the ``f\/irst'', ``second'' and ``third'' sum. Firstly, we
observe that the expression
$W_{ijkl}\Delta^{\frac{n}{2}-3}\nabla^i_tW^{tjkl}$ (along with its
coef\/f\/icient) arises when all $\frac{n}{2}-2$ derivatives
$\partial_\infty$  hit precisely one of the two (ambient) curvature
tensors~-- we then use the formula~(\ref{robingr1}). Secondly, we
observe that the f\/irst sum arises when all derivatives
$\partial_\infty$ are forced to hit either of the two (ambient) curvature
tensors, and moreover each curvature tensor must be hit by at
least one derivative. Thirdly, the second sum arises when exactly
one derivative $\partial_\infty$ hits one of the metric tensors
$\tilde{g}^{AB}$ (recall that the only non-zero components of
$\partial_\infty \tilde{g}^{ab}$ are
$\partial_\infty\tilde{g}^{ab}=-2P^{ab}$\footnote{When $1\le
a,b\le n$.} and $\partial_\infty g^{\infty\infty}=-2$; notice that
the f\/irst term may be discarded since it gives rise to terms with
at least three factors). Fourthly, the third sum arises when exactly two
derivatives $\partial_\infty$ hit  metric terms
$\tilde{g}^{AB}$ (a dif\/ferent term each).
\end{proof}

  Using (\ref{aggelopoulou}) in
conjunction with the formulas (\ref{korusxades}) and
(\ref{korusxades4}), we can derive (\ref{motlagh}). We can
explicitly construct a divergence ${\rm div}_iT^i(g)$  so that for every~$x$:
\begin{gather}
W_{ijkl}\Delta^{\frac{n}{2}-3}\nabla^i_tW^{tjkl}-{\rm div}_iT^i(g)
\nonumber\\
\qquad{}=(-1)^{\frac{n}{2}}\nabla^{(\frac{n}{2}-3)}_{r_1\dots
r_{\frac{n}{2}-3}}\nabla^sW_{sjkl}\otimes
\big(\nabla^{(\frac{n}{2}-3)}\big)^{r_1\dots
r_{\frac{n}{2}-3}}\nabla_tW^{tjkl},\label{pre.int}
\\
\Delta^x\nabla^{st}W_{sjtl}\otimes\Delta^{\frac{n}{2}-4-x}\nabla_{s't'}W^{s'jt'l}-{\rm div}_iT^i(g)
\nonumber\\
\qquad{}=(-1)^{\frac{n}{2}}\frac{1}{2}\nabla^{(\frac{n}{2}-3)}_{r_1\dots
r_{\frac{n}{2}-3}}\nabla^sW_{sjkl}\otimes
\big(\nabla^{(\frac{n}{2}-3)}\big)^{r_1\dots
r_{\frac{n}{2}-3}}\nabla_tW^{tjkl},\label{pre.int2}
\\
\Delta^x\nabla^sW_{sjkl}\otimes\Delta^{\frac{n}{2}-3-x}\nabla_tW^{tjkl}
-{\rm div}_iT^i(g)
\nonumber\\
\qquad{}= (-1)^{\frac{n}{2}-1}
\nabla^{(\frac{n}{2}-3)}_{r_1\dots
r_{\frac{n}{2}-3}}\nabla^sW_{sjkl}\otimes
\big(\nabla^{(\frac{n}{2}-3)}\big)^{r_1\dots
r_{\frac{n}{2}-3}}\nabla_tW^{tjkl},\label{pre.int3}
\\
 \Delta^x\nabla^{as}W_{sjkl}\otimes\Delta^{\frac{n}{2}-2-x}\nabla_{at}W^{tjkl}
-{\rm div}_iT^i(g)
\nonumber\\
\qquad{} = (-1)^{\frac{n}{2}}
\nabla^{(\frac{n}{2}-3)}_{r_1\dots
r_{\frac{n}{2}-3}}\nabla^sW_{sjkl}\otimes
\big(\nabla^{(\frac{n}{2}-3)}\big)^{r_1\dots
r_{\frac{n}{2}-3}}\nabla_tW^{tjkl},\label{pre.int4}
\\
\label{pre.int5}
\sum_{f\in F} a_f C^f(g)-{\rm div}_iT^i(g)=0.
\end{gather}
All the above equations hold modulo complete contractions with at
least three factors. We explain how the divergences for the
above f\/ive equations are constructed. Pick out any Laplacian
$\Delta(={\nabla^a}_a)$,\footnote{A ``Laplacian'' here means
two derivative indices that belong to the same factor
and contract against each other.} appearing in any complete contraction
above, and formally {\it erase} the upper index~${}^a$. The resulting partial
contraction is a Riemannian 1-tensor f\/ield with one free index ${}_a$.
Now, consider the divergence of this 1-tensor f\/ield and subtract
it from the original complete contraction\footnote{Further 
we refer to the
subtraction of this explicitly constructed
divergence as an ``integration by parts''.}. The result (modulo
correction terms with three factors that arise
due to the curvature identity) is a new complete contraction with one fewer Laplacian.
We iterate this step enough times, until in the end we obtain
a complete contraction with no Laplacians. Observe that
up to applying the curvature identity (and thus
introducing correction terms with three factors), the resulting
complete contraction is exactly the one claimed in (\ref{pre.int3}), (\ref{pre.int4}).
In the other cases some additional divergences need
to be subtracted. For (\ref{pre.int}) we consider the upper index $\nabla^i$
in the right factor; we erase this index (thus obtaining a Riemannian
1-tensor f\/ields) and subtract the corresponding
divergence. The result, up to applying the
curvature identity, is the r.h.s.\ of~(\ref{pre.int}).
In the case of (\ref{pre.int2}) we perform additional
``integrations by parts'' by erasing f\/irst the index~${}^s$
in the f\/irst factor, and then the index~${}_{s'}$ in the second factor.
The resulting complete contraction (after the curvature
identity and formula (\ref{korusxades2}) is the r.h.s.\ of~(\ref{pre.int2}).
Finally, to derive (\ref{pre.int5}) we integrate by parts the index ${}^b$ in the
expression
\[
\nabla^{(q+q')}_{t_1\dots t_{q'}r_1\dots r_q}\nabla^{\beta a}W_{abcd}\otimes
{\nabla^{(q+q')}}^{t_1\dots t_{q'}r_1\dots r_q}\nabla^{b\alpha}W_{\alpha\beta cd};
\]
the complete contraction we obtain will contain two dif\/ferentiated Weyl tensors,
one of which is in the form: $\nabla^b_{\dots}\nabla^aW_{abcd}$. This tensor vanishes,
modulo a quadratic  expression in curvatures, by virtue of the curvature identity.
 This concludes our proof of the above f\/ive equations.

 Therefore, using the above formulas we derive that:
\begin{gather}
 \partial^{\frac{n}{2}-2}_\infty
(|\tilde{R}_{ijkl}|^2_{\tilde{g}})-{\rm div}_iT^i(g)= (-1)^{n}\Bigg\{
  2\cdot 2\frac{1}{(n-3)(n-4)\cdots 4}
  \nonumber\\
  \qquad{} +2\sum_{x=1}^{\frac{n}{2}-3}
{{\frac{n}{2}-2}\choose{x}}\left[\frac{1}{(n-3)(n-4)\cdots (n-2x)}\right]
 \left[\frac{1}{(n-3)(n-4)\cdots(4+2x)}\right]
\nonumber\\
\qquad{}
-\sum_{x=0}^{\frac{n}{2}-3} 2\cdot 4{{\frac{n}{2}-2}\choose{1}}\cdot
 {{\frac{n}{2}-3}\choose{x}}\cdot
 \left[\frac{1}{(n-3)(n-4)\cdots (n-2x-2)}\right]
\nonumber \\
 \qquad{}\times\left[\frac{1}{(n-3)(n-4)\cdots (4+2x)}\right]
  +\sum_{x=0}^{\frac{n}{2}-4}
2{{\frac{n}{2}-2}\choose{1}}2\cdot2\cdot{{\frac{n}{2}-3}\choose{1}}\cdot
{{\frac{n}{2}-4}\choose{x}}
\nonumber\\
\qquad{}\times\left[\frac{1}{(n-3)(n-4)\cdots
(n-4-2x)}\right]\cdot\left[\frac{1}{(n-3)(n-4)\cdot (4+2x)}\right]\Bigg\}
\nonumber\\
\qquad{}\times \nabla^{(\frac{n}{2}-3)}_{r_1\dots
r_{\frac{n}{2}-3}}\nabla^sW_{sjkl}
(\nabla^{(\frac{n}{2}-3)})^{r_1\dots
r_{\frac{n}{2}-3}}\nabla_tW^{tjkl}.\label{aggelopoulou'}
 \end{gather}
The above holds modulo terms of length at least 3. Now, if we can
show that the constant $\{\cdots\}$ is strictly positive, we will
have shown (\ref{motlagh}).

We limit attention to the terms in the f\/irst sum with $x=y+2$, the
terms in the second sum with $x=y+1$ and in the third sum with
$x=y$ (for some given $y$, $0\le y\le \frac{n}{2}-4$). We observe
that:
\begin{gather*}
2{{\frac{n}{2}-2}\choose{y+2}}\left[\frac{1}{(n-3)(n-4)\cdots
(n-2(y+2))}\right] \left[\frac{1}{(n-3)(n-4)\cdots(4+2(y+2))}\right]
\nonumber\\
\qquad{} -2\cdot 4{{\frac{n}{2}-2}\choose{1}}\cdot
 {{\frac{n}{2}-3}\choose{y+1}}\cdot
 \left[\frac{1}{(n-3)(n-4)\cdots (n-2(y+1)-2)}\right]
\nonumber \\
\qquad{}\times\left[\frac{1}{(n-3)(n-4)\cdots (4+2(y+1))}\right]
\nonumber\\
 \qquad{} +2{{\frac{n}{2}-2}\choose{1}}2\cdot2\cdot{{\frac{n}{2}-3}\choose{1}}\cdot
{{\frac{n}{2}-4}\choose{y}} \left[\frac{1}{(n-3)(n-4)\cdots
(n-4-2y)}\right]
\nonumber\\
\qquad{}\times \left[\frac{1}{(n-3)(n-4)\cdot (4+2y)}\right]
\nonumber\\
\qquad{} =2\cdot \frac{(\frac{n}{2}-2)\cdots
\left(\frac{n}{2}-3-y\right)}{y!}\left[\frac{1}{(n-3)(n-4)\cdots (n-4-2y)}\right]
\nonumber\\
\qquad{}\times \left[\frac{1}{(n-3)(n-4)\cdots
(4+2y)}\right]\left\{\frac{(6+2y)(4+2y)}{(y+2)(y+1)}-4\frac{4+2y}{y+1}+4\right\}=0.
 \end{gather*}

Now, the only terms we have not taken into account are the
term $2\cdot 2\frac{2}{(n-3)(n-4)\cdots 4}$, the term in the f\/irst
sum with $x=1$ and the term in the second sum with $x=0$. But we
observe that those three terms add up to a positive number:
\begin{gather*}
4\frac{1}{(n-3)(n-4)\cdots
4}+2\left(\frac{n}{2}-2\right)\cdot\frac{1}{n-3}\frac{1}{(n-3)(n-4)\cdots 6}
\nonumber\\
\qquad{} -2\cdot 4\left(\frac{n}{2}-2\right)\frac{1}{n-3}\frac{2}{(n-3)^2(n-4)\cdots
4}=4\frac{1}{(n-3)(n-4)\cdots 4}>0.
\end{gather*}
Hence, using the above two formulas we derive that the constant in
(\ref{aggelopoulou'}) is strictly positive and thus we derive our
claim. This concludes the proof of Lemma~\ref{theatro}.
\end{proof}

\section{The proof of Lemmas \ref{postpone1st}, \ref{violi}, \ref{postpone2}, \ref{tistexnes}}
\label{appA}

 The three subsections below correspond to the cases $s<\sigma-2$, $s=\sigma-2$, $s=\sigma-1$.
  We also prove Lemma~\ref{tistexnes} in the subsections
that deal with the cases $s=\sigma-2$, $s=\sigma-1$.

\subsection[The proof of Lemmas \ref{postpone1st}, \ref{violi}, \ref{postpone2}
when $s<\sigma -2$]{The proof of Lemmas \ref{postpone1st}, \ref{violi}, \ref{postpone2}
when $\boldsymbol{s<\sigma -2}$}

{\bf Brief discussion.} Recall the discussion
regarding these lemmas in the introduction.
 Notice that since $s<\sigma-2$ we have
  $L^{*}=J^{*}=\varnothing$ by def\/inition; thus in this case
we only have to show Lemma~\ref{postpone1st}. Therefore,
  we will show that there exists a linear
  combination of tensor f\/ields, (indexed in $H$ below) so that:
\begin{gather*}
 \sum_{f\in F^{*}} a_f C^f_{g}(\psi_1,\dots,\phi_s)
-{\rm div}_i\sum_{h\in H} a_h C^{h,i}_{g}(\psi_1,\dots,\phi_s)=
 \sum_{f\in F^{OK}} a_f C^f_{g}(\psi_1,\dots,\phi_s),
\end{gather*}
where the tensor f\/ields indexed in $F^{OK}$ have all the properties of the ``cumbersome''
 tensor f\/ields, but moreover have  at most $\sigma-3$ factors in the form $R$ or $\Delta\psi_h$.
That will prove Lemma~\ref{postpone1st} in this
setting.

{\bf Rigorous discussion.} We observe that the complete
contractions in $I^s_{g}(\psi_1,\dots ,\psi_s)$ that have
$|\Delta|\ge\sigma -2$ must be  contractions $C^f$ indexed in  $F^{\sigma
-s}\bigcup F^{\sigma -s-1}\bigcup F^{\sigma -s-2}$.  We index the
complete contractions $C^f$ with $|\Delta|=\sigma -1$ in the sets
  $F^{\sigma -s}_{|\Delta|=\sigma-1}$, $F^{\sigma
-s-1}_{|\Delta|=\sigma -1}$, $F^{\sigma -s-2}_{|\Delta|=\sigma
-1}$ and the ones with $|\Delta|=\sigma -2$ in the sets $F^{\sigma
-s}_{|\Delta|=\sigma-2}$,  $F^{\sigma-s -1}_{|\Delta|=\sigma -2}$,
$F^{\sigma -s-2}_{|\Delta|=\sigma -2}$, respectively.
(Recall that the upper labels stand for the number of factors
 $\nabla^{(p)}{\rm Ric}$ or $R$ in the contractions indexed in~$F$;
 the lower labels stand for the number of factors~$R$
  plus the number of factors~$\Delta\psi_h$.)
 We claim:

\begin{lemma}\sloppy
\label{reduce} There is a linear combination of acceptable vector
fields of length $\sigma$, say $\sum_{h\in H} a_h
C^{h,i}_{g}(\psi_1,\dots ,\psi_s)$,  so that in the notation
above:
\begin{gather}  \sum_{f\in F^{\sigma
-s}_{|\Delta|=\sigma-1}\bigcup F^{\sigma -1}_{|\Delta|=\sigma
-1}\bigcup F^{\sigma -s-2}_{|\Delta|=\sigma -1}\bigcup F^{\sigma
-s}_{|\Delta|=\sigma-2}\bigcup F^{\sigma -s-1}_{|\Delta|=\sigma
-2}\bigcup  F^{\sigma -s-2}_{|\Delta|=\sigma -2}}a_f C^f_g(\psi_1,\dots,\psi_s)
\nonumber\\
\qquad{}-{\rm div}_i
\sum_{h\in H} a_h C^{h,i}_{g}(\psi_1,\dots ,\psi_s)=
 \sum_{f\in F^{q>0,\delta>\mu}_{|\Delta|\le\sigma -3}} a_f
C^f_{g}(\psi_1,\dots ,\psi_s),\label{kaneis}
\end{gather}
where the right hand side stands for a generic linear combination
of acceptable complete contractions in the form \eqref{linisymric}
with $q>0$ factors $\nabla^{(p)}{\rm Ric}$ or $R$ and with  $\delta >\mu$ and
$|\Delta|\le \sigma -3$.
\end{lemma}

 We observe that if we can show the above,
we will then clearly have proven the
remaining case for our Lemma~\ref{postpone1st}. So, the rest of this
subsection is devoted to showing the above Lemma~\ref{reduce}.

\begin{proof}[Proof of Lemma \ref{reduce}.]
 We distinguish two cases and prove them separately. Either
$\sigma <\frac{n}{2}-1$ or $\sigma =\frac{n}{2}-1$. We begin with
the f\/irst case.

{\bf Proof of Lemma \ref{reduce} in the case $\boldsymbol{\sigma<\frac{n}{2}-1}$.}
The proof consists of three steps, which are spelled out along the proof.

  We observe that each complete contraction in the form
(\ref{linisymric}) with length $\sigma<\frac{n}{2}$ and weight~$-n$ and $|\Delta|=\sigma-1$ (i.e.\ a total of $\sigma -1$ factors $\Delta\psi_h$ or $R$) must
have the one factor that is {\it not} of the form $\Delta\psi_h$
or $R$  being in the form:
\begin{gather}
\label{larissa1} \nabla^{a_1\dots a_{\frac{n-2(\sigma
-1)}{2}}}\nabla_{a_1\dots a_{\frac{n-2(\sigma -1)}{2}}}\psi_h,
\end{gather}
 or in the form
\begin{gather}
\label{larissa2} \nabla^{a_1\dots a_{\frac{n-2(\sigma
-1)}{2}}}\nabla_{a_1\dots a_{\frac{n-2(\sigma
-1)}{2}-2}}{\rm Ric}_{a_{\frac{n-2(\sigma -1)}{2}-1}a_{\frac{n-2(\sigma
-1)}{2}}}.
\end{gather}
 Therefore all the complete contractions with $|\Delta|= \sigma -1$
must have $\delta\ge \mu+4$.

{\bf First step in the proof of Lemma \ref{reduce}.}
Let   $F_{\sigma-1}=F^{\sigma-s}_{|\Delta|=\sigma-1}\bigcup
F^{\sigma-s-1}_{|\Delta|=\sigma-1} \bigcup F^{\sigma-s-2}_{|\Delta|=\sigma-1}$.
We claim that there is a linear combination of
vector f\/ields in the form (\ref{linisymric}) with length $\sigma$,
say $\sum_{h\in H} a_h C^{h,i}_{g}(\psi_1,\dots ,\psi_s)$ so that
modulo complete contractions of length $\ge\sigma +1$:
\begin{gather}
 \sum_{f\in F_{\sigma -1}} a_f
C^f_{g}(\psi_1,\dots ,\psi_s) -{\rm div}_i \sum_{h\in H} a_h
C^{h,i}_{g}(\psi_1,\dots ,\psi_s)
\nonumber\\ \qquad{} =\sum_{u=1}^s ({\rm const})_u C^u_{g}(\psi_1,\dots ,\psi_s)+
\sum_{w=1}^s\sum_{q=1}^{w-1}({\rm const})_{q,w} C^{w,q}_{g}(\psi_1,\dots
,\psi_s)
\nonumber\\ \qquad{} +({\rm const})_{*} C^{*}_{g}(\psi_1,\dots ,\psi_s)
+\sum_{z\in Z} a_z C^z_{g}(\psi_1,\dots ,\psi_s),\label{twentysix}
\end{gather}
where $C^{*}_{g}(\psi_1,\dots ,\psi_s)$ stands for the complete
contraction:
\begin{gather*}
 {\rm contr}\Big(\nabla^{a_1\dots a_{\frac{n-2(\sigma
-1)}{2}-1}}\nabla_{a_1\dots a_{\frac{n-2(\sigma
-1)}{2}-3}}{\rm Ric}_{a_{\frac{n-2(\sigma -2)}{2}-1}a_{\frac{n-2(\sigma
-1)}{2}-1}}
\nonumber\\
\qquad{} \otimes \nabla^{ab}{\rm Ric}_{ab}\otimes R^{\sigma
-s-2}\otimes\Delta\psi_1\otimes\dots \otimes\Delta\psi_s\Big),
\end{gather*}
whereas each $C^u_{g}(\psi_1,\dots ,\psi_s)$ stands for the
complete contraction:
\begin{gather*}
{\rm contr}\Big(\nabla^{a_1\dots a_{\frac{n-2(\sigma
-1)}{2}-1}}\nabla_{a_1\dots a_{\frac{n-2(\sigma
-1)}{2}-3}}{\rm Ric}_{a_{\frac{n-2(\sigma -2)}{2}-1}a_{\frac{n-2(\sigma
-1)}{2}-1}}\otimes {\nabla^{ab}}_{ab}\psi_u
\nonumber\\
\qquad{} \otimes R^{\sigma
-s-1}\otimes\Delta\psi_1\otimes\dots\otimes
\hat{\Delta\psi_u}\otimes\Delta\psi_s\Big),
\end{gather*}
 and each $C^{w,q}_{g}(\psi_1,\dots ,\psi_s)$ stands for the complete contraction:
\begin{gather*}
 {\rm contr}\Big(\nabla^{a_1\dots a_{\frac{n-2(\sigma
-1)}{2}-1}}\nabla_{a_1\dots a_{\frac{n-2(\sigma
-1)}{2}-1}}\psi_w\otimes {\nabla^{ab}}_{ab}\psi_q
\nonumber\\
\qquad{} \otimes R^{\sigma
-s}\otimes\Delta\psi_1\otimes\dots\otimes
\hat{\Delta\psi_q}\otimes\dots\otimes\hat{\Delta\psi_w}\otimes\dots\otimes
\Delta\psi_s\Big).
\end{gather*}

Finally, $\sum_{z\in Z} a_z C^z_{g}(\psi_1,\dots ,\psi_s)$ stands
for a generic linear combination of complete contractions in the
form~(\ref{linisymric}) with $|\Delta|\le \sigma -3$ and
$\delta\ge\mu +2$. Moreover, we claim that each of the complete
contractions $C^{*}$, $C^{w,q}$, $C^u$ above has $\delta\ge\mu +4$.

{\bf Proof of (\ref{twentysix}).} We present the proof of this claim in detail, as the same
argument will be used repeatedly in many other instances. We have
already observed that we can write out:
\[
\sum_{f\in F_{\sigma -1}} a_f
C^f_{g}(\psi_1,\dots ,\psi_s)= ({\rm const})_1 C^I_{g}(\psi_1,\dots
,\psi_s)+\sum_{y=1}^s ({\rm const})_y C^{II,y}_{g}(\psi_1,\dots
,\psi_s).
\]
 Here $C^I$ is the complete contraction with one factor
 as in (\ref{larissa2})  and $\sigma -s-1$ factors $R$ and $s$ factors
$\Delta\psi_1,\dots ,\Delta\psi_s$. $C^{II,y}$ is the complete
contraction with one factor as in (\ref{larissa1}) and $\sigma -s$
factors $R$ and $s-1$ factors $\Delta\psi_1,\dots
,\hat{\Delta\psi_y},\dots \Delta\psi_s$. In particular, $C^I$
 is in the form:
\begin{gather}
 {\rm contr}\Big(\nabla^{a_1\dots
a_{\frac{n-2(\sigma -1)}{2}}}\nabla_{a_1\dots a_{\frac{n-2(\sigma
-1)}{2}-2}}{\rm Ric}_{a_{\frac{n-2(\sigma -1)}{2}-1}a_{\frac{n-2(\sigma
-1)}{2}}}
\nonumber\\
\qquad {}\otimes R^{\sigma
-s-1}\otimes\Delta\psi_1\otimes\dots\otimes\Delta\psi_s\Big).\label{italida}
\end{gather}

 Now, if we can show that we can f\/ind vector f\/ields
\[
\sum_{h\in H_1} a_h C^{h,i}_{g}(\psi_1,\dots ,\psi_s), \qquad
\sum_{h\in H_2} a_h C^{h,i}_{g}(\psi_1,\dots ,\psi_s),
\] so that:
\begin{gather}
  C^I_{g}(\psi_1,\dots ,\psi_s) -{\rm div}_i \sum_{h\in H_1} a_h
C^{h,i}_{g}(\psi_1,\dots ,\psi_s)
\nonumber\\
\qquad{}=\sum_{u=1}^s ({\rm const})_u C^u_{g}(\psi_1,\dots ,\psi_s)+
\sum_{w=1}^s\sum_{q=1}^{w-1}({\rm const})_{q,w} C^{w,q}_{g}(\psi_1,\dots
,\psi_s)
\nonumber\\ \qquad{}+ ({\rm const})_{*} C^{*}_{g}(\psi_1,\dots ,\psi_s)
+\sum_{z\in Z} a_z C^z_{g}(\psi_1,\dots ,\psi_s)\label{twentysixa}
\end{gather}
and
\begin{gather}
  C^{II,y}_{g}(\psi_1,\dots ,\psi_s) -{\rm div}_i \sum_{h\in H_2} a_h
C^{h,i}_{g}(\psi_1,\dots ,\psi_s)
\nonumber\\ \qquad{}= \sum_{u=1}^s ({\rm const})_u C^u_{g}(\psi_1,\dots ,\psi_s)+
\sum_{w=1}^s\sum_{q=1}^{w-1}({\rm const})_{q,w} C^{w,q}_{g}(\psi_1,\dots
,\psi_s)
\nonumber\\
\qquad{}+ ({\rm const})_{*} C^{*}_{g}(\psi_1,\dots ,\psi_s)
+\sum_{z\in Z} a_z C^z_{g}(\psi_1,\dots ,\psi_s),\label{twentysixb}
\end{gather}
then  clearly (\ref{twentysix}) will follow. We f\/irst prove
(\ref{twentysixa}).

{\bf Proof of (\ref{twentysixa}).} We def\/ine $C^{I,a_1}_{g}(\psi_1,\dots ,\psi_s)$ to stand for the
vector f\/ield that arises from~$C^I$  by erasing the index~${}^{a_1}$ in the f\/irst factor (see~(\ref{italida})) and making the index ${}_{a_1}$ that is then left into
a free index. We then calculate:
\begin{gather*}
C^{I}_{g}(\psi_1,\dots ,\psi_s)-{\rm div}_{a_1}
C^{I,a_1}_{g}(\psi_1,\dots
,\psi_s)\nonumber
\\
\qquad{} =-(\sigma-s-1)C^{I,A}_{g}(\psi_1,\dots
,\psi_s)-\sum_{h=1}^s C^{I,h}_{g}(\psi_1,\dots ,\psi_s),
\end{gather*}
where $C^{I,A}$ stands for the complete contraction that arises in
${\rm div}_{a_1} C^{I,a_1}_{g}(\psi_1,\dots ,\psi_s)$ when
$\nabla_{a_1}$ hits one of the factors $R$ and $C^{I,h}$ stands
for the complete contraction that arises when $\nabla_{a_1}$ hits
the factor~$\Delta\psi_h$. Given $C^{I,A}$, we then def\/ine
$C^{I,A,a_1}_{g}(\psi_1,\dots ,\psi_s)$ to stand for the vector
f\/ield that arises  by erasing the index ${}_{a_1}$ in its
f\/irst factor (note it is a derivative index) and making the index~${}_{a_1}$ in the factor $\nabla_{a_1}R$ into a free index. Also,
for each~$h$, $1\le h\le s$, we def\/ine $C^{I,h,a_1}_{g}$ to stand
for the vector f\/ield that arises from $C^{I,h}_g(\psi_1,\dots,\psi_s)$
 by erasing the index~${}_{a_1}$ in its f\/irst factor (note it is a derivative index) and
making the index~${}_{a_1}$ in the factor
$\nabla_{a_1}\Delta\psi_h$ into a free index. We then calculate:
\begin{gather}
\big[C^I_{g}(\psi_1,\dots ,\psi_s) -{\rm div}_{a_1}
C^{I,a_1}_{g}(\psi_1,\dots ,\psi_s)\big]
\nonumber\\
\qquad{}- (\sigma -s-1){\rm div}_{a_1}C^{I,A,a_1}_{g}(\psi_1,\dots
,\psi_s)-{\rm div}_{a_1}\sum_{h=1}^s C^{I,h,a_1}_{g}(\psi_1,\dots
,\psi_s)
\nonumber\\
\qquad{} =(\sigma -s-1) C^{*}_{g}(\psi_1,\dots ,\psi_s)
+\sum_{u=1}^s C^u_{g}(\psi_1,\dots ,\psi_s)+ \sum_{z\in Z}
a_z C^z_{g}(\psi_1,\dots ,\psi_s).\label{sobara}
\end{gather}
(We are using the notation of (\ref{twentysix}).) (\ref{sobara})  just follows
by the def\/initions. This proves (\ref{twentysixa}).

{\bf Proof of  (\ref{twentysixb}).} We denote by
$C^{II,y,*}_{g}(\psi_1,\dots ,\psi_s)$ the complete contraction:
\begin{gather*}
{\rm contr}\Big(\Delta^{\frac{n-2(\sigma
-1)}{2}-1}\psi_y\otimes\Delta R\otimes R^{\sigma
-s-1}\otimes\Delta\psi_1\otimes\dots\otimes
\hat{\Delta\psi_y}\otimes\dots\otimes\Delta\psi_s\Big).
\end{gather*}

Then, by complete analogy with the previous case, we
explicitly construct a divergence, ${\rm div}_i \sum_{h\in H} a_h
C^{h,i}_{g}(\psi_1,\dots ,\psi_s)$  from
$C^{II,y}_{g}(\psi_1,\dots ,\psi_s)$ so that:
\begin{gather} C^{II,y}_{g}(\psi_1,\dots
,\psi_s)- {\rm div}_i \sum_{h\in H} a_h C^{h,i}_{g}(\psi_1,\dots
,\psi_s)
\nonumber\\ \qquad{} = \sum_{y'\ne y} C^{y,y'}_{g}(\psi_1,\dots ,\psi_s)+
C^{II,y,*}_{g}(\psi_1,\dots ,\psi_s)+ \sum_{z\in Z} a_z
C^z_{g}(\psi_1,\dots ,\psi_s).\label{quickly}
\end{gather}

 Now, clearly $C^{II,y,*}_{g}(\psi_1,\dots ,\psi_s)$ is not
in the form $C^u$ or $C^{*}$ above. Also, for each of the complete
contractions $C^{y,y'}_{g}(\psi_1,\dots ,\psi_s)$, we inquire
whether $y>y'$ or $y<y'$. In the f\/irst case, we actually have a
complete contraction in the form $C^{q,w}$ that is allowed in  the
right hand side of~(\ref{twentysixb}). So in that case, we keep
the complete contraction~$C^{y,y'}$.

In the second case,
$C^{y,y'}$ is not a complete contraction in one of the forms on the right hand
side of~(\ref{twentysixb}).
 Now, by repeating the same argument as before,
it follows that we can explicitly construct
 a divergence ${\rm div}_i\sum_{h\in H} a_h
C^{h,i}_{g}(\psi_1,\dots ,\psi_s)$  such that  $y<y'$ and deduce that
modulo complete contractions of length $\ge\sigma +1$:
\begin{gather*}
 C^{y,y'}_{g}(\psi_1,\dots ,\psi_s)- {\rm div}_i\sum_{h\in H} a_h
C^{h,i}_{g}(\psi_1,\dots ,\psi_s)
\nonumber\\
\qquad{} =C^{y',y}_{g}(\psi_1,\dots ,\psi_s) +\sum_{z\in Z} a_z
C^z_{g}(\psi_1,\dots ,\psi_s),
\\
C^{II,y,*}_{g}(\psi_1,\dots ,\psi_s)- {\rm div}_i\sum_{h\in H} a_h
C^{h,i}_{g}(\psi_1,\dots ,\psi_s)
\nonumber\\ \qquad{} =C^{y,*}_{g}(\psi_1,\dots ,\psi_s)
+\sum_{z\in Z} a_z C^z_{g}(\psi_1,\dots ,\psi_s).
\end{gather*}
Thus, we have shown (\ref{twentysixb}) and therefore~(\ref{twentysix}).

{\bf A study of the complete contractions with $\boldsymbol{|\Delta|=\sigma-2}$ in (\ref{kaneis}).}
 Now, we focus on the complete contractions in the index sets
$F^q$, $q>0$ in (\ref{kaneis}) that have $|\Delta|=\sigma -2$.
We have observed that
only complete contractions in $F^{\sigma -s}$, $F^{\sigma -s-1}$
and $F^{\sigma -s-2}$ can have $|\Delta|=\sigma -2$. We have denoted
 the respective index sets by
$F^{\sigma -s}_{|\Delta|=\sigma -2}\subset F^{\sigma -s}$,
$F^{\sigma -s-1}_{|\Delta|=\sigma -2}\subset F^{\sigma -s-1}$
and $F^{\sigma -s-2}_{|\Delta|=\sigma -2} \subset F^{\sigma -s-2}$.
Clearly, since we are dealing with the case $s<\sigma -2$, each complete contraction $C^f$ indexed
 in one of the three sets above must have at least one factor $R$.
 Hence (by the ``Important remark'' in the introduction)
 it follows that each complete contraction $C^f$ above has
$\delta\ge\mu +2$. Therefore, by the def\/inition of
$F^{\sigma -s}_{|\Delta|=\sigma -2}$,  we derive
that we can write out:
\begin{gather*}
\sum_{f\in F^{\sigma -s}_{|\Delta|=\sigma -2}} a_f C^f_{g}
(\psi_1,\dots ,\psi_s)=\sum_{f\in F^{\sigma -s}_{|\Delta|= \sigma
-2,A}} a_f C^f_{g}(\psi_1,\dots ,\psi_s)
\nonumber\\ \qquad{} +\sum_{f\in F^{\sigma -s}_{|\Delta|=
\sigma -2,B}} a_f C^f_{g}(\psi_1,\dots ,\psi_s) + \sum_{f\in
F^{\sigma -s}_{|\Delta|= \sigma -2,C}} a_f C^f_{g}(\psi_1,\dots
,\psi_s).
\end{gather*}
Here $\sum_{f\in F^{\sigma -s}_{|\Delta|= \sigma -2,A}} a_f
C^f_{g}(\psi_1,\dots ,\psi_s)$ stands for a linear combination of
complete contractions in the form:
\begin{gather}
{\rm contr}\Big(\nabla^{a_1\dots a_t}\nabla^{(p)}_{r_1\dots r_p}{\rm Ric}_{ij}\otimes
\nabla^{b_1\dots b_s}\nabla^{(q)}_{z_1\dots z_q}{\rm Ric}_{i'j'}\otimes
R^{\sigma -s-2}\otimes\Delta\psi_1\otimes\dots\otimes\Delta\psi_s\Big),\label{catalyst1}
\end{gather}
with $\delta\ge\mu  +2$ (we observe that all the indices
${}_{r_1},\dots ,{}_{r_p},{}_i,{}_j$ in the f\/irst
 factor that are {\it not} involved in an internal contraction must
 contract against an index ${}_{z_1},\dots ,{}_{z_q},{}_{i'},{}_{j'}$ in the second factor
and vice versa).
 $\sum_{f\in F^{\sigma -s}_{|\Delta|=
\sigma -2,B}} a_f C^f_{g}(\psi_1,\dots ,\psi_s)$ stands for a
generic linear combination of complete contractions in the form:
\begin{gather}
 {\rm contr}\Big(\nabla^{a_1\dots a_t}\nabla^{(p)}_{r_1\dots
r_p}{\rm Ric}_{ij}\otimes \nabla^{b_1\dots b_y}\nabla^{(q)}_{z_1\dots
z_q}\psi_h\otimes R^{\sigma
-s-2}\nonumber\\
\qquad{} \otimes\Delta\psi_1\otimes\dots\otimes\hat{\Delta\psi_h}\otimes
\dots \otimes\Delta\psi_s\Big),\label{catalyst2}
\end{gather}
with $\delta\ge\mu  +2$  and where $y+q\ge 2$
 (but the second factor is not in the form $\Delta\psi_h$, by def\/inition)
 and where as above all the indices ${}_{r_1},\dots ,{}_{r_p},{}_i,{}_j$ in the f\/irst
 factor that are {\it not} involved in an internal contraction must
 contract against an index ${}_{z_1},\dots ,{}_{z_q}$ in the second factor
and vice versa.

 Finally, $\sum_{f\in F^{\sigma -s}_{|\Delta|=
\sigma -2,C}} a_f C^f_{g}(\psi_1,\dots ,\psi_s)$ stands for a
generic linear combination of complete contractions in the form:
\begin{gather}
 {\rm contr}\Big(\nabla^{a_1\dots a_t}\nabla^{(p)}_{r_1\dots
r_p}\psi_h\otimes \nabla^{b_1\dots b_y}\nabla^{(q)}_{z_1\dots
z_q}\psi_{h'}\otimes R^{\sigma
-s-2}
\nonumber \\ \qquad{} \otimes\Delta\psi_1\otimes\dots\otimes\hat{\Delta\psi_h}\otimes
\dots \otimes\hat{\Delta\psi_{h'}}\otimes\dots \otimes\Delta\psi_s\Big),\label{catalyst3}
\end{gather}
with $\delta\ge\mu  +2$ and where $t+p$, $y+q\ge 2$ (and neither of the f\/irst two factors
is of the form $\Delta\psi_b$)  and
 all the indices ${}_{r_1},\dots ,{}_{r_p}$
in the f\/irst
 factor that are {\it not} involved in an internal contraction must
 contract against an index ${}_{z_1},\dots ,{}_{z_q}$ in the second factor
and vice versa.

 Similarly, we see that we can write out:
\begin{gather*}
 \sum_{f\in F^{\sigma -s-1}_{|\Delta|=\sigma -2}} \!\! a_f C^f_{g}
(\psi_1,\dots ,\psi_s)=\sum_{f\in F^{\sigma -s-1}_{|\Delta|=
\sigma -2,\alpha}} \!\! a_f C^f_{g}(\psi_1,\dots ,\psi_s)
 +\sum_{f\in F^{\sigma -s-1}_{|\Delta|=
\sigma -2,\beta}} \!\! a_f C^f_{g}(\psi_1,\dots ,\psi_s),
\end{gather*}
where $\sum_{f\in F^{\sigma -s-1}_{|\Delta|= \sigma -2,\alpha}} a_f
C^f_{g}(\psi_1,\dots ,\psi_s)$ stands for a generic linear
combination of complete contractions in the form:
\begin{gather*}
 {\rm contr}\Big(\nabla^{a_1\dots a_t}\nabla^{(m)}_{r_1\dots r_m}R_{ijkl}\otimes
\nabla^{b_1\dots b_s}\nabla^{(q)}_{z_1\dots z_q}\psi_h\otimes
R^{\sigma
-s-1}\nonumber\\
\qquad{} \otimes\Delta\psi_1\otimes\dots\otimes\hat{\Delta\psi_h}\otimes
\dots \otimes\Delta\psi_s\Big),
\end{gather*}
with $\delta\ge\mu  +2$ and where $s+q\ge 2$ and where as above we
observe that all the indices ${}_{r_1},\dots
,{}_{r_m},{}_i,{}_j,{}_k,{}_l$ in the f\/irst
 factor that are {\it not} involved in an internal contraction must
 contract against an index ${}_{z_1},\dots ,{}_{z_q}$ in the second factor
and vice versa.

  Also, $\sum_{f\in F^{\sigma -s-1}_{|\Delta|=
\sigma -2,\beta}} a_f C^f_{g}(\psi_1,\dots ,\psi_s)$ stands for a
generic linear combination of complete contractions in the form:
\begin{gather}
{\rm contr}\Big(\nabla^{a_1\dots a_t}\nabla^{(m)}_{r_1\dots r_m}R_{ijkl}\otimes
\nabla^{b_1\dots b_s}\nabla^{(q)}_{z_1\dots z_q}{\rm Ric}_{i'j'}\otimes
R^{\sigma -s-2}\otimes\Delta\psi_1\otimes \dots
\otimes\Delta\psi_s\Big),\label{catalyst5}
\end{gather}
 with $\delta\ge\mu  +2$ and where as above we observe that
 all the indices ${}_{r_1},\dots ,{}_{r_p},{}_i,{}_j,{}_k,{}_l$ in
 the f\/irst factor that are {\it not} involved in an internal
 contraction must contract against an index ${}_{z_1},\dots ,{}_{z_q},{}_{i'},{}_{j'}$ in
 the second factor and vice versa.

  Finally, it follows that:
\begin{gather*}
\sum_{f\in F^{\sigma -s-2}_{|\Delta|=\sigma -2}} a_f
C^f_{g} (\psi_1,\dots ,\psi_s)=\sum_{f\in F^{\sigma
-s-2}_{|\Delta|= \sigma -2,\gamma}} a_f C^f_{g}(\psi_1,\dots
,\psi_s),
\end{gather*}
where $\sum_{f\in F^{\sigma -s-2}_{|\Delta|= \sigma -2,\gamma}}
a_f C^f_{g}(\psi_1,\dots ,\psi_s)$ stands for a linear combination
of complete contractions in the form:
\begin{gather}
{\rm contr}\Big(\nabla^{a_1\dots a_t}\nabla^{(p)}_{r_1\dots r_p}R_{ijkl}\otimes
\nabla^{b_1\dots b_s}\nabla^{(q)}_{z_1\dots z_q}R_{i'j'k'l'}
 \otimes
R^{\sigma -s-2}\otimes\Delta\psi_1\otimes\dots\otimes\Delta\psi_s\Big) \label{catalyst6}
\end{gather}
(all the indices ${}_{r_1},\dots
,{}_{r_p},{}_i,{}_j,{}_k,{}_l$ in the f\/irst
 factor that are {\it not} involved in an internal contraction must
 contract against an index ${}_{z_1},\dots ,{}_{z_q},{}_{i'},{}_{j'},{}_{k'},{}_{l'}$ in the second
factor and vice versa).

{\bf Second step in the proof of Lemma \ref{reduce}.}
The second step consists of proving the equations (\ref{salty}),
 (\ref{salty2}),
(\ref{salty3}), (\ref{salty4}), (\ref{salty5}), (\ref{saranta}),
(\ref{saranta2}) below. We f\/irstly claim that
we can subtract a divergence, ${\rm div}_i \sum_{h\in H} a_h
 C^{h,i}_{g}(\psi_1,\dots ,\psi_s)$, from the sublinear combination
 $\sum_{f\in F^{\sigma -s}_{|\Delta|=
\sigma -2,A}} a_f C^f_{g}(\psi_1,\dots ,\psi_s)$ so as to
 obtain an equation, modulo complete contractions of length
$\ge\sigma +1$:
\begin{gather}
  \sum_{f\in F^{\sigma -s}_{|\Delta|= \sigma -2,A}} a_f
C^f_{g}(\psi_1,\dots ,\psi_s)- {\rm div}_i \sum_{h\in H} a_h
C^{h,i}_{g}(\psi_1,\dots ,\psi_s)\label{salty}
\\ \qquad{}  =({\rm const})_1\cdot C^{*}_{g}(\psi_1,\dots ,\psi_s)+
\sum_{f\in F^{\sigma -s-1}_{|\Delta|= \sigma -2,\beta}} a_f
C^f_{g}(\psi_1,\dots ,\psi_s)  +\sum_{z\in Z'} a_z
C^z_{g}(\psi_1,\dots ,\psi_s),\nonumber
\end{gather}
where $C^{*}_{g}(\psi_1,\dots ,\psi_s)$ stands (as before) for the complete
 contraction:
\begin{gather*}
{\rm contr}\Big(\nabla^{a_1\dots
a_{\frac{n}{2}-4}ij}\nabla^{(\frac{n}{2}-4)}_{a_1\dots
a_{\frac{n}{2}-4}}{\rm Ric}_{ij}\otimes \nabla^{i'j'}{\rm Ric}_{i'j'}\otimes
R^{\sigma -s-2}\otimes\Delta\psi_1\otimes\dots\otimes\Delta\psi_s\Big),
\end{gather*}
 while $\sum_{f\in F^{\sigma -s-1}_{|\Delta|= \sigma -2,\beta}}
a_f C^f_{g}(\psi_1,\dots ,\psi_s)$ again stands for a generic
linear combination of complete contractions
 in the form (\ref{catalyst5}) and $\sum_{z\in Z'} a_z C^z_{g}(\psi_1,\dots ,\psi_s)$
stands for a generic linear combination of complete contractions in the
 form (\ref{linisymric}) with length $\sigma$, $\delta\ge \mu +1$ and
$|\Delta|=\sigma -3$.

  Next, consider the linear combination
$\sum_{f\in F^{\sigma -s}_{|\Delta|= \sigma -2,B}} a_f
C^f_{g}(\psi_1,\dots ,\psi_s)$. We then claim that there is a
divergence of a vector f\/ield, ${\rm div}_i \sum_{h\in H} a_h
C^{h,i}_{g}(\psi_1,\dots ,\psi_s)$, so that, modulo complete
 contractions of length $\ge\sigma +1$:
\begin{gather}
 \sum_{f\in F^{\sigma -s}_{|\Delta|= \sigma -2,B}} a_f
C^f_{g}(\psi_1,\dots ,\psi_s)- {\rm div}_i \sum_{h\in H} a_h
C^{h,i}_{g}(\psi_1,\dots ,\psi_s)\label{salty2}
\\  \qquad{} =\sum_{y=1}^s({\rm const})_y\cdot C^{y}_{g}(\psi_1,\dots ,\psi_s)+
\!\!\sum_{f\in F^{\sigma -s-1}_{|\Delta|= \sigma -2,\alpha}}\!\!\! a_f
C^f_{g}(\psi_1,\dots ,\psi_s)  +\sum_{z\in Z'} a_z
C^z_{g}(\psi_1,\dots ,\psi_s),\nonumber
\end{gather}
where the complete contractions $C^{y}_{g}(\psi_1,\dots ,\psi_s)$ are in
the form (\ref{catalyst2}) with every index involved in an
internal contraction, and the linear combination $\sum_{z\in Z'}
a_z C^z_{g}(\psi_1,\dots ,\psi_s)$ is the same as
in~(\ref{salty}).

Finally, we consider the linear combination
$\sum_{f\in F^{\sigma -s}_{|\Delta|= \sigma -2,C}} a_f
C^f_{g}(\psi_1,\dots ,\psi_s)$. Just by swit\-ching two factors,
 we assume that $h>h'$ in (\ref{catalyst3}).  We then claim that
we can construct a~divergence of a vector f\/ield, ${\rm div}_i \sum_{h\in H} a_h
C^{h,i}_{g}(\psi_1,\dots ,\psi_s)$, so that, modulo complete
 contractions of length $\ge\sigma +1$:
\begin{gather}
  \sum_{f\in F^{\sigma -s}_{|\Delta|= \sigma -2,C}} a_f
C^f_{g}(\psi_1,\dots ,\psi_s)- {\rm div}_i \sum_{h\in H} a_h
C^{h,i}_{g}(\psi_1,\dots ,\psi_s)
\nonumber\\ \qquad{} =\sum_{h=1}^s\sum_{h'=1}^{h-1}({\rm const})_{h,h'}\cdot
C^{h,h'}_{g}(\psi_1,\dots ,\psi_s) +\sum_{z\in Z'} a_z
C^z_{g}(\psi_1,\dots ,\psi_s),\label{salty3}
\end{gather}
where the complete contractions $C^{h,h'}_{g}(\psi_1,\dots ,\psi_s)$ are
in the form~(\ref{catalyst3}) with every index involved in an
internal contraction and the linear combination $\sum_{z\in Z'}
a_z C^z_{g}(\psi_1,\dots ,\psi_s)$ is the same as
 above.

{\bf More claims.} We now consider the generic linear combinations
 $\sum_{\!f{\in} F^{\sigma {-}s{-}1}_{|\Delta|{=} \sigma {-}2,\alpha}}\!\!\!\! a_f
C^f_{g}(\psi_1,\!{\dots},\!\psi_s),\!$ $\sum_{f\in F^{\sigma
-s-1}_{|\Delta|= \sigma -2,\beta}} a_f C^f_{g}(\psi_1,\dots
,\psi_s)$. In order to state our next claims, we def\/ine two
complete contractions: $C^{+}_{g}(\psi_1,\dots ,\psi_s)$ will
stand for the complete contraction:
\begin{gather*}
{\rm contr}\Big(\nabla^{a_1\dots
a_{\frac{n-2\sigma}{2}-1}il}\nabla_{a_1\dots
a_{\frac{n-2\sigma}{2}-1}}R_{ijkl}\otimes {\rm Ric}^{jk}\otimes
R^{\sigma-s-2}\otimes\Delta\psi_1\otimes\dots\otimes\Delta\psi_s\Big),
\end{gather*}
and $C^{++,h}_{g}(\psi_1,\dots ,\psi_s)$, $1\le h\le s$ will stand
for the complete contraction:
\begin{gather*}
{\rm contr}\Big(\nabla^{a_1\dots
a_{\frac{n-2\sigma}{2}-1}il}\nabla_{a_1\dots
a_{\frac{n-2\sigma}{2}-1}}R_{ijkl}\otimes
\nabla^{jk}\psi_h\otimes
R^{\sigma-s-1}\\
\qquad{} \otimes\Delta\psi_1\otimes\dots\otimes\hat{\Delta\psi_h}\otimes\dots
\otimes\Delta\psi_s\Big).
\end{gather*}

Moreover, we def\/ine $C^\sharp_{g}(\psi_1,\dots ,\psi_s)$ to
stand for the complete contraction:
\begin{gather*}
{\rm contr}\Big(\nabla^{a_1\dots
a_{\frac{n-2\sigma}{2}-3}il}\nabla_{a_1\dots
a_{\frac{n-2\sigma}{2}-3}}R_{ijkl}\otimes
\nabla^{ad}{{R_a}^{jk}}_d\otimes
R^{\sigma-s-1}\otimes\Delta\psi_1\otimes\dots \otimes\Delta\psi_s\Big).
\end{gather*}

{\bf Next claims.} Consider the linear combination
$\sum_{f\in F^{\sigma -s-1}_{|\Delta|=\sigma -2,\alpha}} a_f
C^f_{g}(\psi_1,\dots ,\psi_s)$. We claim that there is a
divergence of a vector f\/ield, ${\rm div}_i \sum_{h\in H} a_h
C^{h,i}_{g}(\psi_1,\dots ,\psi_s)$, so that, modulo complete
 contractions of length $\ge\sigma +1$:
\begin{gather}
  \sum_{f\in F^{\sigma -s-1}_{|\Delta|=\sigma -2,\alpha}} a_f
C^f_{g}(\psi_1,\dots ,\psi_s)- {\rm div}_i \sum_{h\in H} a_h
C^{h,i}_{g}(\psi_1,\dots ,\psi_s)
\nonumber\\
\qquad{} =  \sum_{y=1}^s({\rm const})'_y\cdot
C^{++,y}_{g}(\psi_1,\dots ,\psi_s) +\sum_{z\in Z'} a_z
C^z_{g}(\psi_1,\dots ,\psi_s),\label{salty4}
\end{gather}
and furthermore, for each $C^{++,y}$, we have $\delta
(C^{++,h})\ge \mu+2$.

On the other hand, we consider the linear combination
$\sum_{f\in F^{\sigma -s-1}_{|\Delta|=\sigma -2,\beta}} a_f
C^f_{g}(\psi_1,\dots ,\psi_s)$.  We then claim that there is a
divergence of a vector
 f\/ield, ${\rm div}_i \sum_{h\in H} a_h
C^{h,i}_{g}(\psi_1,\dots ,\psi_s)$, so that, modulo complete
 contractions of length $\ge\sigma +1$:
\begin{gather}
 \sum_{f\in F^{\sigma -s-1}_{|\Delta|=\sigma -2,\beta}} a_f
C^f_{g}(\psi_1,\dots ,\psi_s)- {\rm div}_i \sum_{h\in H} a_h
C^{h,i}_{g}(\psi_1,\dots ,\psi_s)\label{salty5}
\\
\qquad{} = \sum_{f\in F^{\sigma -s-1}_{|\Delta|=\sigma -2,\gamma}}
a_f C^f_{g}(\psi_1,\dots ,\psi_s)+ ({\rm const})''\cdot
C^{+}_{g}(\psi_1,\dots ,\psi_s)  +\sum_{z\in Z'} a_z
C^z_{g}(\psi_1,\dots ,\psi_s)\nonumber
\end{gather}
and furthermore for $C^{+}$, we have $\delta
(C^{+})\ge \mu+2$. (For the def\/inition of
$\sum_{f\in F^{\sigma -s-1}_{|\Delta|=\sigma -2,\gamma}}\cdots$, see~(\ref{catalyst6}).)

 Finally, we claim that there are divergences of vector
f\/ields,  ${\rm div}_i \sum_{h\in H} a_h C^{h,i}_{g}(\psi_1,\dots
,\psi_s)$  so that:
\begin{gather} C^{+}_{g}(\psi_1,\dots
,\psi_s)-{\rm div}_i \sum_{h\in H} a_h C^{h,i}_{g}(\psi_1,\dots
,\psi_s)\nonumber
\\
\qquad{}= C^\sharp_{g}(\psi_1,\dots ,\psi_s) +\sum_{z\in Z'} a_z
C^z_{g}(\psi_1,\dots ,\psi_s),\label{saranta}
\end{gather}
modulo complete contractions of length $\ge\sigma +1$ and also:{\samepage
\begin{gather}
 C^{++,h}_{g}(\psi_1,\dots
,\psi_s)-{\rm div}_i \sum_{h\in H} a_h C^{h,i}_{g}(\psi_1,\dots
,\psi_s)=\sum_{z\in Z'} a_z C^z_{g}(\psi_1,\dots ,\psi_s),\label{saranta2}
\end{gather}
modulo complete contractions of length $\ge\sigma +1$.}

We will now derive (\ref{salty}),
(\ref{salty2}), (\ref{salty3}), (\ref{salty4}), (\ref{salty5}),
(\ref{saranta}), (\ref{saranta2}). Before outlining the proof,
however, we note that once we show the above equations, we can
derive that there is a~linear combination of acceptable vector
f\/ields, $\sum_{h\in H} a_h C^{h,i}_{g}(\psi_1,\dots ,\psi_s)$, so
that:
\begin{gather}
 \sum_{f\in F^{\sigma }_{|\Delta|=\sigma-1}\bigcup F^{\sigma
-1}_{|\Delta|=\sigma -1}\bigcup F^{\sigma -2}_{|\Delta|=\sigma
-1}\bigcup F^{\sigma }_{|\Delta|=\sigma-2}\bigcup F^{\sigma
-1}_{|\Delta|=\sigma -2}\bigcup  F^{\sigma -2}_{|\Delta|=\sigma
-2}} a_f C^f_{g}(\psi_1,\dots ,\psi_s)
\nonumber\\
\qquad{} -{\rm div}_i \sum_{h\in H}
a_h C^{h,i}_{g}(\psi_1,\dots ,\psi_s)= ({\rm const})_{*}\cdot
C^{*}_{g}(\psi_1,\dots ,\psi_s)
\nonumber\\
\qquad{} +\sum_{y=1}^s
({\rm const})_{y,*}\cdot C^{y,*}_{g}(\psi_1,\dots ,\psi_s)+
\sum_{w=1}^s\sum_{q=1}^{w-1}({\rm const})_{w,q}\cdot
C^{w,q}_{g}(\psi_1,\dots ,\psi_s)
\nonumber\\
\qquad{}+ ({\rm const})_\sharp C^\sharp_{g}(\psi_1,\dots ,\psi_s) +\sum_{z\in
Z'}a_z C^z_{g}(\psi_1,\dots ,\psi_s).\label{portolos}
\end{gather}

{\bf Proof of (\ref{salty}), (\ref{salty2}),
(\ref{salty3}), (\ref{salty4}), (\ref{salty5}), (\ref{saranta}),
(\ref{saranta2}).}
The divergences needed in all the above equations are
constructed ``by hand'', by consecutively picking derivative indices,
 making them into free indices and taking the corresponding divergence.
As the proof is essentially the
same for the equations (\ref{salty}), (\ref{salty2}),
(\ref{salty3}), (\ref{salty4}), (\ref{salty5}), we will only
demonstrate the f\/irst one. Afterwards, we show (\ref{saranta}) and
(\ref{saranta2}).

  We consider a complete contraction $C^f_{g}(\psi_1,\dots
,\psi_s)$ in the form (\ref{catalyst1}). We f\/irst show that modulo
 introducing complete contractions of length $\ge\sigma +1$, we can
 subtract a divergence ${\rm div}_i \sum_{h\in H} a_h C^{h,i}_{g}(\psi_1,
\dots ,\psi_s)$ from $C^f_{g}(\psi_1,\dots ,\psi_s)$ and obtain:
\begin{gather}
  C^f_{g}(\psi_1,\dots ,\psi_s)-{\rm div}_i \sum_{h\in H} a_h
C^{h,i}_{g}(\psi_1,\dots ,\psi_s)= \sum_{v=2}^{\frac{n}{2}-\sigma
-4}({\rm const})_{f,v}\cdot C^{v,\tau}_{g}(\psi_1,\dots ,\psi_s)
\nonumber\\
\qquad{} + \sum_{f\in F^{\sigma
-s}_{|\Delta|= \sigma -2,\beta}} a_f C^f_{g}(\psi_1,\dots ,\psi_s)
+\sum_{z\in Z'} a_z C^z_{g}(\psi_1,\dots ,\psi_s),\label{kourasteis}
\end{gather}
where $C^{v,\tau}_{g}(\psi_1,\dots ,\psi_s)$ stands for a complete
contraction in the form (\ref{catalyst1}) where and all the
indices ${}_{r_1},\dots ,{}_{r_p},{}_i,{}_j$ and ${}_{z_1},\dots
,{}_{z_q},{}_{i'},{}_{j'}$ are contracting against a derivative
index (such a complete contraction has $\delta\ge\mu +2$) and also
the f\/irst factor has $t=v$.

 We observe that once we show the above, we can then repeat
the proof of equation (\ref{quickly}) to this setting to f\/ind a
divergence ${\rm div}_i \sum_{h\in H} a_h C^{h,i}_{g}(\psi_1,\dots
,\psi_s)$ so that for each $v$:
\begin{gather}
  C^{v,\tau}_{g}(\psi_1,\dots ,\psi_s)-{\rm div}_i \sum_{h\in H} a_h
C^{h,i}_{g}(\psi_1,\dots ,\psi_s)
\nonumber\\
\qquad{} = C^{*}_{g}(\psi_1,\dots ,\psi_s)+ \sum_{z\in Z'} a_z C^z_{g}(\psi_1,\dots ,\psi_s),\label{kourasteis2}
\end{gather}
modulo complete contractions of length $\ge\sigma +1$.

  Then, combining equations (\ref{kourasteis}) and
(\ref{kourasteis2}), we will deduce (\ref{salty}).

 We show (\ref{kourasteis}) by an induction. Firstly let
us observe that the maximum value that $\delta$ can have for a
complete contraction in the form (\ref{catalyst1}) is
$\frac{n}{2}+(\sigma-s)$. If $\delta=\frac{n}{2}+(\sigma-s)$ then
all indices appearing in (\ref{catalyst1}) must be involved in an
internal contraction (and thus there is nothing to show). Now, let us suppose that
$C^f_{g}(\psi_1,\dots ,\psi_s)$ has $\delta=P\ge\mu +2$,\footnote{This inequality holds
by virtue of the ``Important remark'' in the introduction, and because of
the fact that $s< \sigma-2$.}
$P<\frac{n}{2}+(\sigma-s)$. We then construct a vector f\/ield
$C^{h,i}_{g}(\psi_1,\dots ,\psi_s)$ so that modulo complete
contractions of length $\ge\sigma +1$:
\begin{gather}
  C^f_{g}(\psi_1,\dots ,\psi_s)-{\rm div}_i C^{h,i}_{g}(\psi_1,\dots
,\psi_s)
\nonumber\\
\qquad{} = (C^f)'_{g}(\psi_1,\dots ,\psi_s)
+ \sum_{f\in F^{\sigma
-s}_{|\Delta|= \sigma -2,\beta}} a_f C^f_{g}(\psi_1,\dots ,\psi_s)
+\sum_{z\in Z'} a_z C^z_{g}(\psi_1,\dots ,\psi_s).\label{kourasteisb}
\end{gather}
Here $(C^f)'_{g}(\psi_1,\dots ,\psi_s)$ is a complete contraction
in the form (\ref{catalyst1}) with $\delta =P+1$. Clearly, if we
can show the above then by iterative repetition (\ref{kourasteis}) will follow.

  The vector f\/ield $C^{h,i}_{g}(\psi_1,\dots
,\psi_s)$ needed for (\ref{kourasteisb}) can be easily
constructed. We observe that since $\sigma<\frac{n}{2}-1$ and the
weight is $-n$, we have that at least one of the f\/irst two factors
in (\ref{catalyst1}) has at least four indices. With no loss of
generality, we assume it is the f\/irst factor, $\nabla^{a_1\dots
a_t}\nabla^{(p)}_{r_1\dots r_p}{\rm Ric}_{ij}$ for which $t+p+2\ge 4$.
Moreover, since we are assuming $\delta<\frac{n}{2}+(\sigma-s)$,
it follows that at least one of the indices ${}_{r_1},\dots
,{}_{r_p},{}_i,{}_j$ is {\it not} contracting against an index
${}_{a_1},\dots ,{}_{a_t}$. That index can either be a derivative
index (say ${}_{r_1}$ with no loss of generality) or--if all the indices
${}_{r_1},\dots ,{}_{r_p}$ are each contracting against one of the
indices ${}_{a_1},\dots ,{}_{a_t}$~-- an internal index (say ${}_i$
with no loss of generality).

  In the f\/irst case, we def\/ine $C^{h,i}_{g}(\psi_1,\dots ,\psi_s)$
to be the vector f\/ield obtained from $C^f$ be erasing the
derivative index ${}_{r_1}$ in the f\/irst factor and making the
index ${}^{r_1}$ in the second factor in (\ref{catalyst1}) into a
free index. We check that for this vector f\/ield, (\ref{kourasteisb})
indeed holds.

  In the second case, we see that since $t+p+2\ge 4$, we can
apply the second Bianchi identity (modulo introducing a complete
contraction that will belong to the linear combination $\sum_{f\in
F^{\sigma -s}_{|\Delta|= \sigma -2,\beta}} a_f
C^f_{g}(\psi_1,\dots ,\psi_s)$) and be reduced to the previous
case, where the index ${}_{r_1}$ is not contracting against one of
the indices ${}_{a_1},\dots ,{}_{a_t}$. Thus, we have proven
(\ref{salty}). The equations~(\ref{salty2}), (\ref{salty3}),
(\ref{salty4}), (\ref{salty5}) follow by essentially the same
reasoning.

 We now show (\ref{saranta}), (\ref{saranta2}). In both cases,
we def\/ine $C^{b,l}_{g}(\psi_1,\dots ,\psi_s)$ to stand for the
vector f\/ield obtained from $C^{+}_{g}(\psi_1,\dots ,\psi_s)$ and
$C^{++,h}_{g}(\psi_1,\dots ,\psi_s)$ respectively, by erasing the
derivative index ${}^i$ in the f\/irst factor. Then, for
(\ref{saranta2}) we def\/ine
\[
\sum_{h\in H} a_h C^{h,i}_{g}(\psi_1,\dots ,\psi_s)
= C^{b,i}_{g}(\psi_1,\dots ,\psi_s).
\]
 Since $C^{++,h}$ has
$\delta\ge \mu +2$, we observe that if $\nabla_i$ in ${\rm div}_i
C^{b,i}_{g}(\psi_1,\dots ,\psi_s)$ hits the f\/irst factor, we
cancel out the complete contraction $C^{++,h}$. When $\nabla_i$
hits the second factor, we get a~complete contraction that is
equal to a complete contraction of length $\sigma +1$ (this is due
to the antisymmetry of the indices~${}_i$,~${}_j$). If it hits one of the
other factors, we get a complete contraction that belongs to the
linear combination $\sum_{z\in Z'} a_z C^z_{g}(\psi_1,\dots
,\psi_s)$.

 The case of (\ref{saranta}) is more complicated. We again
consider the same vector f\/ield $C^{b,i}_{g}(\psi_1$, $\dots ,\psi_s)$
as above. We then observe that modulo complete contractions of
length $\ge\sigma +1$:
\begin{gather}
\label{dance} C^{+}_{g}(\psi_1,\dots ,\psi_s)-{\rm div}_i
C^{b,i}_{g}(\psi_1,\dots ,\psi_s)=C'_{g}(\psi_1,\dots
,\psi_s)+\sum_{z\in Z'} a_z C^z_{g}(\psi_1,\dots ,\psi_s),
\end{gather}
where $C'$ is the complete contraction:
\[
{\rm contr}\big(\Delta^{\frac{n}{2}-1-\sigma}\nabla^lR_{ijkl}\otimes\nabla^{l'}{R^{ijk}}_{l'}\otimes
R^{\sigma
-s-2}\otimes\Delta\psi_1\otimes\dots\otimes\Delta\psi_s\big).
\]

 Notice that since $\sigma
<\frac{n}{2}-1$, we have $t>1$ derivatives on the f\/irst factor. We
then apply the second Bianchi identity and write:
\begin{gather}
 C'_{g}(\psi_1,\dots ,\psi_s)=
{\rm contr}\big(\Delta^{\frac{n}{2}-1-\sigma-1}\nabla^{itl}R_{tjkl}\otimes
\nabla^{l'}{{R_i}^{jk}}_{l'}\nonumber\\
\phantom{C'_{g}(\psi_1,\dots ,\psi_s)=}{} \otimes R^{\sigma-s-2}
\otimes\Delta\psi_1\otimes\dots\otimes\Delta\psi_s\big).\label{ring}
\end{gather}

 We then def\/ine $C^{c,i}_{g}(\psi_1,\dots ,\psi_s)$
 to stand for the vector f\/ield that arises from the right hand
 side of the
 above by erasing the index ${}^i$ in the f\/irst factor and making the
 index~${}_i$ that it contracted against in the second factor into a free
 index. We observe that:
 \begin{gather}
 \label{johncoates}
C'_{g}(\psi_1,\dots ,\psi_s)-{\rm div}_i C^{c,i}_{g}(\psi_1,\dots
,\psi_s)=C^\sharp_{g}(\psi_1,\dots ,\psi_s)+\sum_{z\in Z'} a_z
C^z_{g}(\psi_1,\dots ,\psi_s).
 \end{gather}

 Therefore, combining (\ref{dance}), (\ref{ring}) and (\ref{johncoates})
 we have that the vector f\/ield needed for (\ref{saranta}) is
 precisely:
\[
\sum_{h\in H} a_h C^{h,i}_{g}(\psi_1,\dots ,\psi_s)=
C^{b,i}_{g}(\psi_1,\dots ,\psi_s)-C^{c,i}_{g}(\psi_1,\dots
,\psi_s).
\]

{\bf The third step of the proof of Lemma \ref{reduce}.}
 In view of equation (\ref{portolos}),
by subtracting the divergence ${\rm div}_i \sum_{h\in H} a_h
C^{h,i}_{g}(\psi_1,\dots ,\psi_s)$ from $I^s_{g}(\psi_1,\dots
,\psi_s)$, we have obtained a relation, modulo complete
contractions of length $\ge\sigma +1$:
\begin{gather}
 I^s_{g}(\psi_1,\dots ,\psi_s)- {\rm div}_i \sum_{h\in H} a_h
C^{h,i}_{g}(\psi_1,\dots ,\psi_s)
 = \sum_{l\in L_\mu} a_l
C^{l,\iota}_g(\psi_1,\dots,\psi_s)\nonumber\\
\qquad{} +\sum_{j\in J} a_j
C^j_g(\psi_1,\dots,\psi_s)+
 \sum_{q_1=1}^{\sigma -s-3} \sum_{f\in F^{q_1}} a_f
C^f_{g}(\psi_1,\dots ,\psi_s)
  + ({\rm const})_{*}\cdot
C^{*}_{g}(\psi_1,\dots ,\psi_s)
\nonumber\\
\qquad{} +\sum_{y=1}^s
({\rm const})_{y,*}\cdot C^{y,*}_{g}(\psi_1,\dots ,\psi_s)+
\sum_{w=1}^s\sum_{q=1}^{w-1}({\rm const})_{w,q}\cdot
C^{w,q}_{g}(\psi_1,\dots ,\psi_s)
\nonumber\\
\qquad{} +({\rm const})_\sharp C^\sharp_{g}(\psi_1,\dots ,\psi_s) +\sum_{z\in
Z'}a_z C^z_{g}(\psi_1,\dots ,\psi_s).\label{portolos2}
\end{gather}

We next claim that $({\rm const})_{*}, ({\rm const})_{y,*},
({\rm const})_{w,q}, ({\rm const})_\sharp=0$; if we can prove this, we will then have proven Lemma
\ref{reduce} in this case where $s<\sigma-2$.

We f\/irst prove that $({\rm const})_{*}=0$. We denote the r.h.s.\ of the
above by $Z_g(\psi_1,\dots,\psi_s)$. Clearly,
$\int_{M^n} Z_g (\psi_1,\dots,\psi_s)dV_g=0$.

We now apply the ``main conclusion'' of the super divergence formula
to this integral equation (see \cite{a:dgciI}), deriving an equation
${\rm supdiv}[I^s_{g}(\psi_1,\dots
,\psi_s,\Omega^{\sigma -s})]=0$.\footnote{Recall that
 $I^s_{g}(\psi_1,\dots ,\psi_s,\Omega^{\sigma -s})$ stands
for the linear combination that arises from
$I^s_{g}(\psi_1,\dots ,\psi_s)$ by formally replacing each factor
$\nabla^{(p)}_{r_1\dots r_p}{\rm Ric}_{ij}\ne R$ by
$-\nabla^{(p+2)}_{r_1\dots r_pij}\Omega$ and each factor $R$ by
$-2\Delta\Omega$.} We focus on
the sublinear combination ${\rm supdiv}_{*}[I^s_{g}(\psi_1,\dots
,\psi_s,\Omega^{\sigma -s})]$ of complete contractions with
$\delta=0$, $s$ factors $\nabla\psi_h$ and $\sigma -s-2$ factors
$\nabla\Omega$. Since the ``main conclusion'' of the super divergence formula holds formally,
it follows that:
\begin{gather}
\label{astro1} {\rm supdiv}_{*}[Z_{g}(\psi_1,\dots
,\psi_s,\Omega^{\sigma -s})]=0,
\end{gather}
 modulo complete contractions of
length $\ge\sigma +1$.

Since each $C^f$ in (\ref{portolos2}) has $|\Delta|\le \sigma -3$,
it follows that:
\begin{gather*}
{\rm supdiv}_{*}[Z_{g}(\psi_1,\dots ,\psi_s,\Omega^{\sigma
-s})]=(-1)^{\frac{n}{2}}({\rm const})_{*}\cdot C^x_{g}(\psi_1,\dots
,\psi_s),
\end{gather*}
where $C^x_{g}(\psi_1,\dots ,\psi_s)$ is the expression:
\begin{gather}
{\rm contr}\Big(\nabla^{y_1\dots y_{\sigma -s-2}w_1\dots
w_s} \big(\nabla^{x_1\dots x_{\frac{n}{2}-\sigma
+2}}\Omega\otimes\nabla_{x_1\dots x_{\frac{n}{2}-\sigma
+2}}\Omega\big)\nonumber\\
\qquad{} \otimes\nabla_{y_1}\Omega\otimes\dots\otimes
\nabla_{y_{\sigma
-s-2}}\Omega\otimes\nabla_{w_1}\psi_1\otimes\dots\otimes\nabla_{w_s}\psi_s\Big).\label{eidiko1}
\end{gather}

Since this complete contraction is clearly not zero, we deduce that $({\rm const})_{*}=0$.

  Showing that the other constants in (\ref{portolos2}) are zero follows essentially
the same pattern. We next show that each $({\rm const})_{y,*}=0$. In
order to do that, we again consider the main conclusion of
the super divergence formula applied to
$Z_{g}(\psi_1,\dots ,\psi_s)$, ${\rm supdiv}[Z_{g}(\psi_1,\dots ,\psi_s, \Omega^{\sigma -s})]=0$.
 We pick out the sublinear combination  ${\rm supdiv}_{**,y}[Z_{g}(\psi_1,\dots ,\psi_s,
\Omega^{\sigma -s})]$ of complete contractions with length
$\sigma$, $\delta=0$ and $s-1$ factors $\nabla\psi_h$, $h=1,\dots
,\hat{y},\dots ,s$ and $\sigma -2-(s-1)$ factors $\nabla\Omega$.
Since the main conclusion of the super divergence holds formally, we deduce that:
\begin{gather}
\label{astro2}
{\rm supdiv}_{**,y}[I^s_{g}(\psi_1,\dots
,\psi_s,\Omega^{\sigma -s})]=0,
\end{gather}
 modulo complete contractions of
length $\ge\sigma +1$.

Analogously to the previous case, we deduce that:
\begin{gather*}
{\rm supdiv}_{**,y}[I^s_{g}(\psi_1,\dots
,\psi_s,\Omega^{\sigma -s})]=(-1)^{\frac{n}{2}}({\rm const})_{y,*}\cdot
C^{xx}_{g}(\psi_1,\dots ,\psi_s),
\end{gather*}
where $C^{xx}_{g}(\psi_1,\dots ,\psi_s)$ is the complete
contraction:
\begin{gather}
{\rm contr}\Big(\nabla^{y_1\dots y_{\sigma
-s-1}w_1\dots\hat{w_y}\dots w_s} \big(\nabla^{x_1\dots
x_{\frac{n}{2}-\sigma +2}}\Omega\otimes\nabla_{x_1\dots
x_{\frac{n}{2}-\sigma +2}}\psi_y\big)
\\
\qquad {} \otimes \nabla_{y_1}\Omega\otimes\dots\otimes \nabla_{y_{\sigma
-s-1}}\Omega\otimes\nabla_{w_1}\psi_1\otimes\dots\otimes\hat{\nabla\psi_y}
\otimes\dots\otimes\nabla_{w_s}\psi_s\Big).\label{eidiko2}
\end{gather}
 Therefore, we deduce that each $({\rm const})_{y,*}=0$.

{\sloppy To prove that $({\rm const})_{w,q}=0$, we again consider
$Z_{g}(\psi_1,\dots ,\psi_s, \Omega^{\sigma -s})$,\footnote{See
the ``Main consequence''
of the super divergence formula in Subsection~2.2.3 in~\cite{alexakis}.} for which
$\int_{M^n}Z_{g}(\psi_1,\dots ,\psi_s, \Omega^{\sigma -s})=0$. We apply the
super divergence formula to this equation (see~\cite{a:dgciI}), deriving a
new local equation,
${\rm supdiv}[Z_{g}(\psi_1,\dots ,\psi_s, \Omega^{\sigma -s})]=0$; we pick out
 the sublinear combination
${\rm supdiv}_{***,(w,q)}[I^s_{g}(\psi_1,\dots ,\psi_s, \Omega^{\sigma
-s})]$ of complete contractions with length $\sigma$, $\delta=0$
and $s-2$ factors $\nabla\psi_h$, $h=1,\dots ,\hat{w},\dots,
\hat{q},s$ and $\sigma -s$ factors $\nabla\Omega$. Since the super
divergence formula holds formally, we deduce that:
\begin{gather}
\label{astro3}
 {\rm supdiv}_{***,(w,q)}\big[I^s_{g}(\psi_1,\dots
,\psi_s,\Omega^{\sigma -s})\big]=0,
\end{gather}
modulo complete
contractions of length $\ge\sigma +1$.

}

 We now claim that:
\begin{gather}
\label{aram3} {\rm supdiv}_{***,(w,q)}\big[I^s_{g}(\psi_1,\dots
,\psi_s,\Omega^{\sigma
-s})\big]=(-1)^{\frac{n}{2}}({\rm const})_{\sharp}\cdot
C^{xxx}_{g}(\psi_1,\dots ,\psi_s),
\end{gather}
where $C^{xxx}_{g}(\psi_1,\dots ,\psi_s)$ is the complete
contraction:
\begin{gather}
   {\rm contr}\Big(\nabla^{y_1\dots y_{\sigma
-s-2}z_1\dots\hat{z_w}\dots\hat{z_q}\dots  z_s} \big(\nabla^{x_1\dots
x_{\frac{n}{2}-\sigma +2}}\psi_w\otimes\nabla_{x_1\dots
x_{\frac{n}{2}-\sigma +2}}\psi_q\big)
\nonumber\\
\qquad{} \otimes\nabla_{y_1}\Omega\otimes\dots\otimes
\nabla_{y_{\sigma
-s}}\Omega\otimes\nabla_{w_1}\psi_1\otimes\dots\otimes \hat{\nabla\psi_w}
\otimes  \dots\otimes \hat{\nabla\psi_q}\otimes\dots\otimes\nabla_{w_s}\psi_s\Big).\label{eidiko3}
\end{gather}

 In order to see this, we recall that we have already shown that
$({\rm const})_{*}$,
 $({\rm const})_{II,y,*}, ({\rm const})_{w,q}$ $=0$, and also each
other complete contraction $C^f$ in (\ref{portolos2}) has
$|\Delta|\le\sigma -3$, hence each complete contraction of length
$\sigma$ in each ${\rm Tail}[C^f]$ can have at most $\sigma -3$ factors
$\nabla\Omega$ or $\nabla\psi_h$.  Therefore, we deduce that
$({\rm const})_{(w,q)}=0$.

 Finally, to show $({\rm const})_\sharp=0$, we
again consider $Z_{g}(\psi_1,\dots ,\psi_s,
\Omega^{\sigma -s-2})$,\footnote{See the ``Main consequence''
of the super divergence formula in Subsection~2.2.3 in~\cite{alexakis}.}
for which we have $\int_{M^n}Z_{g}(\psi_1,\dots ,\psi_s,
\Omega^{\sigma -s-2})dV_g=0$. We apply the super divergence formula to this equation,
deriving  a new local equation:  ${\rm supdiv}[I^s_{g}(\psi_1,\dots
,\psi_s, \Omega^{\sigma -s-2})]=0$; we pick out the sublinear
combination ${\rm supdiv}_{****}[Z_{g}(\psi_1,\dots ,\psi_s,
\Omega^{\sigma -s-2})]$ of complete contractions with length
$\sigma$, $\delta=0$ and $s$ factors $\nabla\psi_h$,  and $\sigma
-s-2$ factors $\nabla\Omega$. Since the super divergence formula
holds formally, we deduce that:
\begin{gather*}
{\rm supdiv}_{****}\big[Z_{g}(\psi_1,\dots
,\psi_s,\Omega^{\sigma -s-2})\big]=0,
\end{gather*}
modulo complete contractions of length $\ge\sigma +1$.

 Analogously to the previous case, we deduce that:
\begin{gather}
\label{aram4}
{\rm supdiv}_{****}\big[I^s_{g}(\psi_1,\dots
,\psi_s,\Omega^{\sigma
-s-2})\big]=(-1)^{\frac{n}{2}}({\rm const})_\sharp\cdot
C^{xxxx}_{g}(\psi_1,\dots ,\psi_s),
\end{gather}
where $C^{xxxx}_{g}(\psi_1,\dots ,\psi_s)$ is the complete
contraction:
\begin{gather*}
 {\rm contr}\Big(\nabla^{y_1\dots y_{\sigma -s-2}z_1\dots  z_s}
\big(\nabla^{x_1\dots x_{\frac{n}{2}-\sigma
-2}i'l'}{{R^i}_{jk}}^l\otimes\nabla_{x_1\dots
x_{\frac{n}{2}-\sigma -2}il}{{R_{i'}}^{jk}}_{l'}\big)
\\ \qquad{} \otimes\nabla_{y_1}\Omega\otimes\dots\otimes
\nabla_{y_{\sigma
-s-2}}\Omega\otimes\nabla_{w_1}\psi_1\otimes\dots\otimes\nabla_{w_s}\psi_s\Big).
\end{gather*}

 Hence $({\rm const})_\sharp =0$. We have shown Lemma \ref{reduce} in the case
$\sigma <\frac{n}{2}-1$.

{\bf Proof of Lemma \ref{reduce} in the case
$\boldsymbol{\sigma=\frac{n}{2}-1}$.}
 The case $\sigma =\frac{n}{2}-1$ follows
similarly. In this case, we can apply the same method
of explicitly constructing divergences to show that
there is a linear combination of acceptable vector f\/ields,
$\sum_{h\in H} a_h C^{h,i}_{g}(\psi_1,\dots ,\psi_s)$, so that
(\ref{portolos}) holds,  where in this setting
$C^{*}_{g}(\psi_1,\dots ,\psi_s)$ is in the form:
\begin{gather*}
{\rm contr}\big(\nabla^a{\rm Ric}_{ai}\otimes\nabla^b{\rm Ric}_b^i\otimes R^{\sigma
-s-2}\otimes\Delta\psi_1\otimes\dots\otimes\Delta\psi_s\big),
\end{gather*}
whereas $C^{y,*}_{g}(\psi_1,\dots ,\psi_s)$ is in the form:
\begin{gather*}
{\rm contr}\big(\nabla^a{\rm Ric}_{ai}\otimes\nabla^{bi}_b\psi_y\otimes R^{\sigma
-s-2}\otimes\Delta\psi_1\otimes\dots\otimes \hat{\Delta\psi_y}\otimes\dots
\otimes\Delta\psi_s\big),
\end{gather*}
and $C^{q,w}_{g}(\psi_1,\dots ,\psi_s)$ is in the form:
\begin{gather*}
{\rm contr}\big(\nabla^{ai}_a\psi_q\otimes\nabla^{bi}_b\psi_w\otimes
R^{\sigma
-s-2}\otimes\Delta\psi_1\otimes\dots\otimes \hat{\Delta\psi_q}\otimes\dots\otimes
\hat{\Delta\psi_w}\otimes \dots \otimes\Delta\psi_s\big),
\end{gather*}
 and f\/inally $C^\sharp_{g}(\psi_1,\dots ,\psi_s)$ is
in the form:
\begin{gather*}
{\rm contr}\big(\nabla^aR_{ijka}\otimes\nabla^b{R^{ijk}}_b\otimes R^{\sigma
-s-2}\otimes\Delta\psi_1\otimes\dots\otimes\Delta\psi_s\big).
\end{gather*}

Therefore, by subtracting the divergence ${\rm div}_i \sum_{h\in H} a_h
C^{h,i}_{g}(\psi_1,\dots ,\psi_s)$ from $I^s_{g}(\psi_1,\dots
,\psi_s)$, we again obtain (\ref{portolos2}) (modulo complete
contractions of length $\ge\sigma +1$), with the new notational
conventions of this setting.

Hence, if we could show that $({\rm const})_{*}, ({\rm const})_{y,*},
({\rm const})_{w,q}, ({\rm const})_\sharp=0$, we will then have proven Lemma
\ref{reduce} in the case $\sigma =\frac{n}{2}-1$.

  Now, we show that
$({\rm const})_{*}, ({\rm const})_{II,y}, ({\rm const})_{w,q}=0$. As before, we take
${\rm supdiv}[Z_{g}(\psi_1,\dots ,$ $\psi_s,\Omega^{\sigma -s-2})]$ and
 focus on the same sublinear combinations ${\rm supdiv}_{*}, {\rm supdiv}_{**}, {\rm supdiv}_{***}$
 as in the previous case, and then observe
that the same equations (\ref{astro1}), (\ref{astro2}),
(\ref{astro3}) also hold in this case. The complete contractions
$C^x$, $C^{xx}$, $C^{xxx}$ are the same as in equations
(\ref{eidiko1}), (\ref{eidiko2}), (\ref{eidiko3}), where we set
$\sigma =\frac{n}{2}-1$.

  Now, in order to show that $({\rm const})_\sharp =0$, we apply
 the same method as for the previous case: We consider and
focus on the sublinear combination of complete contractions with
length~$\sigma$, $\delta=0$, and $s$ factors $\nabla\psi_h$,
$\sigma -2-s$ factors $\nabla\Omega$. We denote that sublinear
combination by ${\rm supdiv}_{****}[I^s_{g}(\psi_1,\dots ,\psi_s,
\Omega^{\sigma -s-2})]$. By the same arguments as before, we again
deduce an equation~(\ref{aram4}), where~$C^{xxxx}$ here stands for
the complete contraction:
\begin{gather*}
  {\rm contr}\Big(\nabla^{y_1\dots y_{\sigma -s-2}z_1\dots  z_s}
\big(\nabla^{l'}{R_{ijk}}^l\otimes\nabla_l{R^{ijk}}_{l'}\big)
\\ \qquad{} \otimes\nabla_{y_1}\Omega\otimes\dots\otimes
\nabla_{y_{\sigma
-s-2}}\Omega\otimes\nabla_{w_1}\psi_1\otimes\dots\otimes\nabla_{w_s}\psi_s\Big).
\end{gather*}

We thus deduce that $({\rm const})_\sharp =0$ in this case
 also.
 \end{proof}

\subsection[The proof of Lemmas \ref{postpone1st}, \ref{violi}, \ref{postpone2}
and Lemma \ref{tistexnes} when $s=\sigma-2$]{The proof of Lemmas \ref{postpone1st}, \ref{violi}, \ref{postpone2}
and Lemma \ref{tistexnes} when $\boldsymbol{s=\sigma-2}$}

 We observe that if $\sigma <\frac{n}{2}-1$ then
$\mu\le \frac{n}{2}$. We distinguish subcases, based
on $\sigma$ and $\mu$.
 Our three subcases are  $(\sigma<\frac{n}{2}-1,\mu<\frac{n}{2})$,
  $(\sigma<\frac{n}{2}-1,\mu=\frac{n}{2})$,
and  $\sigma=\frac{n}{2}-1$.

  Firstly, we consider the subcase where $\sigma <\frac{n}{2}-1$ {\it
and} $\mu<\frac{n}{2}$.

\begin{proof}[Proof of Lemma \ref{tistexnes} when $\boldsymbol{s=\sigma-2}$, $\boldsymbol{\sigma<\frac{n}{2}-1}$.]
Consider $P(g)|_{\Theta_{\sigma -2}}$.\footnote{Recall that
$\Theta_{\sigma-2}$ stands for the index set of complete
contractions in $P(g)$ (which are in the form ${\rm contr}(\nabla^{(m)}W\otimes\dots\otimes\nabla^{(m')}W\otimes\nabla^{(a)}P\otimes\dots\otimes\nabla^{(a')}P$)
 with length $\sigma$ and with $\sigma-2$ factors $\nabla^{(u)}P$.} We
pick out the sublinear combination of complete contractions in the form:
\begin{gather*}
{\rm contr}\big(\nabla^{a_1\dots a_t}\nabla^{(m_1)}W_{ijkl}\otimes
\nabla^{b_1\dots b_s}\nabla^{(m_2)}W_{i'j'k'l'}\otimes
(P^a_a)^{\sigma -2}\big).
\end{gather*}

We index those complete contractions in the set
$\Theta^{*}_{\sigma -2}$. We then easily see (by repeating the
explicit constructions from the previous subsection)
that we can subtract
a divergence from $P(g)|_{\Theta^{*}_{\sigma -2}}$ so that, modulo
complete contractions of length $\ge\sigma +1$:
\begin{gather*}
P(g)|_{\Theta^{*}_{\sigma -2}}-{\rm div}_i \sum_{h\in H}
a_h C^{h,i}(g) =({\rm const})_\alpha C^\alpha(g) +\sum_{t\in T} a_t
C^t(g),
\end{gather*}
where  $C^\alpha(g)$ in the complete contraction in the form:
\begin{gather}
\label{goutsi}
{\rm contr}\big(\Delta^{\frac{n}{2}-\sigma -2}
\nabla^{il}W_{ijkl}\otimes\nabla^{i'l'}{{W_{i'}}^{jk}}_{l'}\otimes
(P^a_a)^{\sigma -2}\big),
\end{gather}
while\footnote{If $\sigma=\frac{n}{2}-1$ then $C^\alpha(g)$ is in
the form ${\rm contr}(\nabla^lW_{ijkl}\otimes\nabla_{l'}W^{ijkl'}\otimes
(P^a_a)^{\sigma-2})$.} $\sum_{t\in T} a_t C^t(g)$ stands for a
linear combination of
 complete contractions in the form:
\begin{gather*}
{\rm contr}\big(\nabla^{a_1\dots a_t}\nabla^{(m_1)}W_{ijkl}\otimes
\nabla^{b_1\dots b_s}\nabla^{(m_2)}W_{i'j'k'l'}\otimes \nabla^f
P_{fi}\otimes(P^a_a)^{\sigma -3}\big)
\end{gather*}
(with fewer than $\sigma -2$ factors $P^a_a$ and with $\delta\ge \mu$).

  Therefore, from now on we may assume with no loss of generality that the
sublinear combination $P(g)|_{\Theta^{*}_{\sigma -2}}$ is
 precisely $({\rm const})_\alpha C^\alpha(g)$.
\end{proof}

\begin{proof}[Proof of Lemmas \ref{postpone1st}, \ref{violi}, \ref{postpone2} when $\boldsymbol{s=\sigma-2}$,
$\boldsymbol{\sigma<\frac{n}{2}-1}$, $\boldsymbol{\mu<\frac{n}{2}}$.]{\sloppy
 We again observe that  the complete contractions in
$I^s_{g}(\psi_1,\dots ,\psi_s)$ with $|\Delta|\ge \sigma -2$ will
be
 indexed in the sets~$F^2$,~$F^1$
and\footnote{Recall that $F^q\subset F$ stands for the index set
of complete contractions with $q$ factors $\nabla^{(p)}{\rm Ric}$ or~$R$.} in this case, also in the sets $J$
 and $L$, where for each~$C^l$, $l\in L$ and each~$C^j$, $j\in J$ we
 recall
 that~$C^l$,~$C^j$ are in the form~(\ref{linisymric}) with no factors
$\nabla^{(p)}{\rm Ric}$.

}

 We now again consider $I^s_{g}(\psi_1,\dots ,\psi_s)$, written
out as a linear combination in the form (\ref{polis}).

  We recall the def\/inition of the index set $F^{*}$ (see the notation above Lemma~\ref{postpone1st}). It follows that
  $F^{*}= F_{\sigma-1}^2\bigcup F_{\sigma-2}^2\bigcup F_{\sigma-2}^1$
  (recall that the upper labels count the number of factors
  $\nabla^{(p)}{\rm Ric}$ or $R$ and the lower labels count the value of
  $|\Delta|$). Furthermore, since $P(g)|_{\Theta_s}$ is ``good'',
  it follows that $L_\mu^{*}=\varnothing$, while $\sum_{j\in J^{*}} a_j C^j_g(\psi_1,\dots,\psi_s)=
({\rm const}) C_g(\psi_1,\dots,\psi_s)$, where $C_g(\psi_1,\dots,\psi_s)$ is
  the complete contraction ${\rm contr}(\Delta^{\frac{n}{2}-\sigma-2} \nabla^{il}R_{ijkl}\otimes
  \nabla_{i'l'}R^{i'jkl'}\otimes\Delta\psi_1\otimes\dots\otimes\Delta\psi_s)$
  or ${\rm contr}( \nabla^{i}R_{ijkl}\otimes
  \nabla_{i'}R^{i'jkl}\otimes\Delta\psi_1\otimes\dots\otimes\Delta\psi_s)$.

 {\bf A study of the sublinear combination
 $\boldsymbol{\sum_{f\in F^2_{\sigma -1}}
a_f C^f_{g}(\psi_1,\dots ,\psi_s)}$.} As before, we see that since
$|\Delta|=\sigma -1$, $s=\sigma-2$ for each $C^f$, $f\in F^2_{\sigma
-1}$,\footnote{Recall that the subscript ${}_{\sigma-1}$ means that
$|\Delta|=\sigma-1$
 for the complete contractions indexed in $F^2_{\sigma -1}$.}
it follows that each~$C^f$, $f\in F^2_{\sigma -1}$ must have at least
one factor $R$ and will hence have $\delta(C^f)\ge\mu +2$ (by the
decomposition of the Weyl tensor).

  As before, it follows that:
\begin{gather*}
\sum_{f\in F^2_{\sigma -1}} a_f C^f_{g}(\psi_1,\dots
,\psi_s)=({\rm const})_{@} C^{@}_{g}(\psi_1,\dots , \psi_s)  +
\sum_{u=1}^{\sigma -2} ({\rm const})_u C^u_{g}(\psi_1,\dots ,\psi_s),
\end{gather*}
where here $C^{@}_g$ is the complete contraction:
\begin{gather*}
{\rm contr}\big(\Delta^{(\frac{n}{2}-\sigma)}R\otimes R\otimes\Delta\psi_1\otimes
\dots\otimes\Delta\psi_{\sigma -2}\big),
\end{gather*}
and $C^u_{g}(\psi_1,\dots ,\psi_s)$ is the complete contraction:
\begin{gather*}
 {\rm contr}\big(\Delta^{(\frac{n}{2}-\sigma+1)}\psi_u\otimes
R^2\otimes \Delta\psi_1\otimes\hat{\Delta\psi_u}\otimes\dots
 \otimes\Delta\psi_{\sigma -2}\big).
\end{gather*}

As before, we denote by $\sum_{z\in Z'} a_z C^z_{g}(\psi_1,\dots
,\psi_s)$ a generic linear combination of complete contractions in
the form
 (\ref{linisymric}) with length $\sigma$, $|\Delta|\le\sigma -3$ and
$\delta\ge \mu+1$.
 It follows as in the previous subsection that we can construct a vector f\/ield $\sum_{h\in H}
a_h C^{h,i}_{g}(\psi_1,\dots ,\psi_s)$ so that,  modulo complete
 contractions of length $\ge\sigma +1$:
\begin{gather*}
 \sum_{f\in F^2_{\sigma -1}} a_f C^f_{g}(\psi_1,\dots ,\psi_s)
-{\rm div}_i \sum_{h\in H} a_h C^{h,i}_{g}(\psi_1,\dots ,\psi_s)
\\
\qquad{}= ({\rm const})_{!} C^{!}_{g}(\psi_1,\dots ,
\psi_s)+ \sum_{u=1}^{\sigma -2} ({\rm const})_{!,u}
C^{!,u}_{g}(\psi_1,\dots ,\psi_s)
\\
\qquad{} + \sum_{q=1}^{\sigma -2}\sum_{w=1}^{q-1}
({\rm const})_{!,(q,w)} C^{!,(q,w)}_{g}(\psi_1,\dots ,\psi_s) +
\sum_{z\in Z'} a_z C^z_{g}(\psi_1,\dots ,\psi_s),
\end{gather*}
where $C^{!}_{g}(\psi_1,\dots , \psi_s)$ is the complete
contraction:
\begin{gather*}
{\rm contr}\big(\Delta^{(\frac{n}{2}-\sigma-1)}R\otimes \Delta R\otimes
\Delta\psi_1\otimes
\dots\otimes\Delta\psi_{\sigma -2}\big),
\end{gather*}
$C^{!,u}_{g}(\psi_1,\dots ,\psi_s)$ is the complete contraction:
\begin{gather*}
{\rm contr}\big(\Delta^{(\frac{n}{2}-\sigma-1)}R\otimes
\Delta^2 \psi_u\otimes
R\otimes\Delta\psi_1\otimes\hat{\Delta\psi_u}\otimes
\dots\otimes\Delta\psi_{\sigma -2}\big),
\end{gather*}
$C^{!,(q,w)}_{g}(\psi_1,\dots ,\psi_s)$ is the complete
 contraction:
\begin{gather*}
{\rm contr}\big(\Delta^{(\frac{n}{2}-\sigma)}\psi_q\otimes
R^2\otimes\Delta^2 \psi_w\otimes
\Delta\psi_1\otimes\hat{\Delta\psi_w}\otimes\dots \otimes \hat{\Delta\psi_q}
\otimes \dots\otimes\Delta\psi_{\sigma -2}\big).
\end{gather*}
This follows by the same explicit constructions as for equation (\ref{twentysix}).

{\bf A study of the sublinear combination
 $\boldsymbol{\sum_{f\in F^2_{\sigma -2}}
a_f C^f_{g}(\psi_1,\dots ,\psi_s)}$.} We claim that
$\sum_{f\in F^2_{\sigma -2}}
a_f C^f_{g}(\psi_1,\dots ,\psi_s)$ can be expressed as:
\begin{gather}
 \sum_{f\in F^2_{\sigma -2}} a_f C^f_{g}(\psi_1,\dots ,\psi_s)=
\sum_{v\in V_1} a_v C^v_{g}(\psi_1,\dots ,\psi_s) +
 \sum_{v\in V_2} a_v C^v_{g}(\psi_1,\dots ,\psi_s)\nonumber\\
 \qquad{} +
\sum_{v\in V_3} a_v C^v_{g}(\psi_1,\dots ,\psi_s)+
({\rm const})_{\rm spec}\cdot C^{\rm spec}_{g}(\psi_1,\dots ,\psi_s),\label{papariz}
\end{gather}
where $\sum_{v\in V_1} a_v C^v_{g}(\psi_1,\dots ,\psi_s)$ stands
 for a linear combination of complete contractions in the form:
\begin{gather}
{\rm contr}\big(\nabla^{a_1\dots a_t}\nabla^{(p)}_{r_1\dots r_p}\psi_q\otimes
\nabla^{b_1\dots b_y}\nabla^{(q)}_{z_1\dots z_q}\psi_w\otimes R^2\nonumber\\
\qquad{} \otimes
\Delta\psi_1\otimes\dots\otimes \hat{\Delta\psi_q}\otimes \dots\otimes \hat{\Delta\psi_w}\otimes \dots
\otimes\Delta\psi_s\big),\label{egw}
\end{gather}
with $\delta\ge \mu +2$. $\sum_{v\in V_2} a_v C^v_{g}(\psi_1,\dots ,\psi_s)$ is a linear
combination of complete contractions in the form:
\begin{gather}
\label{egw2}
{\rm contr}\big(\nabla^{a_1\dots a_t}\nabla^{(p)}_{r_1\dots
r_p}\psi_u\otimes \nabla^{b_1\dots b_y}\nabla^{(q)}_{z_1\dots
z_q}{\rm Ric}_{ij}\otimes R\otimes
\Delta\psi_1\otimes\dots\otimes\hat{\Delta\psi_u}\otimes\dots
\otimes\Delta\psi_s\big),\!\!\!
\end{gather}
with $\delta\ge \mu +2$. $\sum_{v\in V_3} a_v C^v_{g}(\psi_1,\dots ,\psi_s)$
is a linear combination of complete contractions in the form:
\begin{gather}
\label{egw3}
{\rm contr}\big(\nabla^{a_1\dots a_t}\nabla^{(p)}_{r_1\dots r_p}{\rm Ric}_{ij}\otimes
\nabla^{b_1\dots b_y}\nabla^{(q)}_{z_1\dots z_q}{\rm Ric}_{i'j'}\otimes
\Delta\psi_1\otimes\dots\otimes\Delta\psi_s\big),
\end{gather}
with $\delta\ge \mu +2$. Finally,
$C^{\rm spec}_{g}(\psi_1,\dots,\psi_s)$
 stands for a complete contraction in the form:
\begin{gather*}
{\rm contr}\big(\Delta^{(\frac{n}{2}-\sigma-2)}\nabla_{ij}R\otimes\nabla^{ij}
R\otimes\Delta\psi_1\otimes\dots\otimes\Delta\psi_s\big),
\end{gather*}
where $\delta\ge \mu +1$.

(\ref{papariz}) follows by the def\/inition of $\sum_{f\in F} \cdots$,
 apart for the claims regarding the numbers of internal contractions in the dif\/ferent factors.
 We check these claims by virtue of our assumptions
on~$I^s_{g}$. Firstly, polarizing the function~$\psi$, we write
out~$I^s_{g}$ as a  linear combination of contractions in the
form:
\begin{gather} {\rm contr}
\big(\nabla^{a_1\dots
a_t}\nabla^{(m)}W_{ijkl}\otimes\dots\otimes\nabla^{b_1\dots
b_u}\nabla^{(m_{\sigma-s})}W_{i'j'k'l'}\nonumber\\
\qquad{} \otimes\nabla^{v_1\dots
v_x}\nabla^{(p_1+2)}\psi\otimes\dots\otimes\nabla^{y_1\dots
y_w}\nabla^{(p_s+2)}\psi\big).\label{onlyweyl}
\end{gather}
Now, using the decomposition of the Weyl tensor we decompose the
above complete contractions. As noted in the introduction, we observe that each complete
contraction with a factor~$R$ must have~$\delta\ge\mu +2$, This
shows our claim for~(\ref{egw}) and~(\ref{egw2}). Moreover, we
recall that all complete contractions $C^l (g)$, $l\in
\Theta_{\sigma -2}$ with $s$ factors $P^a_a$ are in the
 form $C^\alpha$, as in (\ref{goutsi}). Hence, by applying the
 decomposition of the Weyl tensor to the
complete contraction $C^\alpha$, we see that the contribution of
$C^\alpha$ to the linear combination $\sum_{f\in F^2_{\sigma -2}}
a_f C^f$ is precisely the complete contraction~$C^{\rm spec}$ (times a
 constant). For each other complete contraction $C^l$,
$l\in \Theta_{\sigma -2}\setminus \Theta^{*}_{\sigma -2}$, we have
that $C^l_{g}(\psi_1,\dots , \psi_s)$, written in the form~(\ref{onlyweyl}) does not have~$s$
 factors~$\Delta\psi_h$ (by def\/inition). Hence, each complete contraction $C^f$,
$f\in F^2_{\sigma -2}$, in the decomposition of each
$C^l_{g}(\psi_1,\dots ,\psi_s)$, $l\in \Theta_{\sigma -2}\setminus
\Theta^{*}_{\sigma -2}$ must have $\delta\ge\mu +2$. This shows~(\ref{egw3}).

  In view of (\ref{papariz}), and by applying the same ``by hand''
 technique as for equation (\ref{twentysix}), we see that we can f\/ind a
 vector f\/ield
$\sum_{h\in H} a_h C^{h,i}_{g}(\psi_1,\dots ,\psi_s)$ so that
modulo complete contractions of length $\ge\sigma +1$:
\begin{gather}
 \sum_{f\in F^2_{\sigma -2}} a_f C^f_{g}(\psi_1,\dots ,\psi_s)
-{\rm div}_i \sum_{h\in H} a_h C^{h,i}_{g}(\psi_1,\dots ,\psi_s)= ({\rm const})_{!} C^{!}_{g}(\psi_1,\dots ,
\psi_s)
\nonumber\\ \qquad{}  + \sum_{u=1}^{\sigma -2} ({\rm const})_{!,u}
C^{!,u}_{g}(\psi_1,\dots ,\psi_s)
  + \sum_{q=1}^{\sigma -2}\sum_{w=1}^{q-1}
({\rm const})_{!,(q,w)} C^{!,(q,w)}_{g}(\psi_1,\dots ,\psi_s)
\nonumber\\
\qquad{} + \sum_{p\in P^{2,\delta\ge\mu +2}_{\sigma-2}} a_p C^p_{g}
(\psi_1,\dots ,\psi_s)+ \sum_{z\in Z'} a_z C^z_{g}(\psi_1,\dots
,\psi_s),\label{papariz2}
\end{gather}
where $\sum_{p\in P^{2,\delta\ge\mu +2}_{\sigma -2}} a_p C^p_{g}
(\psi_1,\dots ,\psi_s)$ stands for a linear combination of
complete contractions in the form (\ref{linisymric}) with length
$\sigma$, $s=\sigma -2$, $|\Delta|=\sigma -2$ and $q=1$ and also
$\delta\ge\mu +2$. The ``by hand'' construction of the vector f\/ield
for (\ref{papariz2}) is the same as the proof of equation
(\ref{twentysix}), with the slight caveat that we now also have
the
 complete contraction $C^{\rm spec}$.  In particular, since we only know
 that $\delta\ge\mu +1$ for $C^{\rm spec}$, we have to check that we get no
 correction terms of the form~(\ref{linisymric}) with $q=1$. But this
follows by the form of $C^{\rm spec}$, i.e.\ that we have two factors
$\nabla^{(t)}R$, so we never have to introduce correction terms by
virtue of
 the second Bianchi identity.

{\bf A study of the sublinear combination
 $\boldsymbol{\sum_{f\in F^1_{\sigma -2}}
a_f C^f_{g}(\psi_1,\dots ,\psi_s)}$.}
We claim that this sublinear combination will be of the form:
\begin{gather}
\sum_{f\in F^1_{\sigma -2}} a_f C^f_{g}(\psi_1,\dots ,\psi_s)=
\sum_{d\in D} a_d C^d_{g}(\psi_1,\dots ,\psi_s)
\nonumber\\
\qquad{} +({\rm const})_\beta C^\beta_{g}(\psi_1,\dots ,\psi_s)+({\rm const})_\gamma
C^\gamma_{g}(\psi_1,\dots ,\psi_s),\label{ga}
\end{gather}
where $\sum_{d\in D} a_d C^d_{g}(\psi_1,\dots ,\psi_s)$ stands for
a generic linear combination of complete contractions in the form~(\ref{linisymric}) with length $\sigma$, $s=\sigma -2$,
$|\Delta|=\sigma -2$ and $q=1$ and also $\delta\ge\mu +2$ (just
like the linear combination $\sum_{p\in P^{2,\delta\ge\mu
+2}_{\sigma -2}} a_p C^p_{g} (\psi_1,\dots ,\psi_s)$ in~(\ref{papariz2})), while $C^\beta_{g}(\psi_1,\dots ,\psi_s)$ is
the complete contraction:
\begin{gather}
\label{ga2}
{\rm contr}\big(\Delta^{(\frac{n}{2}-\sigma -2)}
\nabla^{il}R_{ijkl}\otimes\nabla^{jk}R\otimes
\Delta\psi_1\otimes\dots\otimes\Delta\psi_s\big),
\end{gather}
with $\delta=\frac{n}{2}\ge\mu+1$, and $C^\gamma_{g}(\psi_1,\dots
,\psi_s)$ is the complete contraction:
\begin{gather}
\label{ga3} {\rm contr}\big(\Delta^{(\frac{n}{2}-\sigma -2)}
\nabla^{jk}R\otimes\nabla^{il}R_{ijkl}\otimes
\Delta\psi_1\otimes\dots\otimes\Delta\psi_s\big),
\end{gather}
with $\delta=\frac{n}{2}\ge\mu+1$.

  Observe that the sublinear combination
\[
({\rm const})_\beta
C^\beta_{g}(\psi_1,\dots ,\psi_s)+({\rm const})_\gamma
C^\gamma_{g}(\psi_1,\dots ,\psi_s)
\] above is exclusively the
contribution of the complete contraction $C^\alpha$ in $I^s$ to
the linear combination $\sum_{f\in F^1_{\sigma -2}} a_f
C^f_{g}(\psi_1,\dots ,\psi_s)$. (The fact that such contractions in
the forms (\ref{ga2}), (\ref{ga3}) do arise in the decomposition
of $C^\alpha$ can be directly checked.)

  Now, we check that the contribution of all the other
 complete contractions in $I^s$ (i.e.\ other than $C^\alpha$) to the
 sublinear combination
$\sum_{f\in F^1_{\sigma -2}} a_f C^f_{g}(\psi_1,\dots ,\psi_s)$ is
indeed $\sum_{d\in D} a_d C^d_{g}(\psi_1$, $\dots ,\psi_s)$, as
described in~(\ref{ga}).

 But this is straightforward to observe.
We only need to check that each $C^l_{g}(\psi_1, \dots ,\psi_s)$
in the form~(\ref{onlyweyl}), with $\delta\ge \mu$
 and with less than~$s$ factors $\Delta\psi_s$ can be written out as:
\begin{gather}
\label{paula} C^l_{g}(\psi_1,\dots ,\psi_s)=\sum_{d\in D} a_d
C^d_{g}(\psi_1, \dots ,\psi_s)+\sum_{q=0}^2\sum_{|\Delta|\le\sigma
-3} \sum_{f\in F^q_{|\Delta|}} a_f C^f_{g}(\psi_1,\dots ,\psi_s),
\end{gather}
where $\sum_{d\in D} a_d C^d_{g}(\psi_1,\dots ,\psi_s)$ is as
 above. (\ref{paula})~just follows by the decomposition of the Weyl tensor.

 Now, we seek to ``get rid'' of the complete contractions $C^\beta$,
$C^\gamma$. We construct vector f\/ields $C^{\beta, i}$, $C^{\gamma, i}$, where
 $C^{\beta,i}$ arises from $C^\beta$ by erasing the index~${}^j$
in the second factor and making the index~${}_j$  that it contracted
 against in the f\/irst factor into a free index~$i$.
 $C^{\gamma,i}$~arises from~$C^\gamma$ by erasing the index~${}^j$
in the f\/irst factor and making the index~${}_j$  that it contracted
 against in the second factor into a free index~${}_i$.

 We then observe that modulo complete contractions of length
$\ge\sigma +1$:
\begin{gather*}
C^\beta_{g}(\psi_1,\dots,\psi_s)-{\rm div}_i C^{\beta,
i}_{g}(\psi_1, \dots ,\psi_s)=\sum_{z\in Z'} a_z
C^z_{g}(\psi_1,\dots ,\psi_s),
\\
 C^\gamma_{g}(\psi_1,\dots,\psi_s)-{\rm div}_i C^{\gamma,
i}_{g}(\psi_1, \dots ,\psi_s)=\sum_{z\in Z'} a_z
C^z_{g}(\psi_1,\dots ,\psi_s).
\end{gather*}

 Next, we claim that we can subtract a divergence
$\sum_{h\in H} a_h C^{h,i}_{g}(\psi_1,\dots ,\psi_s)$ from the
linear combination $\sum_{d\in D} a_d C^d_{g}(\psi_1,\dots
,\psi_s)$ so that:
\begin{gather*}
 \sum_{d\in D} a_d
C^d_{g}(\psi_1,\dots ,\psi_s)- {\rm div}_i \sum_{h\in H} a_h
C^{h,i}_{g}(\psi_1,\dots ,\psi_s)=
({\rm const})_{?}\cdot C^{?}_{g}(\psi_1,\dots
,\psi_s)\\
\qquad{} + \sum_{u=1}^s({\rm const})_{??,u}\cdot
C^{??,u}_{g}(\psi_1,\dots ,\psi_s)+\sum_{z\in Z'} a_z
C^z_{g}(\psi_1,\dots ,\psi_s),
\end{gather*}
where $C^{?}_{g}(\psi_1,\dots ,\psi_s)$ is the complete
contraction:
\begin{gather*}
{\rm contr}\big(\Delta^{(\frac{n}{2}-\sigma-1)}\nabla^{il}R_{ijkl}\otimes
{\rm Ric}^{jk}\otimes\Delta\psi_1 \otimes\dots\otimes\Delta\psi_s\big),
\end{gather*}
 and it has $\delta\ge\mu +2$, while $C^{??,u}_{g}(\psi_1,\dots
 ,\psi_s)$ is the complete
contraction:
\begin{gather*}
{\rm contr}\big(\Delta^{(\frac{n}{2}-\sigma-1)}\nabla^{il}R_{ijkl}\otimes
\nabla^{jk}\psi_u\otimes\Delta\psi_1
\otimes\dots\otimes \hat{\Delta\psi_u}\otimes \dots\otimes\Delta\psi_s\big),
\end{gather*}
 and it has $\delta\ge\mu +2$. This follows by imitating
the proof of the previous case. Then, we explicitly construct a vector
f\/ield $\sum_{h\in H} a_h C^{h,i}_{g}(\psi_1,\dots ,\psi_s)$ so
that:
\begin{gather*} ({\rm const})_{?}\cdot
C^{?}_{g}(\psi_1,\dots ,\psi_s)+ \sum_{u=1}^s({\rm const})_{??,u}\cdot
C^{??,u}_{g}(\psi_1,\dots ,\psi_s)\\
\qquad{} -{\rm div}_i \sum_{h\in H} a_h
C^{h,i}_{g}(\psi_1,\dots ,\psi_s)
=({\rm const})_\sharp \cdot C^\sharp_{g}(\psi_1,\dots
,\psi_s) +\sum_{z\in Z'} a_z C^z_{g}(\psi_1,\dots ,\psi_s),
\end{gather*}
where $C^\sharp_{g}(\psi_1,\dots ,\psi_s)$ stands for the complete
contraction:
\begin{gather*}
{\rm contr}\big(\Delta^{(\frac{n}{2}-\sigma -2)}
\nabla^{il}R_{ijkl}\otimes\nabla^{i'l'}{{R_{i'}}^{jk}}_{l'}\otimes
\Delta\psi_1\otimes\dots\otimes\Delta\psi_s\big),
\end{gather*}
 which has $\delta\ge \mu +1$.

{\sloppy   Now, recall that by our lemma's assumption $L_\mu^{*}=\varnothing$, and
 $\sum_{j\in J^{*}}a_j
C^j_g(\psi_1,\dots,\psi_s)=({\rm const})''C^\sharp_g(\psi_1,\dots,\psi_s)$.
  Therefore, by our study of the three sublinear combinations
$\sum_{f\in F^1_{\sigma-2}}\cdots$,
$\sum_{f\in F^2_{\sigma-2}}\cdots$, $\sum_{f\in F^2_{\sigma-1}}\cdots$
we have shown that we can subtract a divergence
from~$I^s$ so that:
\begin{gather*}
 I^s_{g}(\psi_1,\dots ,\psi_s)- {\rm div}_i \sum_{h\in H} a_h
C^{h,i}_{g}(\psi_1,\dots ,\psi_s)=
  \sum_{l\in L_\mu} a_l C^l_{g}(\psi_1,\dots ,\psi_s)\\
\qquad{} +
\!\sum_{j\in J\setminus J^{*}}\! a_j C^j_{g}(\psi_1,\dots ,\psi_s)
+\sum_{q_1=0}^2 \sum_{v=0}^{\sigma -3}\! \sum_{f\in F^{q_1}_v}\! a_f
C^f_{g}(\psi_1,\dots ,\psi_s) + ({\rm const})_{!}\cdot
C^{!}_{g}(\psi_1,\dots ,\psi_s)
\\
\qquad{} +\sum_{y=1}^s
({\rm const})_{y,!}\cdot C^{y,*}_{g}(\psi_1,\dots ,\psi_s)+
\sum_{w=1}^s\sum_{q=1}^{w-1}({\rm const})_{w,q}\cdot
C^{!,w,q}_{g}(\psi_1,\dots ,\psi_s)
\\  \qquad{} +({\rm const})_\sharp C^\sharp_{g}(\psi_1,\dots ,\psi_s) +\sum_{z\in
Z'}a_z C^z_{g}(\psi_1,\dots ,\psi_s).
\end{gather*}
Here each $C^l$, $C^j$ has $|\Delta|\le \sigma-3$.  We then apply the
exact same proof as in the third step in the previous subsection to show that
$({\rm const})_{!}=0$, $({\rm const})_{y,!}=0$, $({\rm const})_{w,q}=0$,
$({\rm const})_\sharp=0$.

}

We have thus shown our claim in the case
$s=\sigma-2$, $\mu<\frac{n}{2}$.
\end{proof}

\begin{proof}[The proof of Lemmas \ref{postpone1st}, \ref{violi}, \ref{postpone2} when $\boldsymbol{s=\sigma-2}$,
$\boldsymbol{\sigma <\frac{n}{2}\!-\!1}$ and $\boldsymbol{\mu=\frac{n}{2}}$.]
  In this ca\-se, we have that $P(g)|_{\Theta_{\sigma -2}}$ can be written out as follows:
\[
P(g)|_{\Theta_{\sigma -2}}=\sum_{l\in L} a_l C^l(g),
\]
where each $C^l(g)$ is a complete contraction in the form:
\begin{gather}
\label{mitra}
{\rm contr}\big(\Delta^{(\alpha_1)}\nabla^{il}W_{ijkl}\otimes\Delta^{(\alpha_2)}
\nabla_{i'l'}W^{i'jkl'}\otimes\Delta^{(\alpha_3)}P^a_a\otimes\dots\otimes
\Delta^{(\alpha_\sigma)}P^a_a\big).
\end{gather}

 This is true by virtue of the fact that each $C^l(g)$, $l\in
\Theta_{\sigma -2}$ must have $\delta=\frac{n}{2}$. We assume with no
 loss of generality that the last $\sigma -2$ factors are arranged so
 that $\alpha_3\ge\alpha_4\ge\dots\ge\alpha_\sigma$.

 We denote by $L_0\subset L$ the index set of the complete
 contractions $C^l(g)$ with $\alpha_3=\dots =\alpha_\sigma=0$.
We claim that:
\begin{gather}
\label{mitra2} \sum_{l\in L_0} a_l C^l(g)=0.
\end{gather}

 We observe that if we can show this, we will then have shown our claim.
 This is because every complete
 contraction $C^l(g)$, $l\in L\setminus L_0$ will have a factor
$\Delta^{(\alpha)} P^a_a$ with $\alpha>0$, hence each complete contraction
 $C^l$, $C^j$, $C^f$ in~(\ref{polis}) will have a factor
 $\nabla^{(p)}\psi_h$ with $p\ge 3$ and two curvature
  factors that are not~$R$ (scalar curvature). Hence, if we can show~(\ref{mitra2}),
we will have shown the remaining cases of both our Lemmas
 when $\mu=\frac{n}{2}$ and $\sigma<\frac{n}{2}-1$.

 We will, in fact show a more general statement: We denote by
$L_\gamma\subset L$ the subsets of complete contractions in the form
 with $\alpha_3=\gamma$, $\alpha_4=\dots=\alpha_\sigma =0$. We will
 then show that for each $\gamma\ge 0$:
\begin{gather}
\label{mitra3} \sum_{l\in L_\gamma} a_l C^l(g)=0.
\end{gather}

 We show the above by an induction. Let us assume that
(\ref{mitra3}) is known for every $\gamma>\gamma_1$. We will then show
 (\ref{mitra3}) for $\gamma=\gamma_1$.
 By our inductive assumption, we may cross out from $P(g)$ the
sublinear combination $\sum_{\gamma>\gamma_1} \sum_{l\in L_\gamma}
a_l C^l(g)$. Now, we consider $I^s_{g}(\psi_1,\dots ,\psi_s)$. We
pick out the sublinear
 combination of $I^{s,\gamma_1}_{g}(\psi_1,\dots ,\psi_s)$ of
complete contractions in the form:
\begin{gather}
\label{fom1}
{\rm contr}\big(\Delta^{(\alpha_1)}\nabla^{il}W_{ijkl}\otimes\Delta^{(\alpha_2)}
\nabla_{i'l'}W^{i'jkl'}\otimes\Delta^{(\gamma_1+1)}\psi_1\otimes
\Delta\psi_2\otimes\dots\otimes\Delta\psi_{\sigma -2}\big)
\end{gather}
(we make the convention that $\alpha_1\ge \alpha_2$).

Clearly, if we can show that modulo complete contractions of
 length $\ge\sigma +1$:
\begin{gather*}
I^{s,\gamma_1}_{g}(\psi_1,\dots ,\psi_s)=0,
\end{gather*}
we will then have shown our inductive claim.

By the transformation law of the Schouten tensor under conformal
re-scalings (see~\cite{alexakis}), and by virtue of our
inductive hypothesis (\ref{mitra}) for $\gamma>\gamma_1$, we deduce that we can write out:
\begin{gather}
\label{poios3erei}
I^s_{g}(\psi_1,\dots,\psi_s)=I^{s,\gamma_1}_{g}(\psi_1,\dots
,\psi_s)+\sum_{u\in U_1} a_u C^u_{g}(\psi_1,\dots ,\psi_s) +
\sum_{u\in U_2} a_u C^u_{g}(\psi_1,\dots ,\psi_s),\!\!
\end{gather}
where each $C^u$ is in the form (\ref{mitra}), and
for $u\in U_1$ $\alpha_h>0$ for at least one $h>0$, while for each
$u\in U_2$ we have that $\alpha_3<\gamma_1$.

 We now use the silly divergence formula for $I^s_{g}$ (recall this formula from~\cite{a:dgciI}),
by integrating by parts with respect to the derivatives on the
factor $\nabla^{(\alpha)}\psi_1$. We denote the silly divergence
formula for $I^s$ by ${\rm silly}[I^s_{g}]$. We focus on the sublinear
combination ${\rm silly}_{*}[I^s_{g}]$ of complete contractions in
${\rm silly}[I^s_{g}]$  in the form:
\begin{gather}
{\rm contr}\big(\nabla^{t_1\dots t_{\gamma_1+1}}\Delta^{(\alpha_1)}
\nabla^{il}W_{ijkl}\otimes
\nabla_{t_1\dots t_{\gamma_1+1}}\Delta^{(\alpha_2)}
\nabla_{i'l'}W^{i'jkl'}\nonumber\\
\qquad{} \otimes\psi_1\otimes
\Delta\psi_2\otimes \dots\otimes\Delta\psi_{\sigma -2}\big).\label{sil}
\end{gather}

Now, let us write out $I^{s,\gamma_1}_{g}(\psi_1,\dots ,
\psi_s)$ in the form:
\begin{gather*}
 I^{s,\gamma_1}_{g}(\psi_1,\dots ,
\psi_s)=\sum_{(p,q)} a_{(p,q)} C^{(p,q)}_{g}(\psi_1,\dots ,\psi_s),
\end{gather*}
where $C^{(p,q)}$ stands for the complete contraction in the form
(\ref{fom1}) with $\alpha_1=p$, $\alpha_2=q$. We then denote by
$C^{(p,q),\sharp}_{g}$ the complete contraction that arises from
$C^{(p,q)}$ by replacing the factor $\Delta^{(\gamma_1+1)}\psi_1$
by a
 factor $\psi_1$ (with no derivatives) and hitting the factor
$\Delta^{(p)}\nabla^{il}W_{ijkl}$ by
 derivatives $\nabla^{t_1\dots t_{\gamma_1+1}}$ and the factor
$\Delta^{(q)}\nabla_{i'l'}W^{i'jkl'}$ by derivatives
$\nabla_{t_1\dots t_{\gamma_1+1}}$
(we are using the Einstein summation convention). We then make three claims. Firstly:
\begin{gather*}
{\rm silly}_{*}[I^s_{g}]=2^{\gamma_1}\sum_{(p,q)} a_{(p,q)}
C^{(p,q),\sharp}_{g} (\psi_1,\dots ,\psi_s).
\end{gather*}

 Now, for each $C^{(p,q),\sharp}$ above, we denote by
$C^{(p,q),\sharp,\iota}_{g}(\psi_1,\dots ,\psi_s)$ the complete
 contraction that arises from it by replacing the two factors
$\nabla^{t_1\dots t_{\gamma_1+1}}\Delta^{(\alpha_1)}
\nabla^{il}W_{ijkl}$, $\nabla_{t_1\dots t_{\gamma_1+1}}
\Delta^{(\alpha_2)}\nabla_{i'l'}W^{i'jkl'}$ by factors
$\nabla^{t_1\dots t_{\gamma_1+1}}\Delta^{(\alpha_1)}
\nabla^{il}R_{ijkl}$,  $\nabla_{t_1\dots t_{\gamma_1+1}}
\Delta^{(\alpha_2)}\nabla_{i'l'}R^{i'jkl'}$. Our second claim is
that modulo complete contractions of length $\ge\sigma +1$:
\begin{gather*}
 \sum_{(p,q)} a_{(p,q)} C^{(p,q),\sharp,\iota}_{g}
(\psi_1,\dots ,\psi_s)=0.
\end{gather*}

Our third claim is that from the above we can deduce that each
$a_{(p,q)}=0$. If we can show the above three claims, then
 clearly (\ref{mitra3}) with $\gamma=\gamma_1$
 will follow.

We begin with our f\/irst claim. For each complete contraction
$C_{g}(\psi_1,\dots ,\psi_s)$ in~(\ref{poios3erei}), we denote by
${\rm sil}[C_{g}(\psi_1,\dots ,\psi_s)]\cdot \psi_1$ the sublinear
combination of complete contractions in the right hand side of
${\rm silly}[I^s_{g}]$. This notation extends to linear combinations.
The silly divergence formula just tells us that:
\begin{gather*}
{\rm sil}[I^{s,\gamma_1}_{g}(\psi_1,\dots
,\psi_s)]+\sum_{u\in U_1} a_u {\rm sil}[C^u_{g}(\psi_1,\dots ,\psi_s)] +
\sum_{u\in U_2} a_u {\rm sil}[C^u_{g}(\psi_1,\dots ,\psi_s)]=0.
\end{gather*}

If we can show that for $u\in U_1\bigcup U_2$, ${\rm sil}[C^u_{g}]$
 contains no complete contractions in the form~(\ref{sil}), our f\/irst
 claim will then follow.

This claim is trivial for the complete contractions $C^u_{g}$,
$u\in U_1$. We notice that if $C^u_{g}$ has a~factor
$\Delta^{(\alpha)}\psi_h$, $\alpha>1$, $h>1$, then each complete
contraction in ${\rm sil}[C^u_{g}]$ will have a factor
$\nabla^{(p)}\psi_h$, $p\ge 2\alpha>2$.
 On the other hand, for each $u\in U_2$, we observe that each complete contraction in
${\rm sil}[C^u_{g}]$ will have less than $\gamma_1+3$ indices in the
f\/irst factor $\nabla^{(m)}W_{ijkl}$ contracting against indices in
the second
 factor $\nabla^{(m')}W_{ijkl}$. So we indeed see that
${\rm sil}[C^u_{g}]$, $u\in U_2$ has no complete contraction of the form~(\ref{sil}). Therefore we have shown our f\/irst claim.

 Now, we check the second claim. We write out ${\rm silly}[I^s_{g}
(\psi_1,\dots ,\psi_s)]$ as a linear combination of complete
contractions in the form (\ref{linisymric}). We def\/ine
${\rm silly}_{\rm spec} [I^s_{g}(\psi_1,\dots ,\psi_s)]$ to stand for the
 sublinear combination of complete contractions with $s-1$ factors
$\Delta\psi_2,\dots ,\Delta\psi_s$, one factor $\psi_1$ and two
factors  $\nabla^{(m)}R_{ijkl}$, $\nabla^{(m')}R_{i'j'k'l'}$ with
$\gamma_1+3$ particular contractions between them.  (We write
$\nabla^{(m)}R_{ijkl}$, $\nabla^{(m')}R_{i'j'k'l'}$ but there is no
restriction on the particular contractions of the indices
 ${}_i,{}_j,{}_k,{}_l,{}_{i'},{}_{j'},{}_{k'},{}_{l'}$, i.e.\ we can
have ${}_i,{}_k$
 contracting between themselves etc.) Since the silly divergence
 formula must hold formally, we deduce that:
\begin{gather}
\label{italos}
{\rm silly}_{\rm spec}[I^s_{g}(\psi_1,\dots ,\psi_s)]=0,
\end{gather}
modulo complete contractions of length $\ge\sigma +1$.

Moreover, clearly any complete contraction $C_{g}(\psi_1,\dots ,
\psi_s)$ in the form (\ref{mitra}) with either
 a factor $\nabla^{(p)}\psi_h$, $p\ge 3$ or with less than $\gamma_1+3$ pairs
 of particular contractions between the f\/irst two factors
$\nabla^{(m)}W_{ijkl}$, $\nabla^{(m')}W_{ijkl}$ cannot contribute to the
 sublinear combination ${\rm silly}_{\rm spec}[I^s_{g}(\psi_1,\dots ,
\psi_s)]$. Hence, we derive that ${\rm sil}[C^u_{g}]$, $u\in U_1\bigcup U_2$
do not contribute to the
 sublinear combination ${\rm silly}_{\rm spec}[I^s_{g}]$, whereas by virtue of
 the decomposition of the factor $\nabla^{il}W_{ijkl}$, we have that:
\begin{gather*}
 {\rm silly}_{\rm spec}[I^s_{g}(\psi_1,\dots ,\psi_s)]=
\sum_{(p,q)} 2^{\gamma_1+1}\frac{n-3}{n-2}a_{(p,q)}
C^{(p,q),\sharp,\iota}_{g} (\psi_1,\dots ,\psi_s)
 +\sum_{t\in T}
a_t C^t_{g}(\psi_1,\dots, \psi_s),
\end{gather*}
where each is in the form (\ref{linisymric}) with either two factors of
 the form $\nabla^{(p)}R$, $\nabla^{(p')}R$,\footnote{Recall
 that $R$ here is the scalar curvature.}
 or one factor of the form
$\nabla^{(p)}R$ and one factor of the form $\nabla^{(m)}R_{ijkl}$ where
none of
 the indices~${}_i$,~${}_j$,~${}_k$,~${}_l$ are contracting between themselves. In fact,
in order to distinguish these two cases we break the index set~$T$
 into subsets $T_1$, $T_2$ accordingly.

Now, since (\ref{italos}) holds formally we deduce that modulo complete contractions of
length $\ge\sigma +1$:
\begin{gather*}
\sum_{t\in T_1} a_t C^t_{g}(\psi_1,\dots, \psi_s)=0.
\end{gather*}
Then, again since (\ref{italos}) holds formally, we deduce that
modulo complete contractions
 of length $\ge\sigma +1$:
\begin{gather*}
\sum_{t\in T_2} a_t C^t_{g}(\psi_1,\dots,\psi_s)=0.
\end{gather*}

 Therefore, we have shown that modulo complete contractions of length
$\ge\sigma +1$:
\begin{gather*}
\sum_{(p,q)} a_{(p,q)} C^{(p,q),\sharp,\iota}_{g}
(\psi_1,\dots ,\psi_s)=0.
\end{gather*}

 Clearly, since the above holds formally, we can deduce
 that for each dif\/ferent pair $(p,q)$ we must have:
\begin{gather*}
a_{(p,q)} C^{(p,q),\sharp,\iota}_{g} (\psi_1,\dots
,\psi_s)=0,
\end{gather*}
and since each complete contraction $C^{(p,q),\sharp,\iota}_{g}
(\psi_1,\dots ,\psi_s)$ is not identically the zero contraction,
we
 deduce that each $a_{(p,q)}$ must be zero. This was our third claim.
 We have thus shown our inductive statement, and we have proven our claim
 in the case $s=\sigma-2$, $\sigma<\frac{n}{2}-1$, $\mu=\frac{n}{2}$.
\end{proof}

\begin{proof}[The proof of Lemmas \ref{postpone1st}, \ref{violi}, \ref{postpone2}
 when $\boldsymbol{s=\sigma-2}$ and $\boldsymbol{\sigma=\frac{n}{2}-1}$.]
 We observe that \linebreak any complete contraction $C^l(g)$ in $P(g)|_{\sigma}$, when
$\sigma=\frac{n}{2}-1$ must be in the form:
\[
{\rm contr}\big(\nabla^{(m)}W_{ijkl}\otimes\nabla^{(m')}W_{i'j'k'l'}\otimes \nabla^{(p_1)}P_{ab}\otimes\dots\otimes
\nabla^{(p_{\frac{n}{2}-3})}P_{a'b'}\big)
\] must have exactly two
derivatives (i.e.\ $m+m'+p_1+\dots+p_{\frac{n}{2}-3}=2$)
 and moreover
  $\delta\le \frac{n}{2}-1$.\footnote{It follows easily
   that in this case $\mu\le\frac{n}{2}$.} We
distinguish two cases.  Either
 $\mu =\frac{n}{2}-1$ or $\mu< \frac{n}{2}-1$.
We start with the f\/irst case. Then we see that:
\begin{gather*}
P(g)|_{\Theta_{\sigma -2}}=({\rm const})\cdot C(g),
\end{gather*}
where $C(g)$ is the complete contraction:
\begin{gather}
\label{outthere}
{\rm contr}\big(\nabla^lW_{ijkl}\otimes
\nabla_{l'}W^{ijkl'} \otimes (P^a_a)^{\sigma -2}\big).
\end{gather}

 We will show that $({\rm const})=0$.
This follows easily by the silly divergence formula. We
consider $I^s_{g}(\psi_1,\dots ,\psi_s)$. It follows that:
\begin{gather*}
 I^s_{g}(\psi_1,\dots ,\psi_s)=s!\cdot C_{g}(\psi_1,
\dots ,\psi_s),
\end{gather*}
where $C_{g}(\psi_1,\dots ,\psi_s)$ is the contraction:
\begin{gather*}
{\rm contr}\big(\nabla^lW_{ijkl}\otimes \nabla_{l'}W^{ijkl'}
\otimes \Delta\psi_1\otimes\dots\otimes\Delta\psi_{\sigma -2}\big).
\end{gather*}

We consider ${\rm silly}[I^s_{g}]$, with integration by parts with
respect to $\nabla^{(p)}\psi_1$, and we focus on the sublinear
combination:
\begin{gather*}
 {\rm silly}_{+}[I^s_{g}(\psi_1,\dots ,\psi_s)=
({\rm const})'\!\cdot \! {\rm contr}\big(\nabla^{la}W_{ijkl}\otimes
\nabla_{l'a}W^{ijkl'} \otimes \psi_1\otimes
 \Delta\psi_2\otimes\dots\otimes\Delta\psi_{\sigma -2}\big).
\end{gather*}

It follows by the same arguments as before that
${\rm silly}_{+}[I^s_{g}]=0$, modulo complete contractions of length
$\ge\sigma +1$ and that $({\rm const})'=2 \cdot s!({\rm const})$, and hence we
deduce our claim in the case $\mu=\frac{n}{2}-1$.

 Now, the case $\mu<\frac{n}{2}-1$. We notice that in this case, the
sublinear combination of complete contractions $C^l(g)$, $l\in
\Theta_{\sigma -2}$ that have $\sigma -2$ factors $P^a_a$ is
precisely of the form $({\rm const})\cdot C(g)$, where $C(g)$ is in the
form (\ref{outthere}) (since $P(g)|_{\Theta_{\sigma-1}}$ is assumed
``good''). We can
then prove our claim in this case by
 exactly applying the method of the case $\sigma=\frac{n}{2}-1$ when $s<\sigma-2$.
\end{proof}

\subsection[The proof of Lemmas \ref{postpone1st}, \ref{violi}, \ref{postpone2}
 when $s=\sigma-1$, and Lemma \ref{tistexnes} when $s=\sigma-1$]{The proof of Lemmas \ref{postpone1st}, \ref{violi}, \ref{postpone2}
 when $\boldsymbol{s=\sigma-1}$,\\ and Lemma \ref{tistexnes} when $\boldsymbol{s=\sigma-1}$}

 In this case, we follow the same pattern as in the two previous
ones. We begin by a trivial observation regarding the maximum
value that $\mu$ can have\footnote{Recall that $\mu$ stands for
the minimum value of $\delta=\delta_W+\delta_P$ among the complete
contractions $C^l(g)$, $l\in \Theta_s$.}. By virtue of the formula
 $I^{\sigma-1}_{g}(\phi):=\frac{d^{\sigma-1}}{dt^{\sigma-1}}|_{t=0}e^{nt\phi}P(e^{2t\phi}g)$,\footnote{Thus
in particular the terms of length $\sigma-1$ in $I^{\sigma-1}_{g}(\phi)$
arise from the sublinear combination $P(g)_{\Theta_{\sigma-1}}$
by just replacing each factor $\nabla^{(a)}_{r_1\dots r_a}P_{ij}$
by $-\nabla^{(a+2)}_{r_1\dots r_a ij}\phi$.} and since we
observe that any complete contraction in the form~(\ref{drassi})
with length $\sigma$ and $\sigma -1$ factors $P^a_a$ {\it must be
zero, modulo complete contractions of length~$\ge\sigma +1$}, we
observe that necessarily $\mu\le \frac{n}{2}-1$ in this case.

  We again distinguish the two cases $\mu=\frac{n}{2}-1$ and
$\mu<\frac{n}{2}-1$. We start with the second case.

\begin{proof}[Proof of Lemma \ref{tistexnes} when $\boldsymbol{s=\sigma-1}$, $\boldsymbol{\mu<\frac{n}{2}-1}$.]
We
focus on the sublinear combination of complete contractions
$C^l(g)$ in $P(g)|_{\Theta_{\sigma -1}}$ that have $\sigma -2$
factors $P^a_a$. We denote their index set by $\Theta^{+}_{\sigma
-1}\subset \Theta_{\sigma -1}$. By the usual ``manual''
construction of explicit divergences it follows
 that we can subtract a divergence ${\rm div}_i\sum_{h\in H} a_h C^{h,i}(g)$
 from $P(g)|_{\Theta^{+}_{\sigma -1}}$ so that modulo complete
contractions of length $\ge\sigma +1$:
\begin{gather*}
P(g)|_{\Theta^{+}_{\sigma -1}}-{\rm div}_i\sum_{h\in H} a_h
C^{h,i}(g)= ({\rm const})_V\cdot C^V(g)+ \sum_{t\in T} a_t
C^t(g)+\sum_{y\in Y} a_y C^y(g),
\end{gather*}
where $C^V(g)$ is the complete contraction:
\begin{gather*}
{\rm contr}\big(\Delta^{\frac{n}{2}-(\sigma -2)-3}\nabla^{il}W_{ijkl}\otimes
P^{jk}\otimes (P^a_a)^{\sigma -2}\big),
\end{gather*}
and where each $C^y(g)$ is in the form (\ref{drassi}) with at
least two factors $\nabla^{(m)}W_{ijkl}$. Also, each $C^t(g)$ is
in the form (\ref{drassi}) with $\delta\ge \mu$ and also with
less than $\sigma -2$ factors $P^a_a$.

 Now, we observe that $C^V(g)$ has $\delta=\frac{n}{2}-1\ge \mu+1$. Therefore, if we
 def\/ine $C^{V,i}(g)$ to stand for the vector f\/ield that arises from $C^V(g)$ by
erasing ${}^i$ and making ${}_i$ in the f\/irst factor into a~free
index, we will then have:
\begin{gather*}
 C^V(g)-{\rm div}_i C^{V,i}(g)=\sum_{y\in Y} a_y C^y(g)+
\sum_{d\in D} a_d C^d(g).
\end{gather*}
Here each $C^y(g)$ is as above, while each $C^d(g)$ has $\delta\ge
\mu$ and also has strictly less than $\sigma -2$ factors $P^a_a$.
Thus, by subtracting the divergence ${\rm div}_i\Big[\sum_{h\in H} a_h
C^{h,i}(g)+({\rm const})_V\cdot C^{V,i}(g)\Big]$ from $P(g)$, we may assume
 with no loss of generality that $P(g)|_{\Theta_{\sigma -1}}$ contains
 only complete contractions with at most $\sigma -3$ factors $P^a_a$. We will
be using this fact below.
\end{proof}

\begin{proof}[Proof of Lemmas \ref{postpone1st}, \ref{violi}, \ref{postpone2}
when $\boldsymbol{s=\sigma\!-\!1}$, $\boldsymbol{\mu<\frac{n}{2}\!-\!1}$.] We
focus on $I^s_{g}(\psi_1,{\dots},\psi_s)$;\footnote{Recall that
$I^s_g(\phi):=\frac{d^s}{dt^s}|_{t=0} (e^{nt\phi}P^{2t\phi}g)$
 and that $I^s_g(\psi_1,\dots,\psi_s)$
arises from $I^s_g(\phi)$ by just polarizing the function~$\phi$.}
for the global conformal
invariants~$P(g)$ above, and we decompose the Weyl tensors (i.e.\
we write~$I^s_g$ as a linear combination in the form~(\ref{polis})).
 In the notation of the introduction,
we have that $F^{*}=F^1_{\sigma -1}\bigcup F^1_{\sigma-2}$.

 We f\/irst focus our attention on the sublinear combination $F^1_{\sigma -1}$. As
in the previous cases, it follows that we can write out:
\begin{gather*}
\sum_{f\in F^1_{\sigma -1}} a_f C^f_{g}(\psi_1,\dots
,\psi_s)= ({\rm const})_{+}\cdot C^{+}_{g}(\psi_1,\dots ,\psi_s)+
\sum_{u=1}^s a_u C^u_{g}(\psi_1,\dots ,\psi_s),
\end{gather*}
where $C^{+}_{g}(\psi_1,\dots ,\psi_s)$ stands for the complete
contraction:
\begin{gather*}
{\rm contr}\big(\Delta^{(\frac{n}{2}-\sigma)}R\otimes\Delta\psi_1\otimes\dots
\otimes\Delta\psi_s\big),
\end{gather*}
 while $C^u_{g}(\psi_1,\dots ,\psi_s)$ stands for the complete contraction:
\begin{gather*}
{\rm contr}\big(\Delta^{(\frac{n}{2}-\sigma+1)}\psi_u\otimes\Delta\psi_1\otimes
\dots\hat{\Delta\psi_u}\dots\otimes\Delta\psi_s\otimes R\big).
\end{gather*}

 We observe that each of the above complete contraction has
$\delta=\frac{n}{2}+1\ge \mu +2$.
 As before, we can explicitly construct a divergence
${\rm div}_i\sum_{h\in H} a_h C^{h,i}_{g}(\psi_1,\dots ,\psi_s)$ such that:
\begin{gather*}
 \sum_{f\in F^1_{\sigma -1}} a_f C^f_{g}(\psi_1,\dots ,\psi_s)-
{\rm div}_i\sum_{h\in H} a_h C^{h,i}_{g}(\psi_1,\dots ,\psi_s)  =
\sum_{u=1}^s({\rm const})_{*,u} C^{*,u}_{g}(\psi_1,\dots ,\psi_s)
\\ \qquad{} +\sum_{q=1}^s\sum_{w=1}^{q-1} ({\rm const})_{*,(q,w)}
C^{*,(q,w)}_{g}(\psi_1,\dots ,\psi_s)+ \sum_{z\in Z'} a_z
C^z_{g}(\psi_1,\dots ,\psi_s),
\end{gather*}
where $C^{*,u}_{g}(\psi_1,\dots ,\psi_s)$ stands for the complete
 contraction:
\begin{gather*}
{\rm contr}\big(\Delta^{(\frac{n}{2}-\sigma-1)}R\otimes\Delta^{(2)}\psi_u\otimes\Delta
\psi_1\otimes \dots\otimes\hat{\Delta\psi_u}\otimes\dots\otimes\Delta\psi_s\big),
\end{gather*}
 while $C^{*,(q,w)}_{g}(\psi_1,\dots ,\psi_s)$
stands for the complete contraction
\begin{gather*}
{\rm contr}\big(\Delta^{(\frac{n}{2}-\sigma)}\psi_q\otimes\Delta^{(2)}\psi_w
\otimes\Delta
\psi_1\otimes \dots\otimes\hat{\Delta\psi_q}\otimes \dots\otimes\hat{\Delta\psi_w}\otimes\dots
\otimes\Delta\psi_s\big),
\end{gather*}
and $\sum_{z\in Z'} a_z C^z_{g}(\psi_1,\dots ,\psi_s)$ stands for a
 generic linear combination of complete contractions with
$|\Delta|\le\sigma -3$ and $\delta\ge\mu +1$.

 We next focus our attention on the sublinear combination
$\sum_{f\in F^1_{\sigma -2}} a_f C^f_{g}(\psi_1,\dots ,\psi_s)$.
Since $P(g)|_{\Theta_{\sigma -1}}$ has {\it no} complete
contractions
 with $\sigma -2$ factors $P^a_a$, it follows that any
 $C^f_{g}(\psi_1,\dots,\psi_s)$, $f\in F^1_{\sigma -2}$ must have
$\delta\ge \mu +2$. Therefore, as in the previous subsection, we
can construct
 a~vector f\/ield $\sum_{h\in H} a_h C^{h,i}_{g}(\psi_1,\dots ,\psi_s)$
so that, modulo complete contractions of length $\ge\sigma +1$:
\begin{gather*}
\sum_{f\in F^1_{\sigma -2}} a_f C^f_{g}(\psi_1,\dots ,\psi_s)
-\sum_{h\in H} a_h C^{h,i}_{g}(\psi_1,\dots ,\psi_s)=\sum_{u=1}^s({\rm const})'_{*,u} C^{*,u}_{g}(\psi_1,\dots ,\psi_s)
\\  +
\sum_{q=1}^s\sum_{w=1}^{q-1} ({\rm const})'_{*,(q,w)}
C^{*,(q,w)}_{g}(\psi_1,\dots ,\psi_s)
+ \sum_{b\in B} a_b C^b_{g}(\psi_1,\dots ,\psi_s)
+\sum_{z\in Z'} a_z C^z_{g}(\psi_1,\dots,\psi_s).
\end{gather*}
Here the complete contractions $C^{*,u}$, $C^{*,(q,w)}$ are the
same as above, while $\sum_{b\in B} a_b C^b_{g}(\psi_1,\dots
,\psi_s)$ stands for a linear combination of complete contractions
in the form~(\ref{linisymric}) with length $\sigma$, $\delta\ge\mu
+2$, $q=0$ and $|\Delta|=\sigma -2$.

  Finally, we focus on the sublinear combinations in
$I^s_{g}(\psi_1,\dots ,\psi_s)$ which have $q=0$ (i.e.\ we consider
the complete contractions that are indexed in the sets
$L_\mu$, $J$ in (\ref{polis})\footnote{These sets index the complete contractions in
$I^s_g$  with no factors $\nabla^{(p)}{\rm Ric}$.}).  Now, since we have that each $C^l(g), l\in
 \Theta_{\sigma -1}$ has less than $\sigma -2$ factors $P^a_a$, it
follows that all the complete contractions in $L_\mu$, $J$ in
(\ref{polis}) must have $|\Delta|\le\sigma-3$.

Then (by the usual construction), we
can construct a linear combination of vector f\/ields, $\sum_{h\in
H} a_h C^{h,i}_{g}(\psi_1,\dots ,\psi_s)$ so that, modulo complete
 contractions of length $\ge\sigma +1$:
\begin{gather*}
\sum_{b\in B} a_b C^b_{g}(\psi_1,\dots ,\psi_s)- {\rm div}_i \sum_{h\in
H} a_h C^{h,i}_{g}(\psi_1,\dots ,\psi_s)
\\
\qquad{}= \sum_{u=1}^s({\rm const})_{L,u} C^{L,u}_{g}(\psi_1,\dots ,\psi_s)+
\sum_{z\in Z'} a_z C^z_{g}(\psi_1,\dots,\psi_s),
\end{gather*}
where $C^{L,u}_{g}(\psi_1,\dots ,\psi_s)$ is the complete
 contraction:
\begin{gather*}
{\rm contr}\big(\Delta^{(\frac{n}{2}-(\sigma
-2)-3)}\nabla^{il}R_{ijkl}\otimes \nabla^{jk}\psi_u\otimes
\Delta\psi_1\otimes\dots \otimes\hat{\Delta\psi_u}\otimes
\dots\otimes\Delta\psi_s\big),
\end{gather*}
which has $\delta\ge\mu +2$, and where $\sum_{z\in Z'} a_z
C^z_{g}(\psi_1,\dots,\psi_s)$ is as above.

Moreover, we def\/ine $C^{L,u,i}_{g}(\psi_1,\dots ,\psi_s)$ to
 stand for the vector f\/ield that arises from $C^{L,u}_g(\psi_1$, $\dots ,\psi_s)$
  by erasing the index ${}^i$ and making ${}_i$ into a free index. We then observe that, modulo
complete contractions of length $\ge\sigma +1$:
\begin{gather*}
C^{L,u}_{g}(\psi_1,\dots ,\psi_s)-{\rm div}_i
C^{L,u,i}_{g}(\psi_1,\dots , \psi_s)=\sum_{z\in Z'} a_z
C^z_{g}(\psi_1,\dots,\psi_s).
\end{gather*}

In conclusion, we have shown that we can subtract a linear
 combination of divergences from $I^s_{g}(\psi_1,\dots ,
\psi_s)$ so that, modulo complete contractions of length
$\ge\sigma +1$:
\begin{gather*}
I^s_{g}(\psi_1,\dots ,\psi_s)-{\rm div}_i C^{h,i}_{g}(\psi_1,\dots ,
\psi_s)=
\sum_{l\in L_\mu} a_l C^{l,\iota}_g(\psi_1,\dots,\psi_s)+\sum_{j\in J} a_j C^j_g(\psi_1,\dots,\psi_s)
\\\qquad{} +
\sum_{f\in F^1\setminus (F^1_{\sigma -1}\bigcup
F^1_{\sigma -2})} a_f C^f_{g}(\psi_1,\dots ,\psi_s)
 +\sum_{u=1}^s({\rm const})_{*,u} C^{*,u}_{g}(\psi_1,\dots ,\psi_s)\\
 \qquad{} +
\sum_{q=1}^s\sum_{w=1}^{q-1} ({\rm const})_{*,(q,w)}
C^{*,(q,w)}_{g}(\psi_1,\dots ,\psi_s)
 + \sum_{z\in Z'} a_z
C^z_{g}(\psi_1,\dots ,\psi_s).
\end{gather*}
(The terms in the second line are {\it not} generic linear
combinations~-- they stand for the original linear combinations
in $I^s_g(\psi_1,\dots,\psi_s)$ as in (\ref{polis}).)

 Then, using the super divergence formula as in the previous subsection, we deduce
that $({\rm const})_{*,u}=0$ for every $u=1,\dots ,s$ and
$({\rm const})_{*,(q,w)}=0$ for every $q=1,\dots ,s$, $w=1,\dots ,q-1$.

 Therefore, we have proven  Lemmas \ref{theatro}, \ref{violi}, \ref{postpone2}
 when $s=\sigma-1$ and $\mu<\frac{n}{2}-1$.
 \end{proof}

 Now, we consider the subcase where $\mu=\frac{n}{2}-1$. We see that
in this case (modulo introducing correction terms with length
$\sigma$ and two factors $\nabla^{(m)}W_{ijkl}$), we can write out
$P(g)|_{\Theta_{\sigma -1}}$:
\begin{gather*}
P(g)|_{\Theta_{\sigma -1}}=\sum_{\gamma\in \Gamma}
a_\gamma C^\gamma_{g}(\psi_1,\dots,\psi_s)+\sum_{\epsilon\in E}
a_\epsilon C^\epsilon_{g}(\psi_1,\dots ,\psi_s),
\end{gather*}
where each $C^\gamma(g)$ is in the form:
\begin{gather}
\label{spaw}
{\rm contr}\big(\Delta^{(\alpha)}\nabla^{il}W_{ijkl}\otimes\Delta^{(\beta)}
P^{jk}\otimes
\Delta^{(\rho_1)}P^a_a\otimes\dots\otimes\Delta^{(\rho_{\sigma
-2})}P^a_a\big)
\end{gather}
(we are making the convention that $\rho_1\ge\dots\ge\rho_{\sigma -2}$),
while each $C^\epsilon(g)$ is in the form:
\begin{gather}
\label{spaw2}
{\rm contr}\big(\Delta^{(\alpha)}\nabla^{il}W_{ijkl}\otimes\Delta^{(\beta_1)}
\nabla^tP^j_t\otimes\Delta^{(\beta_2)}\nabla^wP^k_w\otimes
\Delta^{(\rho_1)}P^a_a\otimes\dots\otimes\Delta^{(\rho_{\sigma
-3})}P^a_a\big)
\end{gather}
(we are making the convention that $\rho_1\ge\dots\ge\rho_{\sigma -3}$).

 Now, a small observation. For each $C^\gamma(g)$ with $\beta>0$,
 we can construct ``by hand'' a vector f\/ield $C^{\gamma,i}(g)$ so that:
\begin{gather}
\label{oliver}
C^\gamma(g)-{\rm div}_i C^{\gamma,i}(g)=\sum_{y\in Y}
a_yC^y(g)+ \sum_{t\in T} a_t C^t(g),
\end{gather}
where each $C^y(g)$ is a complete contraction in the form
(\ref{spaw2})
 and $\sum_{t\in T} a_t C^t(g)$ is a linear combination of
 complete contractions in the form (\ref{drassi}) with at least two factors $
 \nabla^{(m)}W_{ijkl}$. The above holds modulo complete contractions of length $\ge\sigma +1$.

 In view of (\ref{oliver}), we can assume that each $C^\gamma(g)$
 in $P(g)|_{\Theta_{\sigma -1}}$ has $\beta=0$.

We then def\/ine $\Gamma_\kappa$ to stand for the index set of the
 complete contractions in the form (\ref{spaw}) with $\beta =0$,
$\rho_1=\kappa\ge 0$ and $\rho_2,\dots ,\rho_{\sigma -2}=0$. We
observe that, by def\/inition:
\begin{gather*}
\sum_{\gamma\in \Gamma_\kappa} a_\gamma
C^\gamma(g)=({\rm const})_\kappa \cdot C^\kappa(g),
\end{gather*}
where $C^\kappa(g)$ is in the form (\ref{spaw}) with
$\beta=0$, $\rho_1=\kappa$, $\rho_2=\dots =\rho_{s-2}=0$ and with
$\alpha =\frac{n}{2}-(\sigma -2)-3$.

 We will show that under the above assumptions, for each $\kappa\ge 0$:
\begin{gather}
\label{oliverkappa}
 \sum_{\gamma\in \Gamma_\kappa} a_\gamma
C^\gamma_{g}(\psi_1,\dots, \psi_s)=0,
\end{gather}
modulo complete contractions of length $\ge\sigma +1$.

 If we can show the above, we then see that we will have shown
 Lemmas \ref{theatro}, \ref{violi}, \ref{postpone2} in
this setting.
 This is true since if we consider any
$P(g)$ with $P(g)|_{\Theta_{\sigma -1}}$ that satisf\/ies
(\ref{oliverkappa})  and we then write out
 $I^s_{g}(\psi_1,\dots ,\psi_s)$ by decomposing the Weyl tensors, then each complete
 contraction~$C^{l,\iota}$, $C^f$, $C^j$ in
$I^s$ must have at least two factors $\nabla^{(t)}\psi_h$ with
$t\ge 3$
 and at least one factor either in the form $\nabla^{(m)}R_{ijkl}$, $m\ge 2$
 or in the form $\nabla^{(p)}{\rm Ric}$, $p\ge 2$.

  We show (\ref{oliverkappa}) by an induction. We assume that
(\ref{oliverkappa}) holds for each $\kappa>\kappa_1$ and we will show
 it for $\kappa=\kappa_1$.
 We consider the silly divergence formula for
$I^s_{g}(\psi_1,\dots ,\psi_s)$~-- we denote it by
${\rm silly}[I^s_{g}(\psi_1,\dots ,\psi_s)]$~-- and we focus on the
sublinear combination ${\rm silly}_{*}[I^s_{g}]$ that consists of the
complete contractions in the form:
\begin{gather}
\label{olaria}
{\rm contr}\big(\nabla_{t_1\dots
t_{\kappa_1}}\Delta^{(\frac{n}{2}-(\sigma
-2)-3-\kappa_1)}\nabla^{il}R_{ijkl}\otimes\nabla^{t_1\dots
t_{\kappa_1} jk}\psi_1\otimes\psi_2\otimes\Delta\psi_3
\otimes\dots\otimes\Delta\psi_{s-2}\big).
\end{gather}

We then make two claims. Firstly, that modulo complete
contractions of length $\ge\sigma +1$:
\begin{gather*}
 {\rm silly}_{*}[I^s_{g}(\psi_1,\dots ,\psi_s)]=0,
\end{gather*}
and secondly that:
\begin{gather*}
{\rm silly}_{*}[I^s_{g}(\psi_1,\dots ,\psi_s)]=
2^{\kappa_1}\frac{n-3}{n-2}({\rm const})_{\kappa_1}\cdot
C^{\kappa_1}_{g}(\psi_1,\dots ,\psi_s),
\end{gather*}
where $C^{\kappa_1}_{g}$ is in the form (\ref{olaria}). Since
$C^{\kappa_1}$ is clearly not identically zero, we will then
deduce that $({\rm const})_{\kappa_1}=0$, and we will have shown our
inductive claim.

  We begin with the f\/irst claim. We initially denote by
${\rm silly}_{+}[I^s_{g}]$ the sublinear combination in ${\rm silly}[I^s_{g}]$
that consists of complete contractions with length $\sigma$ and a
factor $\psi_2$, factors $\Delta\psi_3,\dots ,\Delta\psi_{s-2}$, a
factor $\nabla^{(\kappa_1+2)}\psi_1$ {\it with no internal
contractions} and a factor $\nabla^{(m)}R_{ijkl}$ (where now we
are allowing any two indices in this last factor to be contracting
against each other). In other words, we are looking at the
sublinear combination of complete contractions in the form:
\begin{gather}
\label{inje}
{\rm contr}\big(\nabla^{(m)}R_{ijkl}\otimes
\nabla^{(\kappa_1+2)}\psi_1\otimes\psi_2\otimes\Delta\psi_3\otimes\dots\otimes
\Delta\psi_s\big),
\end{gather}
where $\nabla^{(\kappa_1+2)}\psi_1$ has no internal contraction.

 Since the silly divergence formula holds formally, we
deduce that:
\begin{gather*}
{\rm silly}_{+}[I^s_{g}]=0,
\end{gather*}
modulo complete contractions of length $\ge\sigma +1$.

We are now in a position to show our f\/irst claim above. By just
applying the decomposition of the Weyl tensor, we calculate:
\begin{gather}
\label{toutou}
{\rm silly}_{+}[I^s_{g}(\psi_1,\dots
,\psi_s)]={\rm silly}_{*}[I^s_{g}(\psi_1,\dots ,\psi_s)]+ ({\rm const})'\cdot
C'_{g}(\psi_1,\dots ,\psi_s),
\end{gather}
where $C'_{g}(\psi_1,\dots ,\psi_s)$ is the complete contraction:
\begin{gather*}
{\rm contr}\big(\nabla^{t_1\dots t_{\kappa_1
+2}}\Delta^{\rm (power)} R\otimes\nabla_{t_1\dots
t_{\kappa_1+2}}\psi_1\otimes
\psi_2\otimes\Delta\psi_3\otimes\dots\otimes\Delta\psi_s\big).
\end{gather*}
Now since ${\rm silly}_{+}[I^s_{g}(\psi_1,\dots ,\psi_s)]=0$, we may
derive that the second sublinear combination in~(\ref{toutou})
vanishes separately\footnote{Using the operation
${\rm Sub}_\omega$, def\/ined in the
section on technical tools, in the Appendix of~\cite{alexakis}.}. Hence we derive our f\/irst claim, that
${\rm silly}_{*}[I^s_{g}(\psi_1,\dots ,\psi_s)]=0$.

 We now proceed to the second claim.
We would like to understand how the sublinear combination
${\rm silly}_{+}[I^s_{g}]$ is related to $I^s_{g}(\psi_1,\dots
,\psi_s)$. Given
that $\mu=\frac{n}{2}-1$ in our case,  while
${\rm silly}_{+}[I^s_{g}]$ consists of complete contractions with
factors $\Delta\psi_3,\dots ,\Delta\psi_{s-1}$, it follows that a
complete contraction $C_{g}(\psi_1,\dots ,\psi_{s-1})$ in
$I^s_{g}$ (where $C_{g}$ is in the form (\ref{onlyweyl})) can only
contribute to ${\rm silly}_{+}[I^s_{g}]$ if it has factors
$\Delta\psi_3,\dots ,\Delta\psi_{s-1}$. (Other\-wise, it must have
at least one factor $\Delta^y\psi_v$, $v\ge 3$ with $y\ge 2$).
Moreover, since the complete contractions in ${\rm silly}_{+}[I^s_{g}]$
each have a factor $\nabla^{(\kappa_1+2)}\psi_1$ {\it with no
internal contractions}, it follows that a complete contraction
$C_{g}(\psi_1,\dots ,\psi_{s-1})$ in $I^s_{g}$ (where $C_{g}$ is
in the form (\ref{onlyweyl})) can only contribute to
${\rm silly}_{+}[I^s_{g}]$ if it has a factor $\nabla^{(y)}\psi_1$ with
no internal contractions (in the form (\ref{onlyweyl})). Moreover,
since we are considering complete contractions in
$I^s_{g}(\psi_1,\dots ,\psi_s)$ with $\delta=\frac{n}{2}-1$, it
follows that $y=2$, and therefore the two indices~${}_a$,~${}_b$ in
the factor~$\nabla^{(2)}_{ab}\psi_1$ must be contracting against
two internal indices in the one factor $\nabla^{(m)}W_{ijkl}$
(other\-wise, if they contracted against  at least one derivative index, we
must have at least one pair of antisymmetric indices~${}_i$,~${}_j$
or~${}_k$,~${}_l$ in $\nabla^{(m)}W_{ijkl}$ both involved in an
internal contraction, but such a~complete contraction must clearly
be zero modulo complete contractions of length~$\ge\sigma +1$).

{\sloppy Therefore, we have that the complete contractions in
$I^s_{g}(\psi_1,\dots ,\psi_s)$ that contribute to
${\rm silly}_{+}[I^s_{g}(\psi_1,\dots ,\psi_s)]$ must be in the form:
\begin{gather*}
{\rm contr}\big(\Delta^{(\alpha)}\nabla^{il}W_{ijkl}\otimes\nabla^{jk}\psi_1
\otimes\Delta^{(\beta)}\psi_2\otimes\Delta\psi_3\otimes\dots\otimes\Delta\psi_s\big).
\end{gather*}
 We denote the above complete contraction
by $C^{\alpha,\beta}_{g}(\psi_1,\dots ,\psi_s)$. We recall that by
our inductive assumption, each $C^{\alpha,\beta}_{g}$ in $I^s_{g}$
must have~$\beta\le\kappa_1$.

}

 Clearly, we also observe that for each $C^{\alpha,\beta}_{g}$ above,
 the contribution of $C^{\alpha,\beta}_{g}$ to ${\rm silly}_{+}[I^s_{g}]$
 can only be in ${\rm silly}[C^{\alpha,\beta}_{g}]$ if we integrate by
 parts all the $2\beta$ indices $({}^{a_1}, {}_{a_1}),\dots
 ,({}^{a_\beta},{}_{a_\beta})$ in the factor $\Delta^{(\beta)}\psi_2$
and we make the derivatives $\nabla^{a_1},\dots ,\nabla^{a_\beta}$
hit the factor~$\nabla^{(m)}R_{ijkl}$ and the derivatives
$\nabla_{a_1},\dots ,\nabla_{a_\beta}$ hit the factor~$\nabla^{(2)}\psi_1$. We observe that otherwise, we will either
have an internal contraction in the factor $\nabla^{(p)}\psi_1$,
or a factor $\nabla^{(p)}\psi_h$, $p,h\ge 3$, or a factor
$\nabla^{(t)}\psi_2, t\ge 1$, or a~factor~$\nabla^{(y)}\psi_1$
with an internal contraction.

 We denote the complete contraction that we thus obtain by
$\hat{C}^{\alpha,\beta}_{g}$. Moreover, we clearly observe that
unless $\beta=\kappa_1$, $\hat{C}^{\alpha,\beta}_{g}$ does not
belong to ${\rm silly}_{+}[I^s_{g}]$. To summarize, we have seen that
${\rm silly}_{+}[I^s_{g}]$ is the sublinear combination in
$\sum_{\gamma\in \Gamma_{\kappa_1}} a_\gamma
C^\gamma_{g}(\psi_1,\dots, \psi_s)$ that consists of complete
contractions in the form~(\ref{inje}).

We will now demonstrate our second claim. We employ the
decomposition:
\begin{gather*}
\nabla^{il}W_{ijkl}=\frac{n-3}{n-2}
\nabla^{il}R_{ijkl}+ ({\rm const})_1\nabla^2_{jk}R+({\rm const})_2\Delta R
g_{jk}.
\end{gather*}

We deduce that:
\begin{gather*}
{\rm silly}_{*}[I^s_{g}(\psi_1,\dots ,\psi_s)]=
\sum_{\gamma\in \Gamma_{\kappa_1}} a_\gamma \frac{n-3}{n-2}
C^{\gamma,\iota}_{g}(\psi_1,\dots, \psi_s).
\end{gather*}
This shows our second claim.

\subsection*{Acknowledgements}

This work has absorbed
the best part of the author's energy over many years.
This research was partially conducted during
the period the author served as a Clay Research Fellow,
an MSRI postdoctoral fellow,
a Clay Liftof\/f fellow and a Procter Fellow.
The author is immensely indebted to Charles
Fef\/ferman for devoting twelve long months to the meticulous
proof-reading of the present paper. He also wishes to express his
gratitude to the Mathematics Department of Princeton University
for its support during his work on this project.

\newpage

\pdfbookmark[1]{References}{ref}
\LastPageEnding

\end{document}